\documentclass[a4paper, 12pt]{article}
\usepackage{amsmath, amssymb, amsthm, mathrsfs}
\usepackage[colorlinks=true]{hyperref}

\newtheorem{thm}{Theorem}[section]
\newtheorem{lem}[thm]{Lemma}
\newtheorem{cor}[thm]{Corollary}
\newtheorem{prop}[thm]{Proposition}

\theoremstyle{definition}
\newtheorem{rem}{Remark}[section]

\numberwithin{equation}{section}

\allowdisplaybreaks[3]
\title{Well-posedness of the Cauchy Problem \\ for the Kinetic DNLS on $\mathbf{T}$}
\author{Nobu KISHIMOTO$^*$ and Yoshio TSUTSUMI$^{**}$}
\date{$^*$ RIMS, Kyoto University, Kyoto 606-8502, JAPAN \\
$^{**}$ Department of Mathematics, Kyoto University, \\ Kyoto 606-8502, JAPAN}

\begin{document}

\maketitle

\begin{abstract}
We consider the Cauchy problem for the \emph{kinetic derivative nonlinear Schr\"odinger equation} on the torus:
\[ \partial_t u - i \partial_x^2 u = \alpha \partial_x \bigl ( |u|^2 u \bigr ) + \beta \partial_x \bigl [ H \bigl ( |u|^2 \bigr ) u \bigr ] , \quad (t, x) \in [0,T] \times \mathbf{T},\]
where the constants $\alpha,\beta$ are such that $\alpha \in \mathbf{R}$ and $\beta <0$, and $H$ denotes the Hilbert transform. 
This equation has dissipative nature, and the energy method is applicable to prove local well-posedness of the Cauchy problem in Sobolev spaces $H^s$ for $s>3/2$.
However, the gauge transform technique, which is useful for dealing with the derivative loss in the nonlinearity when $\beta =0$, cannot be directly adapted due to the presence of the Hilbert transform.
In particular, there has been no result on local well-posedness in low regularity spaces or global solvability of the Cauchy problem.

In this article, we shall prove local and global well-posedness of the Cauchy problem for small initial data in $H^s(\mathbf{T})$, $s>1/2$.
To this end, we make use of the parabolic-type smoothing effect arising from the resonant part of the nonlocal nonlinear term $\beta \partial_x [H(|u|^2)u]$, in addition to the usual dispersive-type smoothing effect for nonlinear Schr\"odinger equations with cubic nonlinearities.
As by-products of the proof, we also obtain smoothing effect and backward-in-time ill-posedness results.

\end{abstract}


\section{Introduction}

In the present paper, we consider the local and global well-posedness of the following Cauchy problem for the \emph{kinetic derivative nonlinear Schr\"odinger equation} (KDNLS) on the one-dimensional torus $\mathbf{T}=\mathbf{R}/2\pi \mathbf{Z}$.
\begin{alignat}{2}
   \partial_t u - i \partial_x^2 u &= \alpha \partial_x \bigl ( |u|^2 u \bigr ) + \beta \partial_x \bigl [ H \bigl ( |u|^2 \bigr ) u \bigr ] , &\qquad (t, x) &\in \mathbf{R} \times \mathbf{T}, \label{kdnls} \\
   u(0, x) &= u_0(x), & x &\in \mathbf{T}, \label{ic}
\end{alignat}
where $\alpha$, $\beta$ are real constants, $\beta \neq 0$, and $H$ is the Hilbert transform.
This equation models propagation of weakly nonlinear and dispersive Alfv\'en waves in a plasma (see, e.g., \cite{DP} and \cite{MW}).
The nonlocal nonlinear term $\partial _x [H(|u|^2)u]$ takes the resonant interaction between the wave modulation and the ions into account, while it is ignored in the following \emph{derivative nonlinear Schr\"odinger equation} (DNLS).
\begin{align}
   &\partial_t u - i \partial_x^2 u = \alpha \partial_x \bigl ( |u|^2 u \bigr ) , \qquad (t, x) \in \mathbf{R} \times \mathbf{T}. \label{dnls}
\end{align} 
The word ``kinetic'' implies that the collective motion of ions in a plasma is modeled by the Vlasov equation and not by the fluid equation.

The KDNLS equation posed on $\mathbf{R}$ is invariant under the scaling transformation
\begin{equation*}
u_\lambda (t,x) = \lambda ^{\frac12} u (\lambda ^2t, \lambda x),\qquad \lambda >0,
\end{equation*}
thus the critical Sobolev space is $L^2$.
A smooth solution $u$ of \eqref{kdnls} satisfies the following $L^2$ law:
\begin{equation}\label{L2conservation}
\frac{d}{dt}\| u(t) \|_{L^2}^2 = \beta \big\| D_x^{\frac12}\bigl ( |u(t)|^2\bigr ) \big\| _{L^2}^2,
\end{equation}
where $D_x=(-\partial_x^2)^{1/2}$.
This in particular shows that the $L^2$ norm is conserved for the DNLS equation \eqref{dnls}.
In fact, \eqref{dnls} is completely integrable and it has infinitely many conserved quantities, among which the following energy is included:
\begin{equation*}
E(u)=\int \left( |\partial_xu|^2-\frac{3\alpha}{2}|u|^2\mathrm{Im}(\bar{u}\partial_xu)+\frac{\alpha^2}{2}|u|^6\right) \,dx.
\end{equation*}
For the KDNLS equation \eqref{kdnls}, conservation of the $L^2$ norm or the energy itself is no longer true, but \eqref{L2conservation} tells us that the $L^2$ norm is non-increasing when $\beta < 0$, which means that the equation has a dissipative nature in this case.
Indeed, it was shown by the authors \cite{KT20} that a functional on $H^1$ analogous to the energy $E(u)$ for \eqref{dnls} can be used to derive the $H^1$ a priori bound for small solutions of \eqref{kdnls} in the dissipative case $\beta <0$.

We recall some of the well-posedness results for \eqref{dnls} and \eqref{kdnls} in Sobolev spaces $H^s$.
The nonlinear structure of the DNLS equation \eqref{dnls} allows us to apply the energy method, which yields local well-posedness for $s>3/2$ equally in the cases of $\mathbf{R}$ and $\mathbf{T}$.
It seems that this is also the case for the dissipative KDNLS equation \eqref{kdnls}, $\beta<0$.

The main obstacle to go beyond $s=3/2$ is loss of a derivative in the nonlinearity, and we need to use the dispersive smoothing effect of the equations.
Note, however, that the Fourier restriction norm method (see Bourgain \cite{Bour}) is not sufficient by itself to recover the first derivative.
For DNLS, a particularly useful tool is the gauge transform
\begin{equation}\label{ga1}
\begin{split}
&v(t) = e^{-i \frac{\alpha}{2}\varphi [u(t)]} u(t),\\
&\varphi [f](x)=\begin{cases}
\displaystyle\int _{-\infty}^x|f(y)|^2\,dy &(x\in \mathbf{R}),\\
\dfrac{1}{2\pi}\displaystyle\int_0^{2\pi} \int_\theta^x \left( |f(y)|^2 - \frac 1 {2\pi} \| f\|_{L^2}^2 \right) \, dy \, d\theta &(x\in \mathbf{T}),
\end{cases} 
\end{split}
\end{equation}
which converts \eqref{dnls} to a tamer equation with cubic nonlinearity $v^2 \partial_x \bar v$ instead of $\partial_x (|v|^2v)$.
Indeed, the cubic term $v^2 \partial_x \bar v$ can be treated with the Fourier restriction norm, while this is not the case for $|v|^2 \partial_x v$ or $\partial_x (|v|^2v)$.
Using the gauge transform \eqref{ga1}, local well-posedness for \eqref{dnls} in $H^s$, $s\geq 1/2$ was shown by Takaoka \cite{Tak1} for the case of $\mathbf{R}$ and by Herr \cite{H} for the case of $\mathbf{T}$. 
However, this gauge transform does not seem to work for the KDNLS equation \eqref{kdnls} as well as for \eqref{dnls}, because the presence of the Hilbert transform $H$ prevents us from removing the term $H(\bar u \partial_x u) u$.
This term causes the same problem as $|u|^2 \partial_x u$.

In the case of $\mathbf{R}$, another way to recover the first derivative is making use of the local smoothing estimate and its inhomogeneous version:
\begin{gather*}
\big\| D_x^{\frac12}e^{it\partial_x^2}u_0\big\| _{L^\infty _xL^2_t}\leq C\| u_0\|_{L^2},\\
\Big\| D_x\int _0^te^{i(t-\tau )\partial_x^2}F(\tau )\,d\tau \Big\|_{L^\infty_xL^2_t}\leq C\| F\|_{L^1_xL^2_t},
\end{gather*}
together with suitable estimates on maximal functions.
In fact, these estimates enable us to show local well-posedness of \eqref{dnls}--\eqref{ic} for $u_0\in H^s(\mathbf{R})$, $s>1/2$ with small $L^2$ norm by a contraction argument without using the gauge transform (see Kenig \mbox{et al.} \cite{KPV}, Molinet and Ribaud \cite{MR04a}; see also Guo \cite{G} for the extension to $s=1/2$).
Again, this method can not be directly adapted to the KDNLS equation, due to the failure of boundedness of the Hilbert transform on $L^1$.%
\footnote{Peres de Moura and Pastor \cite{MP} tried to overcome this technical difficulty by introducing a different gauge transform.
They claimed small data local well-posedness for \eqref{kdnls} in $H^s(\mathbf{R})$, $s>1/2$, both for $\beta <0$ and $\beta >0$.
However, there seems a gap in their estimate of the anti-derivative term which appears after applying the gauge transform.}
We also note that the argument with the local smoothing estimate is available only in the case of $\mathbf{R}$ and there is no counterpart in the periodic setting.
To the best of our knowledge, there is no result on the well-posedness of the Cauchy problem for the KDNLS equation \eqref{kdnls}--\eqref{ic} in low regularities, in the both cases of $\mathbf{R}$ and $\mathbf{T}$. 

In the present paper, we focus on the periodic and $\beta <0$ case.
In this setting, the dissipative aspect of the equation \eqref{kdnls} appears not only as the decreasing $L^2$ norm \eqref{L2conservation} but more clearly as a \emph{smoothing effect of parabolic type}.
To see this, we first examine the resonant frequency part of the cubic nonlinear term $\partial_x [H(|u|^2) u]$.
We define the resonance function $R$ as follows.
\begin{align*}
   R(k_1, k_2, k_3) &= (\tau + k^2) - (\tau_1 + k_1^2) + (\tau_2 + k_2^2) - (\tau_3 + k_3^2) \\
   &= 2(k_1-k_2)(-k_2+k_3),
\end{align*}
where
\[
   k_1,k_2,k_3\in \mathbf{Z},\qquad \tau = \tau_1 - \tau_2 + \tau_3, \quad k = k_1 - k_2 + k_3.
\]
A triplet of frequencies $(k_1, k_2, k_3)$ is said to be resonant if $R(k_1, k_2, k_3) = 0$.
Basically, we can gain regularity of at most $|R|^{1/2}$ for the non-resonant frequency part by means of the Fourier restriction norm, whereas such a dispersive smoothing effect is not available for the resonant part.
We put $N = \partial_x [H(|u|^2) u]$ and take the Fourier transform of $N$ in the spatial variable to have
\[
\hat N(k) = ik \sum_{k_1-k_2+k_3 = k} \Bigl ( -i \,\mathrm{sgn}(k_1-k_2) \hat u(k_1) \bar {\hat u}(k_2) \Bigr ) \hat u(k_3).
\]
Here and hereafter, we use the convention that $\mathrm{sgn}(0)=0$.
Then, the resonant part of $N$ is extracted as
\begin{align*}
&k \sum_{\begin{subarray}{c} k_1-k_2+k_3 = k \\ -k_2+k_3=0 \end{subarray}}\mathrm{sgn}(k_1-k_2) \hat u(k_1) \bar {\hat u}(k_2) \hat u(k_3)\\
&\quad = k \sum _{k_3\in \mathbf{Z}} \mathrm{sgn}(k-k_3)|\hat{u}(k_3)|^2\hat{u}(k)\\
&\quad = \Big( \sum _{k_3\in \mathbf{Z}}|\hat{u}(k_3)|^2\Big) |k| \hat{u}(k) \ + \ \sum _{k_3\in \mathbf{Z}} k \Big( \mathrm{sgn}(k-k_3)-\mathrm{sgn}(k)\Big) |\hat{u}(k_3)|^2\hat{u}(k).
\end{align*}
On one hand, the second sum in the last line has no derivative loss, since $\mathrm{sgn}(k-k_3)-\mathrm{sgn}(k)=0$ unless $|k_3|>|k|$.
Hence, this term can be treated in $H^s$ for $s\geq 1/2$ even without smoothing effect.
On the other hand, the remaining part
\[ 
\Big( \sum _{k_3\in \mathbf{Z}}|\hat{u}(k_3)|^2\Big) |k|\hat{u}(k) \ = \ \frac{\| u\|_{L^2}^2}{2\pi}|k|\hat{u}(k)
\]
is regarded as the resonant part including loss of the first derivative.
%
Therefore, we take the Fourier transform of \eqref{kdnls} (assuming $\alpha =0$, for simplicity) in the spatial variable to have
\begin{align}
   &\partial_t \hat u(t, k) + i k^2 \hat u(t,k) - \beta \frac{\| u(t) \|_{L^2}^2}{2\pi}|k|\hat{u}(t,k) = \hat{F}_1[u(t)](k)+\hat{F}_2[u(t)](k), \label{reseq} 
\end{align}
where the nonlinear parts $\hat{F}_1[u]$ and $\hat{F}_2[u]$ consist of the resonant terms with no derivative loss and of the non-resonant terms with loss of the first derivative, respectively.
If $\beta < 0$, then the third term on the left side of \eqref{reseq} has a parabolic-type dissipative effect on the system, unless $u(t)=0$.
Although $\hat{F}_2[u]$ does not contain resonant frequencies, the standard Fourier restriction norm method is not strong enough to recover the derivative loss in $\hat{F}_2[u]$ like the DNLS equation \eqref{dnls}.
Instead of applying a gauge transform, we shall take advantage of the dissipative effect of the third term on the left side of \eqref{reseq} besides the dispersive smoothing effect to compensate the derivative loss and establish the low-regularity well-posedness of the Cauchy problem \eqref{kdnls}--\eqref{ic}.
In \cite{Tsug}, Tsugawa first observed that some resonant part of the nonlinearity exhibits the parabolic type smoothing for a class of semilinear dispersive equations, showing the regularizing effect of smooth solution via an energy type argument.
He did not, however, use this smoothing property for the proof of existence of solution in low regularity, because the main object in \cite{Tsug} was to find a structural condition on the nonlinearity leading to parabolic nature of the equation.
In the present paper, we take advantage of the parabolic nature of the equation \eqref{kdnls} for the construction of solution in low regularity.
Indeed, we exploit the parabolic type smoothing effect emanating from the explicit resonant term on the left hand side of \eqref{reseq}, which yields the gain of additional $1/2$ derivative, to give a proof based on the contraction argument for the unique existence of solution.
Note that the contraction argument often does not work in the case the resonant frequency term can not be regarded as part of the inhomogeneous term, namely, in the case the influence of the resonant frequency term can not be represented by the Duhamel integral.

The observation on the role of the resonant frequency part in the nonlinearity has also been used to show the ill-posedness of the Cauchy problem (see, e.g., \cite{Tsug} for the fifth order semilinear dispersive equation and \cite{KT18,KT20b} for the third order dispersion NLS equation with Raman scattering term).
While the backward ill-posedness of the parabolic equation plays a crucial role in \cite{Tsug}, the ellipticity from the resonant part of nonlinearity exhibits the ill-posedness in both time directions in \cite{KT18,KT20b}.
In the present paper, we also obtain the backward ill-posedness from the parabolic type smoothing effect.

%
%

The equation \eqref{reseq} indicates that the effective ``linear'' part of the KDNLS equation \eqref{kdnls} is given by $\partial_t u - i\partial_x^2 u -\frac{\beta}{2\pi}\| u(t)\|_{L^2}^2D_xu$.
However, the $L^2$ norm is not conserved for \eqref{kdnls}, and the $t$-dependent coefficient is unfavorable in our argument.
We then replace the ``linear'' part with
\[
   \bigl (\partial_t  - i \partial_x^2  + \mu D_x \bigr ) u,\qquad \mu = - \frac{\beta}{2\pi} \| u_0 \|_{L^2}^2,
\]
at the cost of modifying the nonlinear part by $\frac{\beta}{2\pi}( \| u(t)\|_{L^2}^2-\| u_0\|_{L^2}^2)D_xu$, which is expected to remain small for a short time.\footnote{%
The idea of treating the difference of a non-conserved quantity at two different times as part of the nonlinearity has appeared in the work by Takaoka and the second author~\cite{TT04} on the modified KdV equation and other related works~\cite{NTT10,MPV19,MT18}.
In these results, such difference part needed to be treated in a less perturbative way using a subtle smoothing property of solution derived from the dispersion and the nonlinear structure of the equation.
In our problem, the estimate of the difference part falls within the scope of a contraction argument, thanks to the parabolicity of the linear part (see Lemma~\ref{lem:L2diff} and Proposition~\ref{prop:estN1} below).
}
To construct a solution of \eqref{kdnls}--\eqref{ic} forward in time, we assume $\beta <0$ and $u_0\neq 0$, so that $\mu >0$.
Then, taking the dissipative property of this dispersive-parabolic equation into account, we define a variant of the Fourier restriction norm space as follows:
\begin{gather*}
   \| f \|_{X_\mu ^{s,b}} = \Bigl ( \sum_{k \in \mathbf{Z}} \int_{\mathbf{R}} \langle k \rangle^{2s} \bigl \langle i(\tau + k^2) + \mu |k|\bigr \rangle^{2b} |(\mathcal{F}_{t,x} f)(\tau, k)|^2 \, d\tau \Bigr )^{1/2} ,\\
   X_\mu ^{s,b} = \bigl \{ f \in \mathscr{S}'(\mathbf{R}^2) \big | \ f(t, x + 2\pi) = f(t, x), \quad \| f \|_{X_\mu^{s,b}} < \infty \bigr \} ,
\end{gather*}
where $\langle a \rangle = (1 + |a|^2)^{1/2}$ for $a \in \mathbf{C}$.
This kind of spaces was introduced by Molinet and Ribaud \cite{MR} to study the KdV-Burgers equation
\[ 
\partial_t u + \partial_x^3 u - \partial_x^2 u + u \partial_x u = 0
\]
(see also Molinet and Vento \cite{MV1, MV2}).
The feature of our problem is that the dissipation originates in the nonlinearity of the equation.
In particular, the coefficient of the dissipative term depends on the $L^2$ norm of initial data, and so does the corresponding space $X^{s,b}_\mu$.

\begin{rem} \label{rem:linear1}
When we take $U_\mu (t) = e^{it \partial_x^2 - \mu |t| D_x}$ as the evolution operator associated with 
the linear part $(\partial_t  - i \partial_x^2  + \mu D_x)u$, we note by Lemma~\ref{lem:linear1} below that for $u_0 \in H^s$, $\psi(t) U_\mu (t) u_0$ belongs to $X_\mu^{s,1/2}$, where $\psi$ is a time cut-off function in $C_0^\infty(\mathbf{R})$.
So, the function space $X_\mu^{s,1/2}$ is suitable for $\psi (t) U_\mu (t) u_0$.
But we also note that 
\[
   \psi (t) U_\mu(t) u_0 \not \in X_\mu^{s, b}, \quad b > 1/2.
\]
The latter follows immediately from the Fourier transform of the Poisson kernel:
\begin{align*}
   \mathcal{F}_{t,x} \left [ e^{it\partial_x^2- \mu |t| D_x} \right ] = \frac {2\mu |k|} {(\tau + k^2)^2 + \mu^2 k^2}, \quad k \neq 0.
\end{align*}
This is in contrast to the case $\mu =0$, where $\psi (t)e^{it\partial_x^2}u_0\in X^{s,b}_0$ for any $b\in \mathbf{R}$.
We also present the so-called linear estimate in $X_\mu^{s,1/2}$ for the Duhamel term without sharp cut-off $\chi_{\mathbf{R}_+}(t)$ (see Lemma~\ref{lem:inhom} in Section~\ref{sec:linearestimates}).
Indeed, the linear estimate involving sharp cut-off in \cite{MR} is not very convenient, because the sharp cut-off multiplier is not bounded from $H^{1/2}(\mathbb{R})$ to $H^{1/2}(\mathbb{R})$ (see, e.g., Strichartz \cite[Theorem 3.1 on page 1054]{Str67}).
In order to make it clear what influence the modification of $e^{it \partial_x^2 -\mu t D_x}$ to $e^{it \partial_x^2 -\mu |t| D_x}$ has on the linear estimates, we state the linear estimates in a slightly general form and give a proof in Section~\ref{sec:appendix} of the appendix (see Lemmas \ref{lem:linear1-abs} and \ref{lem:inhom-abs}), which covers general symbols with dispersion and dissipation.
\end{rem}

Nonlinear hyperbolic and dispersive equations have several important properties in common. 
Namely, initial value problems for these equations can be solved towards positive time and negative time and most of them have conservation laws. 
It is natural and important what influence the addition of dissipative effects to these conservative systems has on regularity and asymptotic behavior of solutions. 
Indeed, the research on this issue for nonlinear hyperbolic systems has a long history (see, e.g., Liu and Zeng~\cite{LiuZ}). 
Nevertheless, little is known about nonlinear dissipative effects for nonlinear dispersive equations while there are many papers on nonlinear dispersive equations with linear dissipative terms such as the KdV-Burgers equation (see, e.g., \cite{MR,MV1,MV2}). 
This paper is one of attempts to make a close investigation into this issue.

The first result of the present paper is the following theorem on small-data local well-posedness of the Cauchy problem \eqref{kdnls}--\eqref{ic} in low-regularity Sobolev spaces.
\begin{thm} \label{local}
Let $\alpha \in \mathbf{R}$ and $\beta < 0$ be fixed, and let $s > 1/2$.
Then, the Cauchy problem \eqref{kdnls}--\eqref{ic} is locally well-posed for small data in $H^s(\mathbf{T})$.
More precisely, the following holds.
\begin{enumerate}
\item Existence: There exists $\eta =\eta (s)>0$ such that for any $u_0 \in H^s(\mathbf{T})$ with $\| u_0 \|_{H^s}\leq \eta$ the Cauchy problem has a solution $u \in C([0, 1]; H^s(\mathbf{T}))\cap L^2((0,1);H^{s+\frac12}(\mathbf{T}))$ on the time interval $[0,1]$.
\item Uniqueness: If $u_0\neq 0$, $\| u_0\|_{H^s}\leq \eta$, and $\tilde{u}\in Z^{s}_0(T)\cap L^2((0,T);H^{s+\frac12}(\mathbf{T}))$ is a solution to the Cauchy problem on $[0,T]$ for some $T\in (0,1]$, then $\tilde{u}$ coincides with the solution $u$ constructed in (i).
Here, the Banach space $Z^{s}_0(T)$ will be defined in Section~\ref{sec:fs}.
\item Continuity of the solution map: The solution map $u_0\mapsto u$ is continuous as a map from $\{ u_0\in H^s(\mathbf{T})\,|\,\| u_0\|_{H^s}\leq \eta \}$ into $C([0, 1]; H^s(\mathbf{T}))$.
\item Persistence of regularity: Let $u_0\in H^{s}(\mathbf{T})$ satisfy $\| u_0\| _{H^s}\leq \eta$, and assume that $u_0\in H^{\tilde{s}}(\mathbf{T})$ for some $\tilde{s}>s$.
Then, the solution $u$ constructed in (i) belongs to $C([0,1];H^{\tilde{s}}(\mathbf{T}))$.
\item Smoothing effect: The solution $u$ constructed in (i) is smooth for $t>0$; i.e., $u\in C^\infty ((0,1]\times \mathbf{T})$.
\end{enumerate}
\end{thm}

We prove the above theorem based on a contraction argument in a variant of the space $X^{s,b}_\mu$ (the space $Z^s_\mu (1)\simeq Z^s_0(1)\cap L^2((0,1);H^{s+\frac12}(\mathbf{T}))$ to be defined in Section~\ref{sec:fs}).
Precisely, we first apply the spatial translation 
\[ \mathcal{T}_\alpha \ :\ u(t,x)\quad \mapsto \quad u\Big( t,x-\frac{\alpha}{\pi}\| u(0)\|_{L^2}^2t\Big) ,\]
which converts the equation \eqref{kdnls} into the following ``renormalized'' equation:
\begin{equation}\label{r-kdnls}
\begin{split}
\partial_t u - i \partial_x^2 u = \alpha \partial_x \Big[ \Big( |u|^2 -\frac{1}{\pi} \| u(0)\|_{L^2}^2\Big) u\Big] + \beta \partial_x \bigl [ H \bigl ( |u|^2 \bigr ) u \bigr ] ,& \\
(t, x) \in [0,1]\times \mathbf{T}.&
\end{split}
\end{equation}
Then, we solve the Cauchy problem \eqref{r-kdnls}--\eqref{ic} by a contraction argument, and apply the inverse transform $\mathcal{T}_\alpha^{-1}=\mathcal{T}_{-\alpha}$ to obtain the solution of the Cauchy problem for the original equation.
Note that $\mathcal{T}_\alpha$ is a nonlinear map, but it is bicontinuous on $C([0,T];H^s)$ for any $s\geq 0$ and $T>0$, which is easily verified by the dominated convergence theorem.
We will see that $\mathcal{T}_\alpha$ also maps $Z^{s}_0(T)\cap L^2((0,T);H^{s+\frac12}(\mathbf{T}))$ onto itself (see Section~\ref{sec:fs}).

Our contraction argument for the renormalized equation \eqref{r-kdnls} differs in several ways from a standard Fourier restriction norm method. 
First, since the $X^{s,b}_\mu$ norm depends on the $L^2$ norm of initial data, proof of continuous dependence requires an additional care if we consider such a perturbation of initial data that changes the $L^2$ norm.
The precise statement of the result on the continuous dependence (including Theorem~\ref{local} (iii)) is given in the next proposition.  
\begin{prop}\label{prop:cont}
Let $\alpha \in \mathbf{R}$, $\beta <0$, and $s>1/2$.
Denote by $S$ and $\overline{S}=\mathcal{T}_\alpha ^{-1}\circ S$ the solution map of the Cauchy problem for the renormalized equation \eqref{r-kdnls}-\eqref{ic} (given in the proof of Theorem~\ref{local}, see Proposition~\ref{prop:contraction} below) and that for the original KDNLS equation \eqref{kdnls}--\eqref{ic} (whose existence is claimed in Theorem~\ref{local} (i)), respectively.
(Note that $S$ and $\overline{S}$ are defined on $B_{H^s}(\eta ):=\{ u_0\in H^s(\mathbf{T})\,|\,\| u_0\|_{H^s}\leq \eta \}$ and give the solution in $C([0, 1]; H^s(\mathbf{T}))\cap L^2((0,1);H^{s+\frac12}(\mathbf{T}))$.)
Then, the following holds.
\begin{enumerate}
\item For any $0<\varepsilon \leq \eta$, the map $S$ is Lipschitz continuous on $B_{H^s}(\eta )\cap \{ \| u_0\|_{L^2}\geq \varepsilon \}$ with respect to the topologies of $H^s$ and $C_tH^s_x\cap L^2_tH^{s+\frac12}_x$.
In particular, the same is true for $\overline{S}$ when $\alpha =0$, and even when $\alpha \neq 0$, $\overline{S}$ is Lipschitz continuous in the above sense on $B_{H^s}(\eta )\cap \{ \| u_0\|_{L^2}=c\}$ for each $0<c\leq \eta$.
\item The map $S$ is continuous at the origin with respect to the topologies $H^s$ and $C_tH^s_x$.
In particular, $S$ is continuous in this sense on the whole domain $B_{H^s}(\eta )$, and the same is true for $\overline{S}$ (by the continuity of $\mathcal{T}_{\alpha}^{-1}$).
\item Both $S$ and $\overline{S}$ are not uniformly continuous with respect to the topologies $H^s$ and $C_tH^s_x$ in any $H^s$-neighborhood of the origin.
\end{enumerate}
\end{prop}

The second point we need to take care of is that the parabolic smoothing effect ceases to work as soon as the solution becomes zero.
Indeed, since we do not know uniqueness backward in time, it is not straightforward to exclude the possibility of a nontrivial solution vanishing in finite time.
Then, the statements of uniqueness, persistence of regularity, and smoothing effect in the theorem will be first shown up to the time when the solution may become zero.
Once we obtain smoothness of the solution up to such time, the following \emph{a priori} estimates for smooth solutions are applicable, from which we know \emph{a posteriori} that the nontrivial solution will never become zero in finite time.

\begin{prop}[{\cite[Theorem~1]{KT20}}] \label{prop:apriori}
Let $\alpha ,\beta \in \mathbf{R}$ and $u:[0,T]\times \mathbf{T}\to \mathbf{C}$ be a smooth solution to the KDNLS equation \eqref{kdnls}.
\begin{enumerate}
\item $L^2$ equality: It holds that
\[ \| u(t)\| _{L^2}^2-\beta \int _0^t\big\| D_x^{\frac{1}{2}}(|u(\tau )|^2)\big\| _{L^2}^2\,d\tau \, =\, \| u(0)\| _{L^2}^2,\quad t\in [0,T]. \]
\item $H^1$ a priori bound and $L^2$ lower bound: 
We further assume that $\beta <0$ and that $\| u(0)\|_{L^2}$ is sufficiently small.
Then, there exists $C>0$ such that the following holds:
\begin{align*}
\| u(t)\|_{H^1}&\leq C\| u(0)\|_{H^1}e^{C\| u(0)\|_{L^2}^2},\\
\| u(t)\|_{L^2}&\geq \| u(0)\| _{L^2}\exp \Big[ -C\| u(0)\|_{H^1}e^{C\| u(0)\|_{L^2}^2}t^{\frac{1}{2}}\Big] ,\quad t\in [0,T].
\end{align*}
\end{enumerate}
\end{prop}

This proposition was shown by a simple energy-type argument. 
These \mbox{a priori} bounds and the energy method are employed to prove further properties of the Cauchy problem \eqref{kdnls}--\eqref{ic}, which we collect in the next theorem:
\begin{thm}\label{global}
Let $\alpha ,\beta \in \mathbf{R}$.
The following holds for the Cauchy problem \eqref{kdnls}--\eqref{ic}.
\begin{enumerate}
\item $L^2$ decay, uniqueness for $u_0=0$: 
Let $s>1/2$ and $u\in C([0,T];H^s(\mathbf{T}))$ be a solution to \eqref{kdnls}.
Then, $u$ satisfies the same $L^2$ equality as in Proposition~\ref{prop:apriori} (i).

In particular, if $\beta \leq 0$, any solution of the Cauchy problem in the class $C([0,T];H^s(\mathbf{T}))$ does not increase its $L^2$ norm.
Consequently, there is no nontrivial solution in $C([0,T];H^s(\mathbf{T}))$ starting from $u_0=0$.
\item Unconditional uniqueness for smooth solutions:
Assume $\beta \leq 0$, and let $s>3/2$.
Then, the solution to the Cauchy problem with $u_0\in H^s(\mathbf{T})$ is unique in the class $C([0,T];H^s(\mathbf{T}))$.
\item Backward ill-posedness (non-existence) for small smooth data:
Let $s>3/2$.
Assume that $u_0\in H^s(\mathbf{T})\setminus C^\infty (\mathbf{T})$ satisfies $\| u_0\| _{H^{s'}}< \eta (s')$ for some $s'\in (1/2,s]$, where $\eta (s')$ is the threshold size for local existence given in Theorem~\ref{local} (i). 
If $\beta <0$ (resp.~$\beta >0$), then for any $T>0$ there exists no backward (resp.~forward) solution to the Cauchy problem in the class $C([-T,0];H^s(\mathbf{T}))$ (resp.~$C([0,T];H^s(\mathbf{T}))$).
\item Global existence for small low-regularity data:
Assume $\beta <0$, and let $s>1/2$.
Then, there exists $\tilde{\eta}=\tilde{\eta}(s)>0$ (smaller than $\eta (s)$ given in Theorem~\ref{local} (i)) such that for any $u_0\in H^s(\mathbf{T})$ with $\| u_0\|_{H^s}\leq \tilde{\eta}$, the solution $u$ to the Cauchy problem on $[0,1]$ constructed in Theorem~\ref{local} can be extended to a global solution in $C([0,\infty );H^s(\mathbf{T}))\cap C^\infty ((0,\infty )\times \mathbf{T})$ satisfying
\[ \sup _{t\in [0,\infty )}\| u(t)\|_{H^{s_0}}\leq \eta (s_0),\quad
s_0:=\begin{cases}
\frac{1}{3-2s}\in (\tfrac{1}{2},s) &(\frac12 <s<1),\\
\ 1 &(s\geq 1).
\end{cases}
\]
\end{enumerate}
\end{thm}
The first two statements (i), (ii) are consequences of the $L^2$ energy argument for a solution and the difference of two solutions, respectively. 
For (ii), the assumption $\beta \leq 0$ is essential, because it ensures that the term including loss of regularity is non-positive and hence can be dropped in the energy estimate.
(iii) follows immediately from (ii) and the smoothing effect (Theorem~\ref{local} (v)).
As for (iv), it is easily shown when $s\geq 1$ by Theorem~\ref{local} and the a priori $H^1$ bound in Proposition~\ref{prop:apriori} (ii).
Even for $s\in (1/2,1)$, we can obtain a small a priori bound in $H^{s_0}$ (not in $H^s$, though) for positive time by using the smoothing effect (Theorem~\ref{local} (v)) and interpolating the $L^2$ and the $H^1$ bounds given in Proposition~\ref{prop:apriori}.

\begin{rem}\label{rem:lwp}
Here, we collect some problems left open in the present paper.
\begin{enumerate}
\item An important open problem is to remove the assumptions on smallness of initial data.
For the local result, it is expected that the smallness condition can be completely removed.
For global existence, it is natural to expect that the smallness assumption in the $H^s$ norm of initial data can be relaxed to that in the $L^2$ norm, since the equation \eqref{kdnls} is $L^2$-critical.
We remark that for the integrable DNLS equation (the case $\beta =0$) posed on $\mathbf{R}$, global existence has recently been shown in $H^{1/6}(\mathbf{R})$ without any size restriction (see Bahouri and Perelman \cite{BP}, Harrop-Griffiths \mbox{et al.} \cite{HKV}).  
\item Our argument to show well-posedness relies essentially on the dissipation arising in the resonant nonlinear part, which is available when $\beta <0$ and has strength $\sim |\beta |\| u(t)\|_{L^2}^2D_x$.
Indeed, the constants are not uniform as $|\beta |\to 0$ in various estimates we will show below.
Although the DNLS equation ($\beta =0$) is also possible to solve (but by different means, such as gauge transform or complete integrability), we do not know whether the solution of KDNLS with $\beta <0$ converges to that of DNLS as $\beta \to 0$.
This ``singular limit'' problem should be of both mathematical and physical interests.
Note that in Theorem~\ref{local} the threshold size $\eta (s)$ does in fact depend on $|\beta |$ and approaches to zero as $|\beta |\to 0$. 
\item It is also crucial to consider the problem in the periodic setting, and a similar argument is not likely to work in the case of $\mathbf{R}$.
The short-time Fourier restriction norm method (see Ionescu \mbox{et al.} \cite{IKT}) may be useful in the case of $\mathbf{R}$.
\end{enumerate}
\end{rem}

The rest of the paper is organized as follows.
In Section~\ref{sec:prel}, we introduce various function spaces and give the relevant linear estimates.
The required estimates for nonlinear terms are proved in Section~\ref{sec:nonlinear}.
Using these estimates, we prove Theorem~\ref{local} and Proposition~\ref{prop:cont} in Section~\ref{sec:local} and Section~\ref{section:continuity}, respectively.
In Section~\ref{sec:proof2}, Theorem~\ref{global} is proved.
Finally, in Appendix~\ref{sec:appendix} we give a proof to some linear estimates stated in Section~\ref{sec:prel}.


\section{Preliminaries}
\label{sec:prel}

In this section, we first introduce an integral equation which is equivalent to the Cauchy problem \eqref{kdnls}--\eqref{ic}, and then define the function spaces and present several lemmas which will be used in the proof of the main results.
We use the notation $\mathcal{F}_x,\mathcal{F}_t,\mathcal{F}_{t,x}$ for the Fourier transformation with respect to the indicated variables, and $\hat{f}$ denotes one of $\mathcal{F}_x[f(x)],\mathcal{F}_x[f(t,x)]$ and $\mathcal{F}_t[f(t)]$, while $\tilde{f}$ is used only for $\mathcal{F}_{t,x}[f(t,x)]$.
We use the following definition of $\mathcal{F}$:
\begin{alignat*}{2}
\mathcal{F}_x[f](k)&:=\frac{1}{2\pi}\int _0^{2\pi}e^{-ikx}f(x)\,dx,&\qquad k&\in \mathbf{Z}; \\
\mathcal{F}_t[g](\tau )&:=\frac{1}{2\pi}\int _{\mathbf{R}}e^{-i\tau t}g(t)\,dt,& \tau &\in \mathbf{R}.
\end{alignat*}

\subsection{Formulation of the problem}

We decompose the nonlinear term $N=\partial_x [H(|u|^2)u]$ as follows:
\begin{align*}
\hat{N}(k)
&=ik\bigg[ \sum _{\begin{smallmatrix} k_1-k_2+k_3=k\\ k_2=k_3\end{smallmatrix}} +\sum _{\begin{smallmatrix} k_1-k_2+k_3=k\\ k_2\neq k_3\end{smallmatrix}}\bigg] (-i)\,\mathrm{sgn}(k_1-k_2)\hat{u}(k_1)\bar{\hat{u}}(k_2)\hat{u}(k_3) \\
&=k\Big( \sum _{k_3\in \mathbf{Z}}\mathrm{sgn}(k-k_3)|\hat{u}(k_3)|^2\Big) \hat{u}(k)\\
&\qquad +k\sum _{\begin{smallmatrix} k_1-k_2+k_3=k\\ k_2\neq k_1,k_3\end{smallmatrix}}\mathrm{sgn}(k_1-k_2)\hat{u}(k_1)\bar{\hat{u}}(k_2)\hat{u}(k_3) \\
&=\Big( \sum _{k_3\in \mathbf{Z}}|\hat{u}_0(k_3)|^2\Big) |k|\hat{u}(k) +\Big( \sum _{k_3\in \mathbf{Z}}|\hat{u}(k_3)|^2-\sum _{k_3\in \mathbf{Z}}|\hat{u}_0(k_3)|^2\Big) |k|\hat{u}(k)\notag \\
&\quad +\Big( \sum _{k_3\in \mathbf{Z}}k\big[ \mathrm{sgn}(k-k_3)-\mathrm{sgn}(k)\big] |\hat{u}(k_3)|^2\Big) \hat{u}(k) \\
&\quad +k\sum _{\begin{smallmatrix} k_1-k_2+k_3=k\\ k_2\neq k_1,k_3\end{smallmatrix}}\mathrm{sgn}(k_1-k_2)\hat{u}(k_1)\bar{\hat{u}}(k_2)\hat{u}(k_3). 
\end{align*}
Similarly, we have
\begin{align*}
&\mathcal{F}_x\big[ \partial_x (|u|^2u) \big] (k) \\
&\quad = ik\bigg[ \sum _{\begin{subarray}{c} k_1-k_2+k_3=k\\ k_1=k_2\neq k_3\end{subarray}} + \sum _{\begin{subarray}{c} k_1-k_2+k_3=k\\ k_1\neq k_2=k_3\end{subarray}} - \sum _{\begin{subarray}{c} k_1-k_2+k_3=k\\ k_1=k_2=k_3\end{subarray}} + \sum _{\begin{subarray}{c} k_1-k_2+k_3=k\\ k_2\neq k_1,k_3\end{subarray}} \bigg] \hat{u}(k_1)\bar{\hat{u}}(k_2)\hat{u}(k_3) \\
&\quad =2ik\Big( \sum _{k_1\in \mathbf{Z}}|\hat{u}_0(k_1)|^2\Big) \hat{u}(k) +2ik\Big( \sum _{k_1\in \mathbf{Z}}|\hat{u}(k_1)|^2-\sum _{k_1\in \mathbf{Z}}|\hat{u}_0(k_1)|^2\Big) \hat{u}(k)\\
&\qquad -ik|\hat{u}(k)|^2\hat{u}(k) + ik \sum _{\begin{subarray}{c} k_1-k_2+k_3=k\\ k_2\neq k_1,k_3\end{subarray}} \hat{u}(k_1)\bar{\hat{u}}(k_2)\hat{u}(k_3) .
\end{align*}

Based on this observation, we define the following operators:
\begin{equation}\label{eq:nonlinearity}
\begin{split}
N_1[u,v;u_3]&:=\big( \| u\| _{L^2}^2-\| v\|_{L^2}^2\big) \Big( \frac{\alpha}{\pi}\partial_x u_3 + \frac{\beta}{2\pi}D_xu_3\Big) ,\\
N_2[u_1,u_2,u_3]&:=\mathcal{F}^{-1}_k\bigg[ -i\alpha k \hat{u}_1(k)\bar{\hat{u}}_2(k) \hat{u}_3(k)\\
&\qquad\quad +\beta \Big\{ \sum _{k'\in \mathbf{Z}}k\big[ \mathrm{sgn}(k-k')-\mathrm{sgn}(k)\big] \hat{u}_1(k')\bar{\hat{u}}_2(k')\Big\} \hat{u}_3(k)\bigg] ,\\
N_3[u_1,u_2,u_3]&:=\mathcal{F}^{-1}_k\bigg[ ik\!\!\sum _{\begin{smallmatrix} k_1-k_2+k_3=k\\ k_2\neq k_1,k_3\end{smallmatrix}}\!\!\big[ \alpha -i\beta \,\mathrm{sgn}(k_1-k_2)\big] \hat{u}_1(k_1)\bar{\hat{u}}_2(k_2)\hat{u}_3(k_3)\bigg] .
\end{split}
\end{equation}
The equation \eqref{kdnls} is rewritten as
\begin{align*}
\partial _tu=i\partial _x^2u&+\frac{\alpha}{\pi}\| u_0\| _{L^2}^2 \partial_x u +\frac{\beta}{2\pi}\| u_0\| _{L^2}^2D_xu\\
&+N_1[u,u_0;u]+N_2[u,u,u]+N_3[u,u,u] ,
\end{align*}
and the renormalized equation \eqref{r-kdnls} is rewritten as
\begin{equation}\label{renormalized}
\partial _tu=i\partial _x^2u+\frac{\beta}{2\pi}\| u_0\| _{L^2}^2D_xu +N_1[u,u_0;u]+N_2[u,u,u]+N_3[u,u,u].
\end{equation}
We will construct a solution to the Cauchy problem for the renormalized equation \eqref{renormalized}--\eqref{ic}, and then obtain a solution to the original Cauchy problem \eqref{kdnls}--\eqref{ic} by applying $\mathcal{T}_\alpha ^{-1}=\mathcal{T}_{-\alpha}$.

\subsection{Function spaces}
\label{sec:fs}

Let $\alpha \in \mathbf{R}$ and $\beta <0$.
We consider the integral equation associated with \eqref{renormalized}:
\begin{align}\label{eq:integral}
u(t)&=U_\mu (t)u_0+\int _0^t U_\mu (t-t')\big\{ N_1+N_2+N_3\big\} (t')\,dt' ,\quad t\in [0,T],
\end{align}
where 
\[ \mu :=\frac{|\beta |}{2\pi}\| u_0\|_{L^2}^2,\qquad U_\mu (t)\phi =\mathcal{F}^{-1}_k\big[ e^{-ik^2t-\mu|k||t|}\hat{\phi}(k)\big] ,\qquad t\in \mathbf{R}.\]
The modified evolution operator $U_\mu(t)$ is introduced by Molinet and Ribaud \cite{MR} in the setting of the KdV-Burgers equation (see also Molinet and Vento \cite{MV1, MV2}), which moderates the behavior of the solution in negative time.
Define the Banach spaces $X_\mu ^{s,b}$, $Y_\mu ^{s,b}$ corresponding to the linear propagator $U_\mu (t)$ for $\mu \geq 0$ by the restriction of spatially periodic and temporally tempered distributions $\{ f \in \mathscr{S}'(\mathbf{R}^2) \ | \ f(t, x + 2\pi) = f(t, x) \}$ to those of which the following norms are finite:
\begin{align*}
\| u\| _{X_\mu ^{s,b}}&:=\Big\| \langle k\rangle ^s\left\langle i(\tau +k^2)+\mu |k|\right\rangle ^b\tilde{u}(\tau ,k)\Big\| _{L^2_{\tau ,k}(\mathbf{R}\times \mathbf{Z})},\\
\| u\| _{Y_\mu ^{s,b}}&:=\Big\| \langle k\rangle ^s\big\| \left\langle i(\tau +k^2)+\mu |k|\right\rangle ^b\tilde{u}(\tau ,k)\big\| _{L^1_\tau (\mathbf{R})} \Big\| _{\ell ^2_k(\mathbf{Z})},
\end{align*}
where $\langle \cdot \rangle :=(1+|\cdot |^2)^{1/2}\sim 1+|\cdot |$.
When $b=0$, these norms are independent of $\mu$, so we often omit the subscript and write them as $X^{s,0}$, $Y^{s,0}$.
It is easy to see that $Y^{s,0}\hookrightarrow \{ u\in C(\mathbf{R};H^s)\,|\,u(t)\to 0~\text{in $H^s$}~(t\to \pm \infty )\}$.
Then, define the Banach spaces $Z^s_\mu$ and $N^s_\mu$ by
\[ Z^s_\mu :=X^{s,\frac{1}{2}}_\mu \cap Y^{s,0},\qquad N^{s}_\mu :=X^{s,-\frac{1}{2}}_\mu \cap Y^{s,-1}_\mu .\]
Finally, we define the restriction space $Z^s_\mu (T)$ for $T>0$ by
\begin{gather*}
Z^s_\mu (T):=\big\{ v|_{[0,T]\times \mathbf{T}}\,\big|\, v\in Z^s_\mu \big\} ,\\
\| u\| _{Z^s_\mu (T)}:=\inf \big\{ \| v\| _{Z^s_\mu}\,\big|\, v\in Z^s_\mu ,\, u=v|_{[0,T]\times \mathbf{T}}\big\} .
\end{gather*}
We see $Z^s_\mu (T)\hookrightarrow C([0,T];H^s)$, and $Z^s_\mu (T)\hookrightarrow L^2((0,T);H^{s+\frac{1}{2}})$ if $\mu >0$.
In fact, $Z^s_\mu (T)$ defines an equivalent Banach space to $Z^s_0(T)\cap L^2((0,T);H^{s+\frac12})$ for any $\mu >0$.
We also note that the spatial translation $T_\nu :u(t,x)\mapsto u(t,x-\nu t)$ ($\nu \in \mathbf{R}$ is a constant) maps $Z^s_\mu (T)$ into itself if $\mu >0$, because $\mathcal{F}_{t,x}[T_\nu u](\tau ,k)=\tilde{u}(\tau +\nu k,k)$ and $\langle i(\tau -\nu k+k^2)+\mu |k|\rangle \sim \langle i(\tau +k^2)+\mu |k|\rangle$.

The following property of the restriction space will be used to show uniqueness.
\begin{lem}\label{lem:vanishing}
Let $s\in \mathbf{R}$, $\mu \geq 0$ and $T>0$.
If $u\in Z^s_\mu (T)$ and $u|_{t=0}=0$, then 
\[ \lim _{T'\to 0}\big\| u| _{[0,T']\times \mathbf{T}}\big\| _{Z^s_\mu (T')}=0.\]
\end{lem} 

\noindent
{\bf Proof}.
We employ almost the same proof as for \cite[Claim (3.11)]{K09}.
We outline the argument below.
First, by approximation with smoother functions in $X^{s,1}_\mu (T)$ and the estimate of $U_\mu(t)\phi$ in $Z^s_\mu$ (see Lemma~\ref{lem:linear1} below), the lemma is reduced to showing the same claim for $X^{s,1}_\mu$ instead of $Z^s_\mu$.
But this is essentially the same as to show $\| f|_{[0,T']}\|_{H^1([0,T'])}\to 0$ for $f\in H^1(\mathbf{R})$ satisfying $f(0)=0$.
This is now elementary; for instance, if we set $g(t):=f(t)\chi _{[0,T']}(t)+f(2T'-t)\chi_{(T',2T']}(t)$, then $g \in H^1(\mathbf{R})$ and it is an extension of $f|_{[0,T']}$, $\| g\| _{H^1}^2=2(\| f\| _{L^2(0,T')}^2+\| f'\|_{L^2(0,T')}^2)\to 0$ as $T'\to 0$.
\hfill $\square$


\subsection{Linear estimates}\label{sec:linearestimates}

Here, we collect the linear estimates in the spaces defined above.
For the lemmas in this section, all the implicit constants do not depend on $s\in \mathbf{R}$, $\mu \geq 0$, and $0<T\leq 1$.
\begin{lem} \label{L1}
Let $\rho : \, \mathbf{Z} \mapsto [0,\infty )$.
Then,
\[ \sum_{k \in \mathbf{Z}} \Bigl ( \int_{\mathbf{R}} \frac { \big( 1+\rho (k)\big) ^{1/2}  |\tilde f(\tau, k)| } {1 + |\tau + k^2| + \rho (k) } \ d\tau \Bigr )^2  
~\lesssim~ \sum_{k \in \mathbf{Z}} \int_{\mathbf{R}}  |\tilde f(\tau, k)|^2\ d\tau. \]
In particular, we have
\[ \| u\| _{N^{s+\frac{1}{2}}_\mu} \lesssim \big( 1+\mu ^{-1}\big) ^{\frac{1}{2}}\| u\| _{X^{s,0}}\]
for any $s\in \mathbf{R}$ and $\mu >0$.
\end{lem}

\noindent
{\bf Proof}.
We follow the argument as used in the case of $|\lambda - k^3| > \frac 1 2 k^2$ for the proof of \cite[Lemma~7.42]{Bour}.
By the Cauchy-Schwarz inequality
, we have
\begin{align*}
&   \sum_{k \in \mathbf{Z}} \Bigl ( \int_{\mathbf{R}} \frac { \big( 1+\rho (k)\big) ^{1/2}  |\tilde f(\tau, k)| } {1 + |\tau + k^2| + \rho (k) } \ d\tau \Bigr )^2     \\
&   \leq \sum_{k \in \mathbf{Z}} \Bigl ( \int_{\mathbf{R}} \frac { 1+\rho (k) } { \bigl (1 + |\tau + k^2| + \rho (k) \bigr )^2 } \ d\tau \Bigr )  \Bigl ( \int_{\mathbf{R}} |\tilde f(\tau, k)|^2 \ d\tau \Bigr )  \\
&   \leq  \sum_{k \in \mathbf{Z}} \Bigl ( \Bigl [ - \frac {1+\rho} { 1 + \tau + k^2 + \rho } \Bigr ]_{-k^2}^\infty + \Bigl [ \frac {1+\rho} { 1 - \tau - k^2 + \rho } \Bigr ]_{-\infty}^{-k^2} \Bigl )  \Bigl ( \int_{\mathbf{R}}  |\tilde f(\tau, k)|^2 \ d\tau \Bigr )  \\
&   \leq 2 \sum_{k \in \mathbf{Z}} \int_{\mathbf{R}} | \tilde f(\tau, k)|^2 \ d\tau .
\end{align*}
The latter claim follows from the inequality $1+|k|\leq (1+\mu ^{-1})(1+\mu |k|)$.
This completes the proof of Lemma \ref{L1}.
\hfill $\square$


\bigskip
Let $\psi \in C^\infty _0(\mathbf{R})$ be a bump function satisfying $\psi =1$ on $[-1,1]$, $\mathrm{supp}\;\psi \subset [-2,2]$, and let $\psi _T(t)=\psi (t/T)$.
The next two lemmas were essentially proved by Molinet and Ribaud \cite{MR}, though they considered the estimate of the Duhamel integral with time cut-off $\chi_{\mathbf{R}_+}(t)$.
In Appendix~\ref{sec:appendix}, we restate these lemmas in a slightly general form as Lemmas~\ref{lem:linear1-abs} and \ref{lem:inhom-abs} and give a proof for completeness.

\begin{lem}[{cf.~\cite[Proposition~2.1]{MR}}]\label{lem:linear1}
For any $s\in \mathbf{R}$ and $\mu \geq 0$, we have 
\begin{equation*}
\| \psi (t)U_\mu (t)\phi \| _{Z^s_\mu}\lesssim \| \phi \| _{H^s}. 
\end{equation*}
In particular, there exists $C_0>0$ independent of $s,\mu,T$ such that
\[ \| U_\mu (t)\phi \| _{Z^s_\mu (T)}\leq C_0\| \phi \|_{H^s},\qquad 0<T\leq 1.\]
\end{lem}

\begin{lem}[{cf.~\cite[Proposition~2.3(a)]{MR}}]\label{lem:inhom}
For any $s\in \mathbf{R}$ and $\mu \geq 0$, we have
\begin{equation*}
\Big\| \psi (t)\int _0^t U_\mu (t-t')F(t')\,dt' \Big\| _{Z_\mu ^s}\lesssim \| F\| _{N_\mu ^{s}}. 
\end{equation*}
\end{lem}

\bigskip
Next, we recall the following estimates concerning the localization to the interval $[0,T]$.
These estimates are well-known for $\mu =0$, and the same argument works for $\mu >0$.
We omit the proof.
\begin{lem}\label{lem:linear2}
Let $s\in \mathbf{R}$ and $\mu \geq 0$.
For $0<T\leq 1$ and $0\le b\le \frac{1}{2}$, we have
\begin{gather*}
\| \psi _T(t)u\| _{X^{s,b}_\mu}\lesssim T^{\frac{1}{2}-b}\| u\| _{Z^{s}_\mu},\qquad 
\| \psi _T(t)u\| _{Y^{s,0}}\lesssim \| u\| _{Y^{s,0}}.
\end{gather*}
In particular, multiplication by $\psi _T(t)$ is a bounded operator on $Z^{s}_\mu$ uniformly in $s$, $\mu$, and $0<T\leq 1$.
\end{lem}


The next lemma follows immediately from Lemma~\ref{lem:inhom} and Lemma~\ref{L1}.
\begin{lem}\label{lem:nonlinear1}
For any $s\in \mathbf{R}$, $\mu >0$ and $0<T\leq 1$, we have
\[ \Big\| \psi (t)\int _0^tU_\mu (t-t')F(t')\,dt'\Big\| _{Z^{s+\frac12}_\mu}\lesssim (1+\mu ^{-1})^{\frac{1}{2}}\| F\| _{X^{s,0}}.\]
\end{lem}


\section{Nonlinear estimates}
\label{sec:nonlinear}

In this section, we estimate nonlinear terms $N_1$, $N_2$, and $N_3$ defined in \eqref{eq:nonlinearity}.
Hereafter, we fix $\alpha \in \mathbf{R}$, $\beta <0$ and assume that $u_0\neq 0$.
We will focus on the case of small initial data, so we also assume $0<\mu =(2\pi )^{-1}|\beta |\| u_0\|_{L^2}^2\lesssim 1$.
Again, the implicit constants in the estimates below are independent of $\mu$, $T$, and the regularity indices $s$ and $\sigma$, unless otherwise mentioned.

Let us begin with the next lemma:
\begin{lem}\label{lem:L2diff}
Let $T\in (0,1]$.

(i) For any  $s\in \mathbf{R}$ and $\alpha \in [0,1]$, we have
\begin{equation*}
\sup _{t\in [0,T]}\left| \| U_\mu (t)u_0\|_{H^s}^2-\| u_0\| _{H^s}^2\right| \lesssim T^{\alpha}\mu ^{\alpha}\| u_0\| _{H^{s+\frac{\alpha}{2}}}^2.
\end{equation*}

(ii) Let $u\in C([0,T];L^2)$ and assume $\| u(t)-U_\mu (t)u_0\| _{L^\infty ((0,T);L^2)}\leq \delta \| u_0\| _{L^2}$ for some $\delta \in (0,1]$.
Then, we have
\begin{equation*}
\sup _{t\in [0,T]}\left| \| u(t)\|_{L^2}^2-\| u_0\| _{L^2}^2\right| \lesssim \left( \delta + T^{\frac{1}{2}}\| u_0\| _{H^{\frac{1}{2}}}\right) \| u_0\| _{L^2}^2.
\end{equation*}
\end{lem} 

\medskip \noindent
\textit{Proof}.
(i) For $t\in [0,T]$ and $k\in \mathbf{Z}$,
\begin{align*}
\left| |\widehat{U_\mu (t)u_0}(k)|^2-|\hat{u}_0(k)|^2\right|
&=\left| \int _0^t\frac{d}{dt}\Big( e^{-2\mu |k|t}\Big) \Big|_{t=t'}\,dt' \right| |\hat{u}_0(k)|^2\\
&\leq \min \left\{ 2T\mu |k||\hat{u}_0(k)|^2,\, |\hat{u}_0(k)|^2\right\} \\
&\leq (2T\mu )^\alpha |k|^\alpha |\hat{u}_0(k)|^2,\qquad \alpha \in [0,1].
\end{align*}
Therefore, 
\begin{align*}
\left| \| U_\mu (t)u_0\|_{H^s}^2-\| u_0\| _{H^s}^2\right| 
&\leq \sum _{k\in \mathbf{Z}}\langle k\rangle ^{2s}\left| |\widehat{U_\mu (t)u_0}(k)|^2-|\hat{u}_0(k)|^2\right| \\
&\leq (2T\mu )^{\alpha}\| u_0\| _{H^{s+\frac{\alpha}{2}}}^2.
\end{align*}

(ii) In view of (i) with $s=0$, $\alpha =1/2$ and interpolation, we have
\[ \sup _{t\in [0,T]} \left| \| U_\mu (t)u_0\| _{L^2}^2-\| u_0\| _{L^2}^2\right| \lesssim T^{\frac{1}{2}}\| u_0\| _{L^2}\| u_0\| _{H^{\frac{1}{4}}}^2\lesssim T^{\frac{1}{2}}\| u_0\| _{L^2}^2\| u_0\| _{H^{\frac{1}{2}}}.\]
It suffices to observe that
\begin{align*}
&\sup _{t\in [0,T]} \left| \| u(t)\| _{L^2}^2-\| U_\mu (t)u_0\| _{L^2}^2\right| \\
&\quad \leq \sup _{t\in [0,T]}\Big( \| u(t)-U_\mu (t)u_0\| _{L^2}+2\| U_\mu (t)u_0\| _{L^2}\Big) \| u(t)-U_\mu (t)u_0\| _{L^2}\\
&\quad \leq \left( \delta \| u_0\| _{L^2}+2\| u_0\| _{L^2}\right) \delta \| u_0\| _{L^2}\lesssim \delta \| u_0\| _{L^2}^2
\end{align*}
under the assumption on $u$.\hfill $\Box$

\bigskip
From Lemma~\ref{lem:L2diff}, we deduce the following:
\begin{prop}\label{prop:estN1}
Let $\sigma \geq 0$, $T\in (0,1]$.

(i) Assume $\| u_j\| _{L^\infty ((0,T);L^2)}\lesssim \| u_0\| _{L^2}$ for $j=1,2$.
Then, we have
\begin{align*}
&\Big\| \int _0^tU_\mu (t-t')N_1[u_1(t'),u_2(t');v(t')]\,dt'\Big\| _{Z^\sigma _\mu (T)}\\
&\quad \lesssim \min \Big\{ \frac{\| v\| _{Z^\sigma _\mu (T)}}{\| u_0\|_{L^2}},\;T^{\frac{1}{2}}\| v\| _{L^\infty ((0,T);H^{\sigma +\frac{1}{2}})}\Big\} \| u_1-u_2\| _{L^\infty ((0,T);L^2)}.
\end{align*}

(ii) Assume $\| u-U_\mu (t)u_0\| _{L^\infty ((0,T);L^2)}\leq \delta \| u_0\| _{L^2}$ for some $\delta \in (0,1]$.
Then, 
\begin{align*}
&\Big\| \int _0^tU_\mu (t-t')N_1[u(t'),u_0;v(t')]\,dt'\Big\| _{Z^\sigma _\mu (T)}\\
&\quad \lesssim \Big( \delta +T^{\frac{1}{2}}\| u_0\| _{H^{\frac{1}{2}}}\Big) \min \Big\{ \| v\| _{Z^\sigma _\mu (T)},\;\| u_0\|_{L^2}T^{\frac{1}{2}}\| v\| _{L^\infty ((0,T);H^{\sigma +\frac{1}{2}})}\Big\} .
\end{align*}
%
\end{prop}

\noindent
\textbf{Proof}.
For given $v\in Z^\sigma _\mu (T)$, we take an extension $v^\dagger \in Z^\sigma_\mu$ of $v$ with $\| v^\dagger\| _{Z^\sigma _\mu}\leq 2\| v\| _{Z^\sigma _\mu (T)}$.
From Lemma~\ref{lem:nonlinear1}, we have
\begin{align*}
&\Big\| \int _0^tU_\mu (t-t')N_1[u_1(t'),u_2(t');v(t')]\,dt'\Big\| _{Z^\sigma _\mu (T)}\\
&\quad \leq \Big\| \psi (t)\int _0^tU_\mu (t-t')\Big( \chi _{[0,T]}(t')N_1[u_1(t'),u_2(t');v^\dagger (t')]\Big) \,dt'\Big\| _{Z^\sigma _\mu}\\
&\quad \lesssim \| u_0\|_{L^2}^{-1}\Big\| \chi _{[0,T]}(t)\Big( \| u_1(t)\|_{L^2}^2-\| u_2(t)\|_{L^2}^2\Big) v^\dagger \Big\| _{X^{\sigma +\frac{1}{2},0}}\\
&\quad \lesssim \| u_0\|_{L^2}^{-1}\| v^\dagger \| _{X^{\sigma +\frac{1}{2},0}}\sup _{t\in [0,T]}\big| \| u_1(t)\|_{L^2}^2-\| u_2(t)\|_{L^2}^2\big| \\
&\quad \lesssim \| u_0\|_{L^2}^{-2}\| v^\dagger \| _{X^{\sigma,\frac{1}{2}}_\mu}\sup _{t\in [0,T]}\big| \| u_1(t)\|_{L^2}^2-\| u_2(t)\|_{L^2}^2\big| \\
&\quad \lesssim \| u_0\|_{L^2}^{-2}\| v\| _{Z^\sigma _\mu (T)}\sup _{t\in [0,T]}\big| \| u_1(t)\|_{L^2}^2-\| u_2(t)\|_{L^2}^2\big| .
\end{align*}
Then, the first bound in (i) is verified by the assumption $\| u_j\| _{L^\infty ((0,T);L^2)}\lesssim \| u_0\| _{L^2}$, and the first bound in (ii) follows from Lemma~\ref{lem:L2diff} (ii).
To show the second bounds, replace the estimate
\[ \| v^\dagger \|_{X^{\sigma +\frac{1}{2},0}}\lesssim \| u_0\| _{L^2}^{-1}\| v\| _{Z^{\sigma}_\mu (T)}\]
in the above argument with
\[
\qquad\qquad \| \chi _{[0,T]}v \| _{X^{\sigma +\frac{1}{2},0}}=\| v\| _{L^2((0,T);H^{\sigma +\frac{1}{2}})}\lesssim T^{\frac{1}{2}}\| v\| _{L^\infty ((0,T);H^{\sigma +\frac{1}{2}})}.\qquad\quad \Box
\]

\bigskip
The next term $N_2$ is easy to treat, since there is no loss of derivative.
\begin{prop}\label{prop:estN2}
Let $\sigma \in \mathbf{R}$, $\theta \geq 0$, $\theta _1\in \mathbf{R}$, and $T\in (0,1]$.
Then, we have 
\begin{align*}
&\Big\| \int _0^tU_\mu (t-t')N_2[u_1(t'),u_2(t'),u_3(t')]\,dt'\Big\| _{Z^\sigma _\mu (T)}\\
&\quad \lesssim \frac{T^{\frac{1}{2}}}{\| u_0\|_{L^2}}\| u_1\| _{L^\infty ((0,T);H^{\theta_1})}\| u_2\| _{L^\infty ((0,T);H^{\theta -\theta _1})}\| u_3\| _{L^\infty ((0,T);H^{\sigma +\frac{1}{2}-\theta})}.
\end{align*}
\end{prop}

\noindent
\textbf{Proof}.
We first observe that 
\begin{align*}
\Big| k\big[ \mathrm{sgn}(k-k')-\mathrm{sgn}(k)\big] \Big| &= 2|k|\chi_{\{ k'>k>0\} \cup \{ k'<k<0\}}(k,k') \leq 2\langle k\rangle ^{1-\theta}\langle k'\rangle ^{\theta}
\end{align*}
for any $\theta \geq 0$, and hence 
\[
|\mathcal{F}_xN_2[u_1,u_2,u_3](k)|\lesssim \| u_1\| _{H^{\theta _1}}\| u_2 \|_{H^{\theta -\theta _1}}\langle k\rangle ^{1-\theta}|\hat{u}_3(k)|
\]
for any $\theta \geq 0$ and $\theta _1\in \mathbf{R}$.
From this bound and Lemma~\ref{lem:nonlinear1}, we have
\begin{align*}
&\Big\| \int _0^tU_\mu (t-t')N_2[u_1(t'),u_2(t'),u_3(t')]\,dt'\Big\| _{Z^\sigma _\mu (T)}\\
&\quad \leq \Big\| \psi (t)\int _0^tU_\mu (t-t')\Big\{ \chi _{[0,T]}(t')N_2[u_1(t'),u_2(t'),u_3(t')]\Big\} \,dt'\Big\| _{Z^\sigma _\mu}\\
&\quad \lesssim \frac{1}{\| u_0\|_{L^2}}\big\| \chi _{[0,T]}N_2[u_1,u_2,u_3]\big\| _{X^{\sigma -\frac{1}{2},0}}\\
&\quad \lesssim \frac{1}{\| u_0\|_{L^2}}\big\| \| u_1(t)\| _{H^{\theta _1}}\| u_2(t)\|_{H^{\theta -\theta _1}}\langle \partial _x\rangle ^{1-\theta}u_3(t)\big\| _{L^2((0,T);H^{\sigma -\frac{1}{2}})}\\
&\quad \lesssim \frac{T^{\frac{1}{2}}}{\| u_0\|_{L^2}}\| u_1\| _{L^\infty ((0,T);H^{\theta _1})}\| u_2\|_{L^\infty ((0,T);H^{\theta -\theta _1})}\| u_3\| _{L^\infty ((0,T);H^{\sigma +\frac{1}{2}-\theta})},
\end{align*}
as desired.
\hfill $\Box$

\bigskip
Let us turn to the nonlinear estimate for the last term $N_3$ in \eqref{eq:nonlinearity}.
This is done by the combination of the modulation estimate for the dispersive equation and the parabolic regularity.
More precisely, in dealing with the high-low type nonlinear interactions we can gain at most $\frac{1}{2}$ spatial derivative by the modulation estimate at the cost of $\frac{1}{2}$ integrability in time.
(This is in a sharp contrast to the nonlinearity $u^2\partial _x\bar{u}$ in DNLS (after a suitable gauge transform), for which $\frac{1}{2}$ time integrability is converted to one derivative in the high-low interactions.)
We then recover the remaining loss of $\frac{1}{2}$ derivative by the parabolic regularity.
Since $-\mu D_x$ is an elliptic operator of first order, we have to use up the remaining $\frac{1}{2}$ integrability in time of the Duhamel form.
That is why we cannot have a small factor $T^\varepsilon$ in the nonlinear estimate and hence need to assume smallness in the main theorem. 
One may think of using full time integrability for the parabolic regularity to gain one derivative.
This, however, would lead to an estimate with the constant proportional to $\mu ^{-1}\sim \| u_0\| _{L^2}^{-2}$, which would be unfavorable in the small data setting. 
We therefore make use of nonlinear smoothing effect for the dispersive equation to minimize the influence of the small coefficient on the parabolic term.

Practically, after applying Lemma~\ref{lem:inhom} we need to estimate the $Y^{s,-1}_\mu$ norm of the nonlinear term.
In the previous result \cite{MR}, there was a room in time integrability and it was allowed to replace the $Y^{s,-1}$ norm roughly with $X^{s,-\frac{1}{2}+\varepsilon}$ for a small $\varepsilon >0$ by using the Cauchy-Schwarz inequality.
In our problem, we apply the argument as used in the proof of \cite[Lemma~7.42]{Bour} to avoid such an $\varepsilon$ loss; see Lemma~\ref{L1} and the argument below for estimating the $\ell ^2_kL^1_\tau$ norm in the case {\bf (II)}.

The main part of the nonlinear estimate for $N_3$ is included in the next lemma.

\begin{lem}\label{lem:nonlinear+}
Let $\sigma \geq 0$, and let $a_j,b_j>0$, $j=1,2,3$, be such that $a_j+b_j>1/2$ for $j=1,2,3$.
Then, there exist constants $\varepsilon =\varepsilon (a_j,b_j)\in (0,1/2)$, $C_1=C_1(a_j,b_j)>0$ independent of $\sigma$ and $C_{2,\sigma}=C_{2,\sigma}(\sigma ,a_j,b_j)>0$ such that we have
\[ \| N_3[\psi _Tu_1,\psi _Tu_2,\psi _Tu_3] \| _{N^{\sigma}_\mu} \leq \frac{C_1+C_{2,\sigma}T^\varepsilon}{\| u_0\| _{L^2}}\sum _{(j,m,n)}\| u_j\| _{Z^{\sigma}_\mu}\| u_m\| _{Z^{a_j}_0}\| u_n\| _{Z^{b_j}_0},\]
where the summation is taken over $(j,m,n)\in \{ (1,2,3),(2,3,1),(3,1,2)\}$.
\end{lem}

\begin{rem}
The above estimate is a refinement of the standard trilinear estimate in $Z^{s}_\mu$ in the following sense.
(i) The function of the second largest spatial frequency is measured in $H^a$ with arbitrarily small $a>0$.
(ii) Two low-frequency functions are measured in $Z^a_0$ instead of $Z^a_\mu$.
This modification will be exploited in the proof of continuous dependence of solutions on initial data
.
(iii) The constant on the right side of the estimate satisfies either that it is independent of $\sigma$ or that it has a small factor $T^\varepsilon$.
This will be important to prove the smoothing property for $t>0$.
\end{rem}

\noindent
\textbf{Proof of Lemma~\ref{lem:nonlinear+}}. 
In the following, we write $u_j$ to denote $\psi_Tu_j$.
Note first that
\begin{align*}
&|\mathcal{F}_{t,x}N_3[u_1,u_2,u_3](\tau ,k)|\\
&\quad \lesssim |k|\sum _{\begin{smallmatrix} k=k_1-k_2+k_3\\ k_2\neq k_1,k_3\end{smallmatrix}}\int _{\tau =\tau _1-\tau _2+\tau _3}|\tilde{u}_1(\tau _1,k_1)\tilde{u}_2(\tau _2,k_2)\tilde{u}_3(\tau _3,k_3)|\,d\nu_\tau ,
\end{align*}
where $d\nu_\tau$ is the surface measure on the hyperplane $\{ (\tau_1, \tau_2, \tau_3) \in \mathbf{R}^3\, | \, \tau = \tau_1 - \tau_2 + \tau_3 \}$.
Let $\zeta := \tau + k^2$ and $\zeta_j := \tau_j + k_j^2$, $j = 1, 2, 3$, then 
\begin{equation}\label{modulation}
\begin{split}
\max \{ |\zeta |,|\zeta _1|,|\zeta _2|,|\zeta _3|\} & \geq \frac{1}{4}|\zeta -\zeta _1+\zeta _2-\zeta _3|=\frac{1}{2}|(k_1-k_2)(k_2-k_3)|\\
& \sim \langle k_1-k_2\rangle \langle k_2-k_3\rangle 
\end{split}
\end{equation}
under the restriction $\tau =\tau _1-\tau _2+\tau _3$, $k=k_1-k_2+k_3$, $k_2\neq k_1,k_3$.
We divide the analysis according to which one of $|\zeta |,|\zeta _j|$ is the biggest.

\bigskip
{\bf (I)}
We first suppose that $|\zeta_1| =\max \{ |\zeta |, |\zeta _j| \}$.
(A similar argument can be applied to the case where $|\zeta _2|$ or $|\zeta _3|$ is the biggest.)
We will show
\[ \| N_3[u_1,u_2,u_3] \| _{N^{\sigma}_\mu} \leq \frac{C+C_{\sigma}T^\varepsilon}{\| u_0\| _{L^2}}\sum _{(j,m,n)}\| u_j\| _{Z^{\sigma}_0}\| u_m\| _{Z^{a_j}_0}\| u_n\| _{Z^{b_j}_0},\]
where (and hereafter) $C$ and $C_\sigma$ stand for a positive constant independent of $\sigma$ and that dependent on $\sigma$, respectively.

In view of Lemma~\ref{L1} and \eqref{modulation}, the desired estimate follows once we prove
\begin{align}
&\bigg\| \langle k\rangle ^{\sigma +\frac{1}{2}}\sum _{k=k_1-k_2+k_3}\int _{\begin{smallmatrix}\tau =\tau _1-\tau _2+\tau _3 \\ |\zeta _1|\gtrsim \langle k_1-k_2\rangle \langle k_2-k_3\rangle\end{smallmatrix}}|\tilde{u}_1(\tau _1,k_1)\tilde{u}_2(\tau _2,k_2)\tilde{u}_3(\tau _3,k_3)|\,d\nu_\tau \bigg\| _{\ell ^2_kL^2_\tau}\notag \\
&\quad \leq \big( C+C_\sigma T^\varepsilon \big) \sum _{(j,m,n)}\| u_j\| _{Z^{\sigma}_0}\| u_m\| _{Z^{a_j}_0}\| u_n\| _{Z^{b_j}_0}.\label{claim-1}
\end{align}
We assume $|k_1|\geq |k_3|$ by symmetry and consider the following two cases separately.

\medskip {\bf Case A}: $\langle k_1\rangle \sim \langle k_2\rangle \sim \langle k_3\rangle$.\quad [e.g., $10\langle k_3\rangle \geq \langle k_2\rangle \geq \frac{1}{10}\langle k_1\rangle$]

\noindent
For given $a_1, b_1>0$ with $a_1+b_1>\frac{1}{2}$, we choose $\varepsilon \in (0,1/2)$ so that $4\varepsilon \leq a_1+b_1-1/2$.
Then, the left side of \eqref{claim-1} is bounded by
\begin{align*}
C_\sigma \bigg\| \sum _{k=k_1-k_2+k_3}\int _{\tau =\tau _1-\tau _2+\tau _3} \frac{\langle \zeta _1\rangle ^{\frac{1}{2}-\varepsilon}\langle k_1\rangle ^\sigma \langle k_2\rangle ^{a_1-2\varepsilon}\langle k_3\rangle ^{b_1-2\varepsilon}}{\langle k_1-k_2\rangle ^{\frac{1}{2}-\varepsilon}\langle k_2-k_3\rangle ^{\frac{1}{2}-\varepsilon}}& \\
\times |\tilde{u}_1(\tau _1,k_1)\tilde{u}_2(\tau _2,k_2)\tilde{u}_3(\tau _3,k_3)|\,d\nu_\tau \bigg\| _{\ell ^2_kL^2_\tau}&
\end{align*}
Applying Minkowski's and Young's inequalities in $\tau$, and then Cauchy-Schwarz inequality in $(k_1,k_2,k_3)$, this is bounded by
\begin{align*}
&C_\sigma \bigg\| \sum _{k=k_1-k_2+k_3}\frac{\langle k_1\rangle ^\sigma \langle k_2\rangle ^{a_1-2\varepsilon}\langle k_3\rangle ^{b_1-2\varepsilon}}{\langle k_1-k_2\rangle ^{\frac{1}{2}-\varepsilon}\langle k_2-k_3\rangle ^{\frac{1}{2}-\varepsilon}} \\
&\qquad \qquad \times \| \langle \zeta _1\rangle ^{\frac{1}{2}-\varepsilon}\tilde{u}_1(\cdot ,k_1)\| _{L^2_\tau} \| \tilde{u}_2(\cdot ,k_2)\|_{L^1_\tau} \| \tilde{u}_3(\cdot ,k_3)\|_{L^1_\tau} \bigg\| _{\ell ^2_k}\\
&\leq C_\sigma \sup _{k\in \mathbf{Z}}\Big[ \sum _{k_2,k_3\in \mathbf{Z}}\frac{1}{\langle k-k_3\rangle ^{1-2\varepsilon}\langle k_2-k_3\rangle ^{1-2\varepsilon}\langle k_2\rangle ^{4\varepsilon}\langle k_3\rangle ^{4\varepsilon}}\Big] ^{\frac{1}{2}} \\
&\qquad \qquad \qquad \qquad \times \| u_1\|_{X^{\sigma ,\frac{1}{2}-\varepsilon}} \| u_2\|_{Y^{a_1,0}} \| u_3\|_{Y^{b_1,0}}.
\end{align*}
Since the sum on the right side is finite and $\| \psi _Tu_1\| _{X^{\sigma ,\frac{1}{2}-\varepsilon}}\lesssim T^\varepsilon \| u_1\|_{Z^{\sigma}_0}$ by Lemma~\ref{lem:linear2}, the above quantity is estimated by $C_\sigma T^\varepsilon \| u_1\| _{Z^\sigma_0}\| u_2\|_{Y^{a_1,0}}\| u_3\|_{Y^{b_1,0}}$, and we have \eqref{claim-1}.

\medskip {\bf Case B}: The other cases.

\noindent
There are four possibilities:
\begin{align*}
\text{(a)} \quad &\langle k_2\rangle \gg \langle k_1\rangle \geq \langle k_3\rangle,\quad [\text{e.g.,}~\langle k_2\rangle \geq 10\langle k_1\rangle]\\
\text{(b)} \quad &\langle k_1\rangle \geq \langle k_3\rangle \gg \langle k_2\rangle,\quad [\text{e.g.,}~\tfrac{1}{10}\langle k_3\rangle \geq \langle k_2\rangle]\\
\text{(c)} \quad &\langle k_1\rangle \gg \langle k_2\rangle \gtrsim \langle k_3\rangle,\quad [\text{e.g.,}~\tfrac{1}{10}\langle k_1\rangle \geq \langle k_2\rangle \geq \tfrac{1}{10}\langle k_3\rangle]\\
\text{(d)} \quad &\langle k_1\rangle \gtrsim \langle k_2\rangle \gg \langle k_3\rangle.\,\quad [\text{e.g.,}~10\langle k_1\rangle \geq \langle k_2\rangle \geq 10\langle k_3\rangle]
\end{align*}
We see that each of them implies one of the following three conditions:
\begin{align*}
\text{(i)}\quad &\langle k_1-k_2\rangle \langle k_2-k_3\rangle \gtrsim \langle k\rangle \langle k_2-k_3\rangle \quad \text{and}\quad \langle k\rangle \lesssim \langle k_1\rangle ,\\
\text{(ii)}\quad &\langle k_1-k_2\rangle \langle k_2-k_3\rangle \gtrsim \langle k\rangle \langle k_1+k_3\rangle \quad \text{and}\quad \langle k\rangle \lesssim \langle k_2\rangle ,\\
\text{(iii)}\quad &\langle k_1-k_2\rangle \langle k_2-k_3\rangle \gtrsim \langle k\rangle \langle k_1-k_2\rangle \quad \text{and}\quad \langle k\rangle \lesssim \langle k_3\rangle .
\end{align*}
In fact, it is easy to see that (a) implies (ii), while (b) or (c) implies (i).
As for (d), we further split the region into $\langle k\rangle \lesssim \langle k_3\rangle$ and $\langle k\rangle \gg \langle k_3\rangle$.
In the former situation we have (iii), while in the latter we have $\langle k\rangle \sim \langle k-k_3\rangle =\langle k_1-k_2\rangle \lesssim \langle k_1\rangle$, and hence (i).

Let us focus on the case (iii), for instance, though the following argument can be applied to (i) and (ii) as well.
The case (iii) is further divided into 
\[ \text{(iii-1)}\quad |k|>2^{1/\sigma}|k_3|,\qquad \text{(iii-2)}\quad |k|\leq 2^{1/\sigma}|k_3|.\]
(When $\sigma =0$, we do not divide the case and follow the argument for (iii-2) below.)

In the case (iii-1), we have $\max \{ |k_1|,|k_2|\} > \frac{1}{2}(2^{1/\sigma}-1)|k_3|$.
Hence, it holds inside the norm on the left side of \eqref{claim-1} that
\begin{align*}
\langle k\rangle ^{\sigma +\frac{1}{2}}
&\leq C_\sigma \frac{\langle \zeta_1\rangle ^{\frac{1}{2}-\varepsilon}}{\langle k_1-k_2\rangle ^{\frac{1}{2}-\varepsilon}\langle k_2-k_3\rangle ^{\frac{1}{2}-\varepsilon}} \langle k\rangle ^{\frac{1}{2}-\varepsilon}\langle k_3\rangle ^{\sigma} \max \{ \langle k_1\rangle ,\, \langle k_2\rangle \} ^\varepsilon \\
&\leq C_\sigma \frac{\langle \zeta_1\rangle ^{\frac{1}{2}-\varepsilon}\langle k_1\rangle ^\varepsilon \langle k_2\rangle ^\varepsilon \langle k_3\rangle^{\sigma}}{\langle k_1-k_2\rangle ^{\frac{1}{2}-\varepsilon}}.
\end{align*}
Application of Minkowski's and Young's inequalities in $\tau$ similar to {\bf Case A} shows that the left side of \eqref{claim-1} is bounded by
\[ C_\sigma \bigg\| \bigg[ \frac{1}{\langle \cdot \rangle ^{\frac{1}{2}-\varepsilon}}\sum _{k_1-k_2=\,\cdot}\langle k_1\rangle ^\varepsilon f_1(k_1)\langle k_2\rangle ^\varepsilon f_2(k_2)\bigg] * \big( \langle \cdot \rangle ^\sigma f_3\big) \bigg\| _{\ell ^2_k},\]
where we have written
\[ f_1(k):=\| \langle \tau +k^2\rangle ^{\frac{1}{2}-\varepsilon} \tilde{u}_1(\tau ,k)\| _{L^2_\tau},\qquad f_j(k):=\| \tilde{u}_j(\tau ,k)\|_{L^1_\tau}\quad (j=2,3). \]
By Young's and H\"older's inequalities in $k$, this is evaluated by
\begin{align*}
&C_\sigma \| \langle \cdot \rangle ^{\sigma}f_3\| _{\ell ^2}\| \langle \cdot \rangle ^{-\frac{1}{2}+\varepsilon}\| _{\ell ^r}\| \langle \cdot \rangle ^{\varepsilon}f_1\| _{\ell ^{p}}\| \langle \cdot \rangle ^{\varepsilon}f_2\| _{\ell ^q}\\
&\quad \leq C_\sigma \| \langle \cdot \rangle ^{\sigma}f_3\| _{\ell ^2}\| \langle \cdot \rangle ^{a_3}f_1\| _{\ell ^2}\| \langle \cdot \rangle ^{b_3}f_2\| _{\ell ^2}\\
&\quad = C_\sigma \| u_3\| _{Y^{\sigma ,0}}\| u_1\|_{X^{a_3,\frac{1}{2}-\varepsilon}}\| u_2\| _{Y^{b_3,0}},
\end{align*}
where for given $a_3,b_3>0$ with $a_3+b_3>\frac{1}{2}$ we have chosen $1\leq p,q<2$ and $r>2$ such that $a_3-\varepsilon >\frac{1}{p}-\frac{1}{2}$, $b_3-\varepsilon >\frac{1}{q}-\frac{1}{2}$, $\frac{1}{2}-\varepsilon >\frac{1}{r}$, and $\frac{1}{p}+\frac{1}{q}+\frac{1}{r}=2$.
This is possible if $\varepsilon$ satisfies $0<\varepsilon <\min \{ a_3,b_3,\frac{1}{2},\frac{1}{3}(a_3+b_3-\frac{1}{2})\}$.
We conclude the desired bound \eqref{claim-1} by applying Lemma~\ref{lem:linear2} to derive a factor $T^\varepsilon$ from the norm of $u_1$.

When (iii-2) holds, we have $\langle k\rangle ^\sigma \leq 2\langle k_3\rangle ^\sigma$, which leads to the bound
\[ \langle k\rangle ^{\sigma +\frac{1}{2}}\leq C\frac{\langle \zeta _1\rangle ^{\frac{1}{2}}}{\langle k_1-k_2\rangle ^{\frac{1}{2}}\langle k_2-k_3\rangle ^{\frac{1}{2}}}\langle k\rangle ^{\frac{1}{2}}\langle k_3\rangle ^{\sigma}\leq C\frac{\langle \zeta _1\rangle ^{\frac{1}{2}}\langle k_3\rangle ^{\sigma}}{\langle k_1-k_2\rangle ^{\frac{1}{2}}}.\]
An argument almost parallel to the case (iii-1) then yields the bound 
\[ \text{L.H.S. of \eqref{claim-1}} \leq C\| u_3\| _{Y^{\sigma ,0}}\| u_1\|_{X^{a_3,\frac{1}{2}}}\| u_2\| _{Y^{b_3,0}},\]
where the factor $T^\varepsilon$ does not appear but the constant can be made independent of $\sigma$.
This finishes the case {\bf (I)}.

\bigskip
{\bf (II)} 
We next suppose that $|\zeta | = \max \{ |\zeta |, |\zeta _j|\}$.
In this case, we will show
\begin{align*}
&\| N_3[u_1,u_2,u_3] \| _{N^{\sigma}_0}\\
&\quad \leq \Big( \frac{C}{\| u_0\|_{L^2}}+\frac{C_{\sigma}T^\varepsilon}{\| u_0\| _{L^2}^{1-2\varepsilon}}\Big) \sum _{(j,m,n)}\| u_j\| _{Z^{\sigma}_\mu}\| u_m\| _{Y^{a_j,0}}\| u_n\| _{Y^{b_j,0}}.
\end{align*}

In view of \eqref{modulation}, it suffices to estimate
\[ \bigg\| \frac{\langle k\rangle ^{\sigma +1}}{\langle \zeta \rangle ^{\frac{1}{2}}}\sum _{\begin{smallmatrix} k=k_1-k_2+k_3\\ |\zeta |\gtrsim \langle k_1-k_2\rangle \langle k_2-k_3\rangle \end{smallmatrix}}\int _{\tau =\tau _1-\tau _2+\tau _3}|\tilde{u}_1(\tau _1,k_1)\tilde{u}_2(\tau _2,k_2)\tilde{u}_3(\tau _3,k_3)|\,d\nu _\tau \bigg\| _{\ell ^2_kL^2_\tau}
\]
and
\[ \bigg\| \frac{\langle k\rangle ^{\sigma +1}}{\langle \zeta \rangle} \sum _{\begin{smallmatrix} k=k_1-k_2+k_3\\ |\zeta |\gtrsim \langle k_1-k_2\rangle \langle k_2-k_3\rangle \end{smallmatrix}}\int _{\tau =\tau _1-\tau _2+\tau _3}|\tilde{u}_1(\tau _1,k_1)\tilde{u}_2(\tau _2,k_2)\tilde{u}_3(\tau _3,k_3)|\,d\nu _\tau \bigg\| _{\ell ^2_kL^1_\tau}.
\] 
By the Minkowski inequality and the Cauchy-Schwarz inequality in $\tau$ (similarly to the proof of Lemma~\ref{L1}), both of the above norms are bounded by
\begin{equation}\label{claim-2}
\begin{split}
&C\ \bigg\| \sum _{k=k_1-k_2+k_3}\frac{\langle k\rangle ^{\sigma +1}}{\langle k_1-k_3\rangle ^{\frac{1}{2}}\langle k_2-k_3\rangle^{\frac{1}{2}}}\\
&\qquad\qquad \times \big\| \int _{\tau =\tau _1-\tau _2+\tau _3}|\tilde{u}_1(\tau _1,k_1)\tilde{u}_2(\tau _2,k_2)\tilde{u}_3(\tau _3,k_3)|\,d\nu _\tau \big\| _{L^2_\tau} \bigg\| _{\ell ^2_k}.
\end{split}
\end{equation}
We consider the same case separation as for {\bf (I)}.

For {\bf Case A}, we choose $\varepsilon \in (0,\frac12 )$ as $3\varepsilon \leq a_1+b_1-\frac12$, so that $\langle k\rangle ^{\sigma +1}\lesssim _\sigma \langle k_1\rangle ^{\sigma +\frac{1}{2}-\varepsilon} \langle k_2\rangle ^{a_1-\varepsilon}\langle k_3\rangle ^{b_1-\varepsilon}$ and, by the Young inequality,
\begin{align*}
\eqref{claim-2}&\leq C_\sigma \bigg\| \sum _{k=k_1-k_2+k_3}\frac{1}{\langle k_1-k_2\rangle ^{\frac{1}{2}}\langle k_2-k_3\rangle^{\frac{1}{2}}\langle k_2\rangle ^{\varepsilon}\langle k_3\rangle ^{\varepsilon}}\\
&\qquad \times \| \langle k_1\rangle ^{\sigma +\frac{1}{2}-\varepsilon}\tilde{u}_1(\cdot ,k_1)\|_{L^2_\tau}\| \langle k_2\rangle ^{a_1}\tilde{u}_2(\cdot ,k_2)\|_{L^1_\tau}\| \langle k_3\rangle ^{b_1}\tilde{u}_3(\cdot ,k_3)\|_{L^1_\tau} \bigg\| _{\ell ^2_k}.
\end{align*}
By the argument as for {\bf (I), Case A} and the estimate from Lemma~\ref{lem:linear2}:
\[ \| \psi _Tu_1\| _{X^{\sigma +\frac{1}{2}-\varepsilon ,0}}\lesssim \| u_0\|_{L^2}^{-1+2\varepsilon} \| \psi _Tu_1\|_{X^{\sigma ,\frac{1}{2}-\varepsilon}_\mu} \lesssim \| u_0\|_{L^2}^{-1+2\varepsilon}T^\varepsilon \| u_1\|_{Z^{\sigma}_\mu},\]
we obtain the desired bound.

For {\bf Case B} (iii-1), it holds that
\[ \frac{\langle k\rangle ^{\sigma +1}}{\langle k_1-k_2\rangle ^{\frac12}\langle k_2-k_3\rangle ^{\frac12}} \leq C_\sigma \frac{\langle k_3\rangle^{\sigma +\frac{1}{2}}}{\langle k_1-k_2\rangle ^{\frac12}}\leq C_\sigma \frac{\langle k_3\rangle^{\sigma +\frac{1}{2}-\varepsilon }\langle k_1\rangle ^\varepsilon \langle k_2\rangle ^\varepsilon}{\langle k_1-k_2\rangle ^{\frac12}},\]
and hence
\begin{align*}
\eqref{claim-2}&\leq C_\sigma \bigg\| \sum _{k=k_1-k_2+k_3}\frac{1}{\langle k_1-k_2\rangle ^{\frac{1}{2}}}\\
&\qquad \times \| \langle k_1\rangle ^{\varepsilon}\tilde{u}_1(\cdot ,k_1)\|_{L^1_\tau}\| \langle k_2\rangle ^{\varepsilon}\tilde{u}_2(\cdot ,k_2)\|_{L^1_\tau}\| \langle k_3\rangle ^{\sigma +\frac{1}{2}-\varepsilon}\tilde{u}_3(\cdot ,k_3)\|_{L^2_\tau} \bigg\| _{\ell ^2_k}.
\end{align*}
Then, the argument as for {\bf (I), Case B} (iii-1) and Lemma~\ref{lem:linear2} implies the claimed estimate.
For (iii-2), we instead apply the inequality with $\sigma$-independent constant
\[ \frac{\langle k\rangle ^{\sigma +1}}{\langle k_1-k_2\rangle ^{\frac12}\langle k_2-k_3\rangle ^{\frac12}} \leq C \frac{\langle k_3\rangle^{\sigma +\frac{1}{2}}}{\langle k_1-k_2\rangle ^{\frac12}}\]
and the estimate $\| u_3\| _{X^{\sigma +\frac{1}{2},0}}\lesssim \| u_0\|_{L^2}^{-1} \| u_3\|_{X^{\sigma ,\frac{1}{2}}_\mu}$.
The cases (i) and (ii) can be treated in the same way.

This completes the proof of Lemma~\ref{lem:nonlinear+}.\hfill $\Box$

\begin{prop}\label{prop:estN3}
Let $\sigma $, $a_j,b_j$ ($j=1,2,3$) be as in Lemma~\ref{lem:nonlinear+}.
We have
\begin{align*}
&\Big\| \int _0^t U_\mu (t-t')N_3[u_1(t'),u_2(t'),u_3(t')]\,dt' \Big\| _{Z^\sigma _\mu (T)}\\
&\quad \leq \frac{C_1+C_{2,\sigma}T^\varepsilon}{\| u_0\| _{L^2}}\sum _{(j,m,n)}\| u_j\| _{Z^{\sigma}_\mu (T)}\| u_m\| _{Z^{a_j}_0(T)}\| u_n\| _{Z^{b_j}_0(T)}
\end{align*}
for any $0<T\leq 1$, where $C_1,C_{2,\sigma}>0$ and $0<\varepsilon <1/2$ are constants depending on $a_j,b_j$, and $C_{2,\sigma}$ depends also on $\sigma$ but $C_1$ does not.
\end{prop}

\noindent
\textbf{Proof}.
As we have seen in the proof of Lemma~\ref{lem:nonlinear+}, the set
\[ \Omega :=\big\{ (k,k_1,k_2,k_3)\in \mathbf{Z}^4\,\big|\,k=k_1-k_2+k_3\big\} \]
can be classified as $\Omega =\bigcup _{j=0}^3\Omega _j$;
\begin{gather*}
\Omega _0:=\big\{ (k,k_1,k_2,k_3)\in \Omega \,\big|\, \langle k_1\rangle \sim \langle k_2\rangle \sim \langle k_3\rangle \big\} ,\\
\Omega _1:=\big\{ (k,k_1,k_2,k_3)\in \Omega \,\big|\, \langle k_1-k_2\rangle \langle k_2-k_3\rangle \gtrsim \langle k\rangle \langle k_2-k_3\rangle ,\; \langle k\rangle \lesssim \langle k_1\rangle \big\} ,\\
\Omega _2:=\big\{ (k,k_1,k_2,k_3)\in \Omega \,\big|\, \langle k_1-k_2\rangle \langle k_2-k_3\rangle \gtrsim \langle k\rangle \langle k_1+k_3\rangle ,\; \langle k\rangle \lesssim \langle k_2\rangle \big\} ,\\
\Omega _3:=\big\{ (k,k_1,k_2,k_3)\in \Omega \,\big|\, \langle k_1-k_2\rangle \langle k_2-k_3\rangle \gtrsim \langle k\rangle \langle k_1-k_2\rangle ,\; \langle k\rangle \lesssim \langle k_3\rangle \big\} .
\end{gather*}
We decompose $N_3$ into $N_{3,1}+N_{3,2}+N_{3,3}$, in which the sum with respect to $(k_1,k_2,k_3)$ is restricted to $\Omega _0\cup \Omega _1$, $\Omega _2\setminus (\Omega_0\cup \Omega_1)$, and $\Omega _3\setminus (\Omega_0\cup \Omega _1\cup \Omega _2)$, respectively.
Then, the proof of Lemma~\ref{lem:nonlinear+} actually verifies that 
\[ \big\| N_{3,j}[\psi _Tu_1,\psi _Tu_2,\psi _Tu_3]\big\| _{N^{\sigma}_\mu}\leq \frac{C+C_{\sigma}T^\varepsilon}{\| u_0\|_{L^2}}\| u_j\| _{Z^{\sigma}_\mu}\| u_m\| _{Z^{a_j}_0}\| u_n\| _{Z^{b_j}_0}\]
for $(j,m,n)\in \{ (1,2,3),(2,3,1),(3,1,2)\}$.

Let us fix $(j,m,n)$, and take extensions $u^\dagger_1,u^\dagger_2,u^\dagger_3$ of $u_1,u_2,u_3$ defined on $[0,T]\times \mathbf{T}$ such that 
\[ \| u^\dagger_j\|_{Z^{\sigma}_\mu}\leq 2\| u_j\|_{Z^{\sigma}_\mu (T)},\quad
\| u^\dagger_m\|_{Z^{a_j}_0}\leq 2\| u_m\|_{Z^{a_j}_0(T)},\quad
\| u^\dagger_n\|_{Z^{b_j}_0}\leq 2\| u_n\|_{Z^{b_j}_0(T)}.\]
(Note that the choice of $u^\dagger_1,u^\dagger_2,u^\dagger_3$ may depend on $(j,m,n)$.)
Using Lemma~\ref{lem:inhom} and the above estimate for $N_{3,j}$, we have
\begin{align*}
&\Big\| \int _0^t U_\mu (t-t')N_{3,j}[u_1(t'),u_2(t'),u_3(t')]\,dt' \Big\| _{Z^\sigma _\mu (1)}\\
&\quad \leq \Big\| \psi (t)\int _0^t U_\mu (t-t')N_{3,j}[\psi _T(t')u^\dagger_1(t'),\psi _T(t')u^\dagger_2(t'),\psi _T(t')u^\dagger_3(t')]\,dt' \Big\| _{Z^\sigma _\mu}\\
&\quad \leq \frac{C+C_{\sigma}T^\varepsilon}{\| u_0\| _{L^2}}\| u^\dagger_j\| _{Z^{\sigma}_\mu}\| u^\dagger_m\| _{Z^{a_j}_0}\| u^\dagger_n\| _{Z^{b_j}_0}\\
&\quad \leq 8\frac{C+C_{\sigma}T^\varepsilon}{\| u_0\| _{L^2}}\| u_j\| _{Z^{\sigma}_\mu (1)}\| u_m\| _{Z^{a_j}_0(T)}\| u_n\| _{Z^{b_j}_0(T)}.
\end{align*}
Summing up in $j$, we obtain the claimed estimate.
\hfill $\Box$


\section{Proof of existence, uniqueness, regularity}
\label{sec:local}

We shall first establish the contraction properties in a suitable complete metric space.
\begin{prop}\label{prop:contraction}
Let $\alpha \in \mathbf{R}$, $\beta <0$ be fixed.
Let $s\geq s_0>1/2$.
Then, there exist $\eta =\eta (s_0)>0$ and $T_0=T_0(s,s_0)\in (0,1]$ with $T_0(s_0,s_0)=1$ such that if $u_0\in H^s(\mathbf{T})\setminus \{ 0\}$, $\| u_0\| _{H^{s_0}}\leq \eta$ and $T\in (0,T_0]$, then (with $\mu :=\frac{|\beta|}{2\pi}\| u_0\|_{L^2}^2$) the map 
\[ u~\mapsto ~U_\mu (t)u_0+\int _0^t U_\mu (t-t')\big\{ N_1[u,u_0;u]+N_2[u,u,u]+N_3[u,u,u]\big\} (t')\,dt'\]
is a contraction on the complete metric space $(\mathscr{B}^s_{u_0}(T),d(\cdot ,\cdot ))$;
\begin{align}
&\mathscr{B}^s_{u_0}(T):=\left\{ u\in Z^{s}_\mu (T)\,\left|\, 
\begin{matrix}
\| u\| _{Z^\sigma _\mu (T)}\le 2C_0\| u_0\| _{H^\sigma}~(\sigma =s,s_0,\frac{s_0}{2}),\\[5pt]
\| u-U_\mu (t)u_0\| _{Z^0_\mu (T)}\leq \delta \| u_0\| _{L^2}
\end{matrix}
\right\} \right., \label{def:B}\\[5pt]
&d(u,\tilde{u}):=\max _{\sigma =s,s_0,\frac{s_0}{2}}\frac{\| u_0\| _{H^{s}}}{\| u_0\| _{H^\sigma}}\| u-\tilde{u}\| _{Z^\sigma _\mu (T)}+\frac{\| u_0\| _{H^s}}{\delta \| u_0\| _{L^2}}\| u-\tilde{u}\| _{Z^0_\mu (T)}, \nonumber
\end{align}
where $C_0>0$ is the constant given in Lemma~\ref{lem:linear1} and $\delta >0$ is a small constant independent of $s_0,s$.

As a consequence, (letting $s_0=s>1/2$) we can define the solution map $S:u_0\mapsto u$ for the renormalized KDNLS equation \eqref{r-kdnls} on $\{ u_0\in H^{s}\,|\,\| u_0\|_{H^{s}}\leq \eta (s) \}$ into $C([0,1];H^{s})\cap L^2((0,1);H^{s+\frac12})$, complementing it by $S(0):=0$.
\end{prop}

\begin{rem}
The small factor $\delta$ is needed only to control $N_1$.
Since the constant in Proposition~\ref{prop:estN1} is independent of the regularity index, we can take $\delta$ independently of $s_0$ and $s$.
It is also worth noticing that we can choose $\eta =\eta (s_0)$ independently of $s$, because the constant $C_1$ in Proposition~\ref{prop:estN3} is independent of the regularity index $\sigma$. 
\end{rem}

\noindent
\textbf{Proof of Proposition~\ref{prop:contraction}}.
In the following, the implicit constants do not depend on $s_0,s$ and $T$, unless otherwise mentioned.
The linear part $U_\mu (t)u_0$ can be estimated by Lemma~\ref{lem:linear1}.

\medskip
\noindent
\textbf{Estimate on $N_1[u,u_0;u]$}.
For simplicity, we write $N_1:=N_1[u,u_0;u]$ and $\tilde{N}_1:=N_1[\tilde{u},u_0;\tilde{u}]$ for $u,\tilde{u}\in \mathscr{B}^s_{u_0}(T)$.

Since $\| f\| _{L^\infty (\mathbf{R};L^2)}\leq \| f\| _{Y^{0,0}}$, $u\in \mathscr{B}^s_{u_0}(T)$ implies $\| u-U_\mu (t)u_0\| _{L^\infty ((0,T);L^2)}\leq \delta \| u_0\|_{L^2}$, and thus $\| u\| _{L^\infty ((0,T);L^2)}\lesssim \| u_0\|_{L^2}$.
By Proposition~\ref{prop:estN1} (ii), we have
\begin{align*}
\Big\| \int _0^t U(t-t')N_1(t')\,dt' \Big\| _{Z^\sigma _\mu (T)}
&\lesssim \Big( \delta + T^{\frac{1}{2}}\| u_0\| _{H^{\frac{1}{2}}}\Big) \| u\| _{Z^\sigma _\mu (T)}\\
&\lesssim \Big( \delta + T^{\frac{1}{2}}\| u_0\| _{H^{\frac{1}{2}}}\Big) \| u_0\| _{H^\sigma}
\end{align*}
for $\sigma =s,s_0,\frac{s_0}{2}$, and
\begin{align*}
\Big\| \int _0^t U(t-t')N_1(t')\,dt' \Big\| _{Z^0_\mu (T)}
&\lesssim \Big( \delta + T^{\frac{1}{2}}\| u_0\| _{H^{\frac{1}{2}}}\Big) T^{\frac{1}{2}}\| u\| _{L^\infty ((0,T);H^{\frac{1}{2}})}\| u_0\| _{L^2}\\
&\lesssim \Big( \delta + T^{\frac{1}{2}}\| u_0\| _{H^{\frac{1}{2}}}\Big) T^{\frac{1}{2}}\| u_0\| _{H^{s_0}}\| u_0\| _{L^2},
\end{align*}
where at the last inequality we have estimated as
\[ \| u\| _{L^\infty ((0,T);H^{\frac{1}{2}})}\leq \| u\| _{L^\infty ((0,T);H^{s_0})}\leq \| u\|_{Z^{s_0}_\mu(T)}\leq 2C_0\| u_0\| _{H^{s_0}}.\]
Hence, if $\delta$ and $T$ satisfies
\begin{equation}\label{cond:N1}
T^{\frac{1}{2}}\| u_0\| _{H^{s_0}}\leq \delta \ll 1,
\end{equation}
then for any $u\in \mathscr{B}^s_{u_0}(T)$ we have
\[ \Big\| \int _0^t U_\mu (t-t')N_1(t')\,dt' \Big\| _{Z^\sigma _\mu (T)}\leq \frac{C_0}{10}\| u_0\| _{H^\sigma},\qquad \sigma =s,s_0,\frac{s_0}{2} \]
and
\[ \Big\| \int _0^t U_\mu (t-t')N_1(t')\,dt' \Big\| _{Z^0_\mu (T)}\leq \frac{\delta}{10}\| u_0\| _{L^2}.\]

Next, we estimate the difference.
Let $u,\tilde{u}\in \mathscr{B}^s_{u_0}(T)$, and note that
\[ N_1-\tilde{N}_1=N_1[u,\tilde{u};u]+N_1[\tilde{u},u_0;u-\tilde{u}].\]
Proposition~\ref{prop:estN1} (i), (ii) imply that
\begin{align*}
&\Big\| \int _0^t U_\mu (t-t')\left[ N_1-\tilde{N}_1\Big] (t')\,dt' \right\| _{Z^\sigma _\mu (T)}\\
&\quad \lesssim \frac{\| u\|_{Z^\sigma _\mu (T)}}{\| u_0\|_{L^2}}\| u-\tilde{u}\|_{L^\infty ((0,T);L^2)}+\left( \delta + T^{\frac{1}{2}}\| u_0\| _{H^{\frac{1}{2}}}\right) \| u-\tilde{u}\| _{Z^\sigma _\mu (T)}\\
&\quad \lesssim \frac{\| u_0\|_{H^\sigma}}{\| u_0\|_{L^2}}\| u-\tilde{u}\|_{Z^0_\mu (T)}+\left( \delta + T^{\frac{1}{2}}\| u_0\| _{H^{\frac{1}{2}}}\right) \| u-\tilde{u}\| _{Z^\sigma _\mu (T)}\\
&\quad \lesssim \left( \delta + T^{\frac{1}{2}}\| u_0\| _{H^{\frac{1}{2}}}\right) \frac{\| u_0\|_{H^\sigma}}{\| u_0\|_{H^s}}d(u,\tilde{u})
\end{align*}
for $\sigma =s,s_0,\frac{s_0}{2}$, and that
\begin{align*}
&\Big\| \int _0^t U_\mu (t-t')\left[ N_1-\tilde{N}_1\right] (t')\,dt' \Big\| _{Z^0_\mu (T)}\\
&\quad \lesssim T^{\frac{1}{2}}\| u\|_{L^\infty ((0,T);H^{\frac{1}{2}})}\| u-\tilde{u}\|_{L^\infty ((0,T);L^2)}+\left( \delta + T^{\frac{1}{2}}\| u_0\| _{H^{\frac{1}{2}}}\right) \| u-\tilde{u}\|_{Z^0_\mu (T)}\\
&\quad \lesssim \left( \delta + T^{\frac{1}{2}}\| u_0\| _{H^{s_0}}\right) \| u-\tilde{u}\|_{Z^0_\mu (T)}\\
&\quad \lesssim \left( \delta + T^{\frac{1}{2}}\| u_0\| _{H^{s_0}}\right) \frac{\delta \| u_0\|_{L^2}}{\| u_0\|_{H^s}}d(u,\tilde{u}).
\end{align*}
Again, under the condition \eqref{cond:N1}, we have
\begin{align*}
\max _{\sigma =s,s_0,\frac{s_0}{2}}\frac{\| u_0\| _{H^s}}{\| u_0\|_{H^\sigma}}\Big\| \int _0^t U_\mu (t-t')\left[ N_1-\tilde{N}_1\right] (t')\,dt' \Big\| _{Z^\sigma _\mu (T)}& \\
\qquad +\frac{\| u_0\| _{H^s}}{\delta \| u_0\|_{L^2}}\Big\| \int _0^t U_\mu (t-t')\left[ N_1-\tilde{N}_1\right] (t')\,dt' \Big\| _{Z^0_\mu (T)}&\leq \frac{1}{10}d(u,\tilde{u}).
\end{align*}
This concludes the estimate on $N_1$.

\bigskip\noindent
\textbf{Estimate on $N_2[u,u,u]$}.

We first apply Proposition~\ref{prop:estN2} with $(\theta ,\theta _1)=(\frac{1}{2},0)$ to see, for any $u\in \mathscr{B}^s_{u_0}(T)$, that
\begin{align*}
&\Big\| \int _0^t U_\mu (t-t')N_2(t')\,dt' \Big\| _{Z^\sigma _\mu (T)}\\
&\quad \lesssim \frac{T^{\frac{1}{2}}}{\| u_0\| _{L^2}}\| u\| _{L^\infty ((0,T);L^2)}\| u\| _{L^\infty ((0,T);H^{\frac{1}{2}})}\| u\| _{L^\infty ((0,T);H^\sigma )}\\
&\quad \lesssim \frac{T^{\frac{1}{2}}}{\| u_0\| _{L^2}}\| u\| _{Z^0_0(T)}\| u\| _{Z^{s_0}_0(T)}\| u\| _{Z^\sigma _0(T)}\\
&\quad \lesssim T^{\frac{1}{2}}\| u_0\|_{H^{s_0}}\| u_0\| _{H^\sigma },\qquad\qquad \sigma =s,s_0,\frac{s_0}{2},0.
\end{align*}
For the difference, we note that
\[ N_2-\tilde{N}_2=N_2[u-\tilde{u},u,u]+N_2[\tilde{u},u-\tilde{u},u]+N_2[\tilde{u},\tilde{u},u-\tilde{u}] .\]
Using Proposition~\ref{prop:estN2} for these terms with $(\theta ,\theta _1)=(\sigma +\frac{1}{2},\sigma)$, $(\sigma +\frac{1}{2},\frac{1}{2})$ and $(\frac{1}{2},0)$, respectively, we obtain
\begin{align*}
&\Big\| \int _0^t U_\mu (t-t')\Big[ N_2-\tilde{N}_2\Big] (t')\,dt' \Big\| _{Z^\sigma _\mu (T)}\\
&\quad \lesssim \frac{T^{\frac{1}{2}}}{\| u_0\|_{L^2}}\| u,\tilde{u}\| _{Z^0_0(T)}\| u,\tilde{u}\| _{Z^{s_0}_0(T)}\| u-\tilde{u}\| _{Z^\sigma _0(T)}\\
&\quad \lesssim T^{\frac{1}{2}}\| u_0\| _{H^{s_0}}\| u-\tilde{u}\| _{Z^\sigma _\mu (T)},\qquad\qquad \sigma =s,s_0,\frac{s_0}{2},0.
\end{align*}
These estimates are sufficient if $\delta$ and $T$ satisfy
\begin{equation}\label{cond:N2}
T^{\frac{1}{2}}\| u_0\| _{H^{s_0}}\ll \delta \leq 1.
\end{equation}

\begin{rem}
All the estimates so far hold in the limiting case $s_0=1/2$ as well.
Moreover, we have not needed smallness of the initial data.
In estimating $N_3$ below, we will have to assume $\| u_0\| _{H^{s_0}}\ll \delta$.
This smallness assumption, however, allows us to take $T=1$ if $s=s_0$.
\end{rem}

\noindent
\textbf{Estimate on $N_3[u,u,u]$}.
In order to control $N_3$, we need to assume $s_0>\frac{1}{2}$ and that $\| u_0\|_{H^{s_0}}$ be sufficiently small.
We also note that the dependence of $\eta ,T$ on $s_0,s$ comes from this part, since the constants $C_1$, $C_{2,\sigma}$, $\varepsilon$ in Proposition~\ref{prop:estN3} depends on $s_0$ and $C_{2,\sigma}$ also on $s$ when choosing $\sigma =s$.

All the required estimates are deduced from Proposition~\ref{prop:estN3} with $a_j=b_j=\frac{s_0}{2}$ ($j=1,2,3$).
For $u\in \mathscr{B}^s_{u_0}(T)$ and $\sigma =s,s_0,\frac{s_0}{2},0$, this (together with interpolation) implies that
\begin{align*}
\Big\| \int _0^t U_\mu (t-t')N_3(t')\,dt' \Big\| _{Z^\sigma _\mu (T)}&\leq \frac{C+C_\sigma T^\varepsilon}{\| u_0\| _{L^2}}\| u\| _{Z^{\frac{s_0}{2}}_0(T)}^2\| u\| _{Z^\sigma _\mu (T)}\\
&\leq \frac{C+C_\sigma T^\varepsilon}{\| u_0\| _{L^2}}\| u_0\| _{H^{\frac{s_0}{2}}}^2\| u_0\| _{H^\sigma}\\
&\leq \big( C+C_\sigma T^\varepsilon \big) \| u_0\| _{H^{s_0}}\| u_0\| _{H^\sigma}.
\end{align*}
(Note again that $C,C_\sigma ,\varepsilon$ depend on $s_0$.)
This is sufficient if we impose the condition
\begin{equation}\label{cond:N3}
\big( C+C_\sigma T^\varepsilon \big) \| u_0\| _{H^{s_0}}\leq \frac{1}{10}\min \{ C_0,\delta ,1\} .
\end{equation}

For the difference, noticing that
\[ N_3-\tilde{N}_3=N_3[u-\tilde{u},u,u]+N_3[\tilde{u},u-\tilde{u},u]+N_3[\tilde{u},\tilde{u},u-\tilde{u}] ,\]
we deduce from Proposition~\ref{prop:estN3} and interpolation that
\begin{align*}
&\Big\| \int _0^t U_\mu (t-t')\Big[ N_3-\tilde{N}_3\Big] (t')\,dt' \Big\| _{Z^\sigma _\mu (T)}\\
&\quad \leq \frac{C+C_\sigma T^\varepsilon}{\| u_0\| _{L^2}}\Big( \| u,\tilde{u}\| _{Z^{\frac{s_0}{2}}_0(T)}^2\| u-\tilde{u}\| _{Z^\sigma _\mu (T)}\\
&\hspace{100pt} +\| u,\tilde{u}\| _{Z^\sigma _\mu (T)}\| u,\tilde{u}\| _{Z^{\frac{s_0}{2}}_0(T)}\| u-\tilde{u}\| _{Z^{\frac{s_0}{2}}_0(T)}\Big) \\
&\quad \leq \big( C+C_\sigma T^\varepsilon \big) \Big( \| u_0\| _{H^{s_0}}\| u-\tilde{u}\| _{Z^\sigma _\mu (T)}+\frac{\| u_0\| _{H^\sigma}\| u_0\| _{H^{\frac{s_0}{2}}}}{\| u_0\| _{L^2}}\| u-\tilde{u}\| _{Z^{\frac{s_0}{2}}_0(T)}\Big) ,\\
&\hspace{300pt} \sigma =s,s_0,\frac{s_0}{2},0.
\end{align*}
The first term on the right-hand side is favorable if \eqref{cond:N3} holds.
The second term is evaluated as
\begin{align*}
\frac{\| u_0\| _{H^\sigma}\| u_0\| _{H^{\frac{s_0}{2}}}}{\| u_0\| _{L^2}}\| u-\tilde{u}\| _{Z^{\frac{s_0}{2}}_0(T)}
&\leq \frac{\| u_0\| _{H^\sigma}\| u_0\| _{H^{\frac{s_0}{2}}}^2}{\| u_0\| _{L^2}\| u_0\| _{H^s}}d(u,\tilde{u})\\
&\leq \| u_0\| _{H^{s_0}}\cdot \frac{\| u_0\| _{H^\sigma}}{\| u_0\| _{H^s}}d(u,\tilde{u}),
\end{align*}
which is also sufficient for our purpose under the condition \eqref{cond:N3}.

\bigskip
The claim follows from the estimates above.
To ensure the conditions \eqref{cond:N1}, \eqref{cond:N2} and \eqref{cond:N3}, we first take $\delta$ small (independent of $s_0,s$), then $\eta =\eta (s_0,\delta )$ small, and $T=T(s_0,s,\varepsilon )$ small. \hfill $\Box$

\bigskip
\noindent
{\bf Proof of Theorem~\ref{local} except (iii)}.
The claim (i) on existence has been shown in Proposition~\ref{prop:contraction}.
Let us show (ii), (iv) and (v) first on the maximal time interval where the solutions under consideration (starting from $u_0\neq 0$) do not become zero.
Namely, for $u_0\in H^{s}\setminus \{ 0\}$ with $\| u_0\|_{H^{s}}\leq \eta$, define
\[ T_*=T_*(u_0):=\begin{cases}
\inf \{ t\in (0,1]\,|\,S(u_0)(t)=0\} &\text{(if $\exists t\in (0,1];\ S(u_0)(t)=0$)}, \\
1 &\text{(otherwise)},\end{cases}
\]
and show (ii), (iv), (v) on $[0,T_*]$.
Note that $T_*>0$ by continuity.

To show (ii), let $u:=S(u_0)\in \mathscr{B}^{s}_{u_0}(1)$, and $\tilde{u}$ be another solution on $[0,T]$ starting from $u_0$ and belonging to $Z^{s}_\mu (T)$.
Noting that $\| u(t)\|_{H^{s}}\leq 2C_0\| u_0\|_{H^{s}}$ on $[0,1]$, we replace the threshold size $\eta$ by $(2C_0)^{-1}\eta$, so that $u(t)$ satisfies the initial data smallness condition in Proposition~\ref{prop:contraction} for any $t\in [0,1]$.
Now, suppose that $u(t_0)=\tilde{u}(t_0)$ for some $0\leq t_0<\min \{ T_*,T\}$, and consider
\begin{align*}
v(t,x)&:=u\Big( t+t_0,\,x +\frac{\alpha}{\pi}\big( \| u_0\| _{L^2}^2-\| u(t_0)\|_{L^2}^2\big) (t+t_0)\Big) \quad (t\in [0,1-t_0]),\\
\tilde{v}(t,x)&:=\tilde{u}\Big( t+t_0,\,x +\frac{\alpha}{\pi}\big( \| u_0\| _{L^2}^2-\| u(t_0)\|_{L^2}^2\big) (t+t_0)\Big) \quad (t\in [0,T-t_0]).
\end{align*}
Then, $v$ and $\tilde{v}$ are solutions to the renormalized equation \eqref{r-kdnls} with the same initial data $v_0(x):=u(t_0,x+\frac{\alpha}{\pi}(\| u_0\| _{L^2}^2-\| u(t_0)\|_{L^2}^2)t_0)$.
Moreover, since $\mu (t_0):=\frac{|\beta |}{2\pi}\|v_0\|_{L^2}^2\neq 0$, we have $v\in Z^{s}_{\mu (t_0)}(1-t_0)$ and $\tilde{v}\in Z^{s}_{\mu (t_0)}(T-t_0)$.
(Recall that for any $\mu >0$ the space$Z^s_\mu$ coincides with $Z^s_0\cap L^2_tH^{s+\frac12}_x$ and is invariant under spatial translations of the form $x\mapsto x-Ct$.)
We see from Lemma~\ref{lem:vanishing} that $\| v-U_{\mu (t_0)}(t)v_0\|_{Z^{s}_{\mu (t_0)}(T')}\to 0$ as $T'\to 0$ and the same is true for $\tilde{v}$.
In view of Lemma~\ref{lem:linear1}, both $v$ and $\tilde{v}$ belong to $\mathscr{B}^{s}_{u(t_0)}(T')$ for sufficiently small $T'>0$.
Therefore, $v(t)=\tilde{v}(t)$ on $[0,T']$, and thus $u(t)=\tilde{u}(t)$ on $[t_0,t_0+T']$.
Repeating this argument, we conclude that $u(t)=\tilde{u}(t)$ on $[0,\min \{ T_*,T\} ]$. 

Persistence of regularity (iv) is verified by a similar argument.
Assume $u_0\in H^{\tilde{s}}$ for some $\tilde{s}>s$.
From Proposition~\ref{prop:contraction}, there is some $T=T(\tilde{s},s)\in (0,1]$ such that we have a solution in $\mathscr{B}^{\tilde{s}}_{u_0}(T)\subset C([0,T];H^{\tilde{s}})$.
Since $\mathscr{B}^{\tilde{s}}_{u_0}(T)\subset \mathscr{B}^{s}_{u_0}(T)$, this solution coincides with the unique solution $S(u_0)$ in $\mathscr{B}^{s}_{u_0}(T)$, concluding that $S(u_0)\in C([0,T];H^{\tilde{s}})$. 
This argument can be repeated from any time $t_0\in (0,T_*)$ by considering the new solution $v$ with $v(0)=v_0$ defined above and noting that $v_0\in H^{\tilde{s}}$, $v_0\neq 0$, $\| v_0\| _{H^s}=\| u(t_0)\|_{H^s}\leq \eta$.

The smoothing property (v) follows from uniqueness and persistence of regularity, together with an observation that for any $s'\geq s$ and $u\in \mathscr{B}^{s'}_{u_0}(T)$, $u(t)$ has higher regularity $H^{s'+\frac12}$ for almost every $t\in (0,T)$.

So far, we have proved (i), and (ii), (iv), (v) on the interval $[0,T_*]$ for the renormalized equation \eqref{r-kdnls}.
Since the spatial translation $\mathcal{T}_{-\alpha}$ does not change the regularity and maps $Z^s_\mu (T)=Z^s_0(T)\cap L^2((0,T);H^{s+\frac12})$ into itself (unless $\mu =0$), we have the same conclusion for the original equation \eqref{kdnls}.
Finally, the $L^2$ lower bound for smooth solutions given in Proposition~\ref{prop:apriori} (ii) ensures that the solution $S(u_0)(t)$ starting from $u_0\neq 0$, which immediately gets smooth, will stay away from zero as long as it remains smooth. 
Therefore, we have $T_*(u_0)=1$.
\hfill $\Box$



\section{Continuity of the solution map}
\label{section:continuity}

%
In this section, we shall give a proof of Proposition~\ref{prop:cont}.

Let $\alpha \in \mathbf{R}$, $\beta <0$ be fixed, and let $s>1/2$.
In the preceding section, we have defined the solution map $S:u_0\mapsto u$ for the renormalized equation \eqref{r-kdnls} from $\{ u_0\in H^s\,|\, \| u_0\|_{H^s}\leq \eta \}$ into $C([0,1];H^s(\mathbf{T}))\cap L^2((0,1);H^{s+\frac12}(\mathbf{T}))$.
Notice that continuity of the map $S$ is not clear at the moment, because both the resolution space $Z^s_\mu (1)$ and the metric $d(\cdot ,\cdot )$ depend on the initial data.

\subsection{Lipschitz continuity away from the origin}

Let us show the following result on Lipschitz continuity of the map $S$ away from the origin.
\begin{prop}\label{prop:continuity}
Let $s>1/2$.
Then, replacing $\eta$ with a smaller one if necessary, we have
%
%
\begin{align*}
&\| S(u_0)-S(v_0)\| _{Z_{\max \{ \mu ,\nu \}}^s(1)}+M(u_0,v_0)^{\frac{1}{2}}\| S(u_0)-S(v_0)\| _{Z_{\max \{ \mu ,\nu \}}^{\frac{s}{2}}(1)}\\
&\qquad +\frac{M(u_0,v_0)}{\gamma}\| S(u_0)-S(v_0)\| _{Z_{\max \{ \mu ,\nu \}}^0(1)}\\
&\quad \lesssim \| u_0-v_0\| _{H^s}+\frac{M(u_0,v_0)}{\gamma}\| u_0-v_0\| _{L^2}
\end{align*}
for any $u_0,v_0\in \{ \phi \in H^s(\mathbf{T})\setminus \{ 0\}\,|\, \| \phi \|_{H^s}\leq \eta \}$, where $\mu :=\frac{|\beta |}{2\pi}\| u_0\| _{L^2}^2$, $\nu :=\frac{|\beta |}{2\pi}\| v_0\|_{L^2}^2$, $\gamma \in (0,1)$ is a small constant, and
\[ M(u_0,v_0):=\max \Big\{ \frac{\| u_0\| _{H^s}}{\| u_0\| _{L^2}},\;\frac{\| v_0\| _{H^s}}{\| v_0\| _{L^2}}\Big\} .\]
%
\end{prop}

\begin{rem}
(i) The constant $\gamma$ plays a role similar to that of $\delta$ in the contraction argument.
It is used to control a certain quantity related to the nonlinear term $N_1$ (see the estimate for \eqref{diff2} in the proof below).

(ii) Since $\| u\| _{C([0,1];H^s)} + \mu ^{\frac12} \| u\| _{L^2((0,1);H^{s+\frac12})}\lesssim \| u\| _{Z_{\mu}^s(1)}$, 
the estimate in Proposition~\ref{prop:continuity} is enough to verify the claim (i) in Proposition~\ref{prop:cont}.
In particular, for any $u_0\neq 0$ with $\| u_0\|_{H^s}\leq \eta$, there is an $H^s$-neighborhood of $u_0$ on which $S$ is Lipschitz continuous in this sense.

(iii) The claim (ii) in Proposition~\ref{prop:cont}, continuity of the solution map in the $H^s$-$C_tH^s_x$ topologies at the origin, follows from the bound $\| S(u_0)\|_{C([0,1];H^s)}\leq \| S(u_0)\| _{Z_\mu ^s(1)}\le 2C_0\| u_0\| _{H^s}$ obtained in the existence result.
It is not clear, however, whether $S$ is also continuous as a map into $L^2((0,1);H^{s+\frac12})$ at the origin.
We observe that this is false for the map $u_0\mapsto U_{\mu (u_0)}(t)u_0$ with $\mu (u_0) =\frac{|\beta |}{2\pi}\| u_0\|_{L^2}^2$.
In fact, fix $u_0\neq 0$ and consider the initial data $\varepsilon u_0$, $0< \varepsilon \ll 1$.
Then, $\mu (\varepsilon u_0)=\varepsilon ^2\mu (u_0)$, and we have
\begin{align*}
\big\| U_{\mu (\varepsilon u_0)}(t)(\varepsilon u_0)\big\| _{L^2((0,1);H^{s+\frac12})}
&=\Big\| e^{-\frac{|\beta|}{2\pi}\varepsilon ^2\| u_0\|_{L^2}^2|k|t}\langle k\rangle ^{s+\frac12}\varepsilon \hat{u}_0(k)\Big\| _{\ell ^2(\mathbf{Z};L^2(0,1))}\\
&\sim \| u_0\|_{L^2}^{-1}\Big\| \big( 1-e^{-\frac{|\beta |}{\pi}\varepsilon ^2\| u_0\|_{L^2}^2|k|}\big) ^{1/2}\langle k\rangle ^s\hat{u}_0(k)\Big\|_{\ell ^2} \\
&\to C\frac{\| u_0\|_{H^s}}{\| u_0\|_{L^2}}\ \neq 0\qquad (\varepsilon \to 0).
\end{align*}

(iv) Even in the $H^s$-$C_tH^s_x$ topologies, the proposition says nothing about Lipschitz or uniform continuity of the solution map on a neighborhood of the origin.
Indeed, we have claimed in Proposition~\ref{prop:cont} (iii) and will prove in Section~\ref{section:notUCD} that the solution map is not uniformly continuous around the origin.
\end{rem}

\medskip
\noindent
\textbf{Proof of Proposition~\ref{prop:continuity}}.
Without loss of generality, we assume $\mu \geq \nu$.
Note that the $Z_\mu ^s(1)$ norm is increasing in $\mu$.


Let $u:=S(u_0)$ and $v:=S(v_0)$.
We have
\begin{align*}
\partial _tu&=i\partial _x^2u-\mu D_xu+N_1[u,u_0,u]+N_2[u,u,u]+N_3[u,u,u], \\
\partial _tv&=i\partial _x^2v-\nu D_xv+N_1[v,v_0;v]+N_2[v,v,v]+N_3[v,v,v]\\
&=i\partial _x^2v-\mu D_xv+N_1[v,u_0;v] \\
&\quad +\frac{\alpha}{\pi}\big( \| u_0\|_{L^2}^2-\| v_0\|_{L^2}^2\big) \partial_xv+N_2[v,v,v]+N_3[v,v,v],
\end{align*}
where $N_1,N_2,N_3$ are the nonlinear terms defined in \eqref{eq:nonlinearity}.
Then, $w:=u-v$ satisfies
\begin{align*}
\partial _tw&=i\partial _x^2w-\mu D_xw+\Big\{ N_1[u,u_0;w]+N_1[u,v;v]-\frac{\alpha}{\pi}\big( \| u_0\|_{L^2}^2-\| v_0\|_{L^2}^2\big) \partial_xv \\
&\qquad\qquad\qquad + (N_2[u,u,u]-N_2[v,v,v]) + (N_3[u,u,u]-N_3[v,v,v]) \Big\} ,\\
w(0&)=u_0-v_0,
\end{align*}
so that
\begin{align}
w(t)&=U_\mu (t)(u_0 -v_0 ) \label{diff0}\\
&\quad +\int _0^t U_\mu (t-t')N_1[u(t'),u_0;w(t')] \,dt' \label{diff1}\\
&\quad +\int _0^t U_\mu (t-t')N_1[u(t'),v(t');v(t')]\,dt' \label{diff2}\\
&\quad -\frac{\alpha}{\pi} \int _0^t U_\mu (t-t')\Big\{ \big( \| u_0\|_{L^2}^2-\| v_0\|_{L^2}^2\big) \partial_xv(t')\Big\} \,dt' \label{diff2'}\\
&\quad +\int _0^t U_\mu (t-t')\Big\{ N_2[w,u,u]+N_2[v,w,u]+N_2[v,v,w]\Big\} \,dt' \label{diff3}\\
&\quad +\int _0^t U_\mu (t-t')\Big\{ N_3[w,u,u]+N_3[v,w,u]+N_3[v,v,w]\Big\} \,dt'. \label{diff4}
\end{align}

Recall that the solutions $u,v$ constructed in the previous section have the following bounds:
\begin{gather*}
\| u\| _{Z_\mu ^\sigma (1)}\leq 2C_0\| u_0 \|_{H^\sigma},\quad \| u-U_\mu (t)u_0\|_{Z_\mu ^0(1)}\leq \delta \| u_0\| _{L^2}, \\
\| v\| _{Z_\nu ^\sigma (1)}\leq 2C_0\| v_0\|_{H^\sigma},\quad \| v-U_\nu (t)v_0\|_{Z_\nu ^0(1)}\leq \delta \| v_0\| _{L^2}
\end{gather*}
for $\sigma =s,\frac{s}{2}$, where $\delta$ is the constant appearing in the definition of $\mathscr{B}^{s}_{u_0}$.
Since $\| u_0\|_{L^2}\geq \| v_0\|_{L^2}$, we have
\[ \left\langle i(\tau +k^2)+\| u_0\| _{L^2}^2|k|\right\rangle ^{\frac{1}{2}}\lesssim \frac{\| u_0\| _{L^2}}{\| v_0\| _{L^2}}\left\langle i(\tau +k^2)+\| v_0\| _{L^2}^2|k|\right\rangle ^{\frac{1}{2}},\]
which implies
\[ \| f\| _{X_\mu ^{\sigma ,\frac{1}{2}}}\lesssim \frac{\| u_0\| _{L^2}}{\| v_0\| _{L^2}} \| f\|_{X_\nu ^{\sigma ,\frac{1}{2}}},\qquad \sigma \in \mathbf{R}.\]
Therefore,
\begin{gather}
\| u\| _{Z_\mu ^s(1)}\lesssim \| u_0\|_{H^s} ,\quad \| v\| _{Z_\mu ^s(1)}\lesssim \frac{\| u_0\| _{L^2}}{\| v_0\| _{L^2}}\| v\| _{Z_\nu ^s(1)}\lesssim \frac{\| u_0\| _{L^2}}{\| v_0\| _{L^2}}\| v_0\|_{H^s},\label{u0-1}\\
\| u\| _{Z_\mu ^{\frac{s}{2}}(1)}\lesssim \| u_0\| _{H^{\frac{s}{2}}},\quad \| v\| _{Z_\mu ^{\frac{s}{2}}(1)}\lesssim \frac{\| u_0\| _{L^2}}{\| v_0\| _{L^2}}\| v\| _{Z_\nu ^{\frac{s}{2}}(1)}\lesssim \frac{\| u_0\| _{L^2}}{\| v_0\| _{L^2}}\| v_0\| _{H^{\frac{s}{2}}},\label{u0-2}\\
\| u\| _{Z_\mu ^0(1)}\lesssim \| u_0\| _{L^2},\quad \| v\| _{Z_\mu ^0(1)}\lesssim \frac{\| u_0\| _{L^2}}{\| v_0\| _{L^2}}\| v\| _{Z_\nu ^0(1)}\lesssim \| u_0\| _{L^2}.\label{u0-3}
\end{gather}
In the following, we assume $\| u_0,v_0\| _{H^s}\leq \eta _1$ with a small constant $\eta _1>0$ to be chosen later.

\bigskip
\noindent
\textbf{Estimate on \eqref{diff0}}.
By Lemma~\ref{lem:linear1}, we have 
\[ \| \eqref{diff0}\| _{Z^\sigma _\mu (1)}\lesssim \| u_0-v_0\| _{H^\sigma},\qquad \sigma \in \mathbf{R}.\]

\bigskip
\noindent
\textbf{Estimate on \eqref{diff1}}.
By Proposition~\ref{prop:estN1} (ii), we have
\[ \| \eqref{diff1} \|_{Z^\sigma _\mu (1)}\leq C\Big( \delta +\| u_0\| _{H^{\frac{1}{2}}}\Big) \| w\| _{Z^\sigma _\mu (1)}\leq \frac{1}{10}\| w\| _{Z^\sigma _\mu (1)},\qquad \sigma \geq 0.\]

\bigskip
\noindent
\textbf{Estimate on \eqref{diff2}}.
By Proposition~\ref{prop:estN1} (i), we have
\begin{align*}
\| \eqref{diff2}\| _{Z^\sigma _\mu (1)}&\lesssim \min \big\{ \| u_0\| _{L^2}^{-1}\| v\| _{Z^{\sigma}_\mu (1)},\;\| v\|_{L^\infty ((0,1);H^{\sigma +\frac{1}{2}})}\big\} \| w\| _{Z^0_0(1)},\qquad \sigma \geq 0.
\end{align*}
In particular, using \eqref{u0-1} and \eqref{u0-2}, we have
\begin{align*}
\| \eqref{diff2}\| _{Z^\sigma _\mu (1)}&\lesssim \frac{\| v_0\| _{H^\sigma}}{\| v_0\| _{L^2}}\| w\| _{Z^0_0(1)}\lesssim \gamma \cdot \frac{M(u_0,v_0)^{\frac{\sigma}{s}}}{\gamma}\| w\| _{Z^0_0(1)},\qquad \sigma =s,\frac{s}{2},
\end{align*}
and also
\begin{align*}
\| \eqref{diff2}\| _{Z^0_\mu (1)}&\lesssim \| v\| _{L^\infty ((0,1);H^{\frac{1}{2}})}\| w\| _{Z^0_0(1)}\lesssim \eta _1\| w\| _{Z^0_0(1)}.
\end{align*}

\bigskip
\noindent
\textbf{Estimate on \eqref{diff2'}}.
This term can be treated in the same manner as the preceding case \eqref{diff2}.
Proposition~\ref{prop:estN1} (i) yields
\begin{gather*}
\| \eqref{diff2'}\| _{Z^\sigma _\mu (1)}\lesssim \gamma \cdot \frac{M(u_0,v_0)^{\frac{\sigma}{s}}}{\gamma}\| u_0-v_0\| _{L^2},\qquad \sigma =s,\frac{s}{2},\\ 
\| \eqref{diff2'}\| _{Z^0_\mu (1)}\lesssim \eta _1\| u_0-v_0\| _{L^2}.
\end{gather*}

\bigskip
\noindent
\textbf{Estimate on \eqref{diff3}}.
We follow the proof for the estimate of $N_2-\tilde{N}_2$ in the preceding section.
In view of Proposition~\ref{prop:estN2}, we have
\begin{align*}
\| \eqref{diff3}\| _{Z^\sigma _\mu (1)}&\lesssim \| u_0\| _{L^2}^{-1}\| u,v\| _{L^\infty ((0,1);L^2)}\| u,v\|_{L^\infty ((0,1);H^{\frac{1}{2}})}\| w\| _{L^\infty ((0,1);H^\sigma )}\\
&\lesssim \eta _1\| w\| _{Z^\sigma _0(1)},\qquad \sigma \geq 0.
\end{align*}

\bigskip
\noindent
\textbf{Estimate on \eqref{diff4}}.
From Proposition~\ref{prop:estN3}, we have
\begin{align*}
\| \eqref{diff4}\| _{Z^\sigma _\mu (1)}
&\lesssim \| u_0\| _{L^2}^{-1}\Big( \| w\| _{Z^\sigma _\mu (1)}\| u,v\| _{Z^{\frac{s}{2}}_0(1)}^2+\| u,v\| _{Z^\sigma _\mu (1)}\| u,v\| _{Z^{\frac{s}{2}}_0(1)}\| w\| _{Z^{\frac{s}{2}}_0(1)}\Big) \\
&\lesssim \| u_0\| _{L^2}^{-1}\Big( \| w\| _{Z^\sigma _\mu (1)}\| u_0,v_0\| _{H^{\frac{s}{2}}}^2+\| u,v\| _{Z^\sigma _\mu (1)}\| u_0,v_0\| _{H^{\frac{s}{2}}}\| w\| _{Z^{\frac{s}{2}}_0(1)}\Big) 
\end{align*}
for $\sigma =s,\frac{s}{2},0$.
Now, by interpolation we see that
\[ \| u_0,v_0\| _{H^{\frac{s}{2}}}^2\leq \eta _1\| u_0\| _{L^2},\]
while the estimates \eqref{u0-1}--\eqref{u0-3} yield that
\[ \| u,v\| _{Z^\sigma _\mu (1)}\lesssim \| u_0\| _{H^\sigma}+\frac{\| u_0\| _{L^2}}{\| v_0\|_{L^2}}\| v_0\| _{H^\sigma}\lesssim \| u_0\| _{L^2}M(u_0,v_0)^{\frac{\sigma}{s}},\qquad \sigma =s,\frac{s}{2},0.\]
Hence, we have
\[ \| \eqref{diff4}\| _{Z^\sigma _\mu (1)}\lesssim \eta _1\| w\| _{Z^\sigma _\mu (1)}+M(u_0,v_0)^{\frac{\sigma}{s}}\| u_0,v_0\| _{H^{\frac{s}{2}}}\| w\| _{Z^{\frac{s}{2}}_0(1)},\qquad \sigma =s,\frac{s}{2},0.\]
To handle the second term on the right-hand side, we observe that
\begin{align*}
&\| u_0,v_0\| _{H^{\frac{s}{2}}}\| w\| _{Z^{\frac{s}{2}}_0(1)}\\
&\lesssim \min \{ \| u_0\| _{H^{\frac{s}{2}}},\,\| v_0\| _{H^{\frac{s}{2}}}\} \| w\| _{Z^{\frac{s}{2}}_0(1)}+\| u_0-v_0\| _{H^{\frac{s}{2}}}\| u,v\| _{Z^{\frac{s}{2}}_0(1)}\\
&\lesssim \min \{ \| u_0\| _{H^{\frac{s}{2}}},\,\| v_0\| _{H^{\frac{s}{2}}}\} \| w\| _{Z^{\frac{s}{2}}_0(1)}+\| u_0-v_0\| _{H^{\frac{s}{2}}}\| u_0,v_0\| _{H^{\frac{s}{2}}}\\
&\lesssim \min \{ \| u_0\| _{H^{\frac{s}{2}}},\,\| v_0\| _{H^{\frac{s}{2}}}\} \Big( \| w\| _{Z^{\frac{s}{2}}_0(1)}+\| u_0-v_0\| _{H^{\frac{s}{2}}}\Big) +\| u_0-v_0\| _{H^{\frac{s}{2}}}^2\\
&\lesssim \eta_1\min \Big\{ \frac{\| u_0\| _{L^2}}{\| u_0\| _{H^s}},\,\frac{\| v_0\| _{L^2}}{\| v_0\| _{H^s}}\Big\} ^{\frac{1}{2}}\Big( \| w\| _{Z^{\frac{s}{2}}_0(1)}+\| u_0-v_0\| _{H^{\frac{s}{2}}}\Big) +\eta_1\| u_0-v_0\| _{L^2}\\
&\lesssim \eta_1M(u_0,v_0)^{-\frac{1}{2}}\Big( \| w\| _{Z^{\frac{s}{2}}_0(1)}+\| u_0-v_0\| _{H^{\frac{s}{2}}}\Big) +\eta _1\| u_0-v_0\| _{L^2}.
\end{align*}
Therefore, we have
\begin{align*}
\| \eqref{diff4}\| _{Z^s_\mu (1)}&\lesssim \eta _1\| w\| _{Z^s_\mu (1)}+\eta _1M(u_0,v_0)^{\frac{1}{2}}\Big( \| w\| _{Z^{\frac{s}{2}}_0(1)}+\| u_0-v_0\| _{H^{\frac{s}{2}}}\Big) \\
&\quad +\eta _1M(u_0,v_0)\| u_0-v_0\| _{L^2},\\
\| \eqref{diff4}\| _{Z^{\frac{s}{2}}_\mu (1)}&\lesssim \eta _1\Big( \| w\| _{Z^{\frac{s}{2}}_\mu (1)}+\| u_0-v_0\| _{H^{\frac{s}{2}}}\Big) \\
&\quad +\eta _1M(u_0,v_0)^{\frac{1}{2}}\| u_0-v_0\| _{L^2},\\
\| \eqref{diff4}\| _{Z^0_\mu (1)}&\lesssim \frac{\eta _1}{\gamma}\cdot \gamma M(u_0,v_0)^{-\frac{1}{2}}\Big( \| w\| _{Z^{\frac{s}{2}}_0(1)}+\| u_0-v_0\| _{H^{\frac{s}{2}}}\Big) \\
&\quad +\eta _1\Big( \| w\| _{Z^0_\mu (1)}+\| u_0-v_0\| _{L^2}\Big) .
\end{align*}

\bigskip
For simplicity, we set
\begin{align*}
W&:=\| w\| _{Z^s_\mu (1)}+M(u_0,v_0)^{\frac{1}{2}}\| w\| _{Z^{\frac{s}{2}}_\mu (1)}+\frac{M(u_0,v_0)}{\gamma}\| w\| _{Z^0_\mu (1)},\\
W_0&:=\| u_0-v_0\| _{H^s}+M(u_0,v_0)^{\frac{1}{2}}\| u_0-v_0\| _{H^{\frac{s}{2}}}+\frac{M(u_0,v_0)}{\gamma}\| u_0-v_0\| _{L^2}.
\end{align*}
Then, combining the estimates obtained so far, we have
\[ W\leq \Big( \frac{1}{10}+C\gamma +C\frac{\eta _1}{\gamma}\Big) W+\Big( C+C\frac{\eta_1}{\gamma}\Big) W_0.\]
Taking $\gamma$ and then $\eta _1$ sufficiently small verifies $W\lesssim W_0$.
Finally, we see that
\[ M(u_0,v_0)^{\frac{1}{2}}\| u_0-v_0\| _{H^{\frac{s}{2}}}\leq \| u_0-v_0\| _{H^s}+M(u_0,v_0)\| u_0-v_0\| _{L^2},\]
which completes the proof.
\hfill $\Box$


\subsection{Failure of uniform continuity}
\label{section:notUCD}

In this subsection, we discuss uniformly continuous dependence of the solutions to \eqref{kdnls} or \eqref{r-kdnls} 
on initial data around the origin.
More precisely, we want to know whether the following property holds or not:
\begin{equation}\label{UCD}\tag{\text{UCD}}
\left\{ \begin{aligned}
~&\text{There exist $r>0$, $T>0$ and a continuous non-decreasing }\\
&\text{function $\varrho :[0,\infty )\to [0,\infty )$ with $\varrho (0)=0$ such that }\\
&\text{if $u_0,\tilde{u}_0\in H^s(\mathbf{T})$, $\| u_0\| _{H^s}\leq r$ and $\| \tilde{u}_0\| _{H^s}\leq r$, }\\
&\text{then the (well-defined) solutions $u,\tilde{u}\in C([0,T];H^s(\mathbf{T}))$ }\\
&\text{to the Cauchy problem with initial data $u_0,\tilde{u}_0$, respectively, }\\
&\text{satisfy $\| u-\tilde{u}\|_{C([0,T];H^s)}\leq \varrho (\| u_0-\tilde{u}_0\| _{H^s})$.}
\end{aligned}\right.
\end{equation}

First, we consider the following ``reduced'' KDNLS equation:
\begin{equation}\label{rkdnls}
\begin{cases}
\partial _tu=i\partial _x^2u+\dfrac{\alpha}{\pi}\| u(0)\| _{L^2}^2\partial_xu +\dfrac{\beta}{2\pi}\| u(0)\|_{L^2}^2D_xu,\quad &(t,x)\in [0,T]\times \mathbf{T},\\
u(0,x)=u_0(x), &x\in \mathbf{T}.
\end{cases}
\end{equation}
Since the reduced equation \eqref{rkdnls} is the ``effective linear part'' of the full equation \eqref{kdnls}, it would be reasonable to expect that \eqref{rkdnls} is a good approximation to \eqref{kdnls} (and also \eqref{rkdnls} with $\alpha =0$ is to \eqref{r-kdnls}).
We shall prove that the property \eqref{UCD} already fails for the reduced equation \eqref{rkdnls}.

%
\begin{prop}\label{prop:notUCD}
Let $\alpha \in \mathbf{R}$ and $\beta \leq 0$ be such that $(\alpha ,\beta )\neq (0,0)$.
Then, the property \eqref{UCD} does not hold for the Cauchy problem \eqref{rkdnls} in $H^s(\mathbf{T})$ if $s>0$.
\end{prop}

\begin{rem}
When $\beta <0$, the result is optimal in such a sense that the assumption on $s$ can not be relaxed.
In fact, the solution map for \eqref{rkdnls} is (globally) Lipschitz continuous as a map from $L^2(\mathbf{T})$ into $C([0,\infty );L^2(\mathbf{T}))$, as we see below.
Note first that the unique solution $u$ is given by
\[ \hat{u}(t,k)=e^{-ik^2t+\omega (k)\| u_0\|_{L^2}^2t}\hat{u}_0(k),\quad t\geq 0,\ k\in \mathbf{Z};\qquad \omega (k) :=i\frac{\alpha}{\pi}k+\frac{\beta}{2\pi}|k|.\]
For two initial data $u_0,v_0$, we may assume $\| u_0\|_{L^2}\geq \| v_0\|_{L^2}>0$, and write the corresponding solutions as $u,v$.
Then, for each $k\in \mathbf{Z}$ fixed, we have
\begin{align*}
\big| \hat{u}(t,k)-\hat{v}(t,k)\big| &\leq \big| e^{-ik^2t+\omega (k)\| u_0\|_{L^2}^2t}\big( \hat{u}_0(k)-\hat{v}_0(k)\big) \big| \\
&\quad +\big| e^{-ik^2t}\hat{v}_0(k) \big( e^{\omega (k)\| u_0\|_{L^2}^2t}-e^{\omega (k)\| v_0\|_{L^2}^2t}\big) \big| .
\end{align*}
For the second term on the right side, we estimate the difference of the exponentials as 
\begin{align*}
&\big| e^{\omega (k)\| u_0\|_{L^2}^2t}-e^{\omega (k)\| v_0\|_{L^2}^2t}\big| \\
&\quad = \Big| \omega (k)\big[ \| u_0\| _{L^2}^2-\| v_0\|_{L^2}^2\big] \int _0^t e^{\omega (k)\| u_0\|_{L^2}^2(t-t')+\omega (k)\| v_0\|_{L^2}^2t'}dt'\Big| \\
&\quad \lesssim |k| \big[ \| u_0\| _{L^2}^2-\| v_0\|_{L^2}^2\big] \int _0^te^{\frac{\beta}{2\pi}\| u_0\| _{L^2}^2|k|(t-t')}dt' \\
&\quad \lesssim \frac{(\| u_0\|_{L^2}+\| v_0\|_{L^2})\| u_0-v_0\|_{L^2}}{\| u_0\|_{L^2}^2}\lesssim \frac{\| u_0-v_0\| _{L^2}}{\| u_0\|_{L^2}}.
\end{align*}
(The implicit constant can be taken as $C\frac{|\alpha |+|\beta |}{|\beta |}$.)
Using this bound, we obtain
\[ \sup _{t\geq 0}\| u(t)-v(t)\| _{L^2}\lesssim \| u_0-v_0\|_{L^2}+\| v_0\| _{L^2}\frac{\| u_0-v_0\| _{L^2}}{\| u_0\|_{L^2}}\lesssim \| u_0-v_0\|_{L^2}.\]
We also remark that the (local) $L^2$-Lipschitz continuity of the solution map holds even for the full equation \eqref{kdnls} if it is restricted to smooth data $u_0\in H^s$, $s>3/2$.
This will be shown by the energy estimates (see the proof of Theorem~\ref{global} (ii) in Section~\ref{sec:proof2}).
\end{rem}

\noindent
\textbf{Proof of Proposition~\ref{prop:notUCD}}.
When $\beta \leq 0$, the Cauchy problem \eqref{rkdnls} has a unique solution $u\in C([0,\infty );H^s(\mathbf{T}))$ for any $u_0\in H^s(\mathbf{T})$ if $s\geq 0$.
In particular, for a positive integer $N$ and complex constants $a_0,b_0$,
the solution $u$ to \eqref{rkdnls} with initial data
\[ u_0(x)=\frac{1}{\sqrt{2\pi}}\big( a_0+b_0e^{iNx}\big) \]
is explicitly given by
\[ u(t,x)=\frac{1}{\sqrt{2\pi}}\big( a_0+b_0e^{-iN^2t+\omega (|a_0|^2+|b_0|^2)Nt}e^{iNx}\big) ,\qquad \omega :=i\frac{\alpha}{\pi}+\frac{\beta}{2\pi}\neq 0.\]
To prove the claim, it suffices to verify the following:
For any $r\in (0,1)$ there exist $\{ t_N\} _{N}\subset (0,1)$, $\{ u_{0,N}\} _N,\{ \tilde{u}_{0,N}\} _N\subset H^s(\mathbf{T})$ satisfying
\begin{gather*}
\| u_{0,N}\| _{H^s}\leq r,\quad \| \tilde{u}_{0,N}\| _{H^s}\leq r\qquad \text{for any $N$},\\
t_N\to 0,\quad \| u_{0,N}-\tilde{u}_{0,N}\|_{H^s}\to 0\qquad \text{as $N\to \infty$},
\end{gather*}
such that the corresponding solutions $u_N,\tilde{u}_N\in C([0,\infty );H^s(\mathbf{T}))$ satisfy
\[ \limsup _{N\to \infty} \| u_N(t_N)-\tilde{u}_N(t_N)\| _{H^s}\ >0.\]

Fix $0<\sigma < \min \{ s,\frac{1}{2}\}$, and set $t_N:=\frac{1}{100|\omega |}N^{-1+2\sigma}$.
For an arbitrary $r\in (0,1)$, define $u_{0,N},\tilde{u}_{0,N}$ as
\begin{gather*}
u_{0,N}(x):=\frac{1}{\sqrt{2\pi}}\big( a_{0,N}+b_{0,N}e^{iNx}\big) ,\qquad \tilde{u}_{0,N}(x):=\frac{1}{\sqrt{2\pi}}\big( \tilde{a}_{0,N}+\tilde{b}_{0,N}e^{iNx}\big) ,\\
a_{0,N}:=N^{-\sigma},\quad \tilde{a}_{0,N}:=2N^{-\sigma},\quad
b_{0,N}=\tilde{b}_{0,N}:=\frac{r}{2}N^{-s}.
\end{gather*}
It is easy to see that $\| u_{0,N}\| _{H^s}\leq r$, $\| \tilde{u}_{0,N}\| _{H^s}\leq r$ for all sufficiently large $N$, $t_N\to 0$ and $\| u_{0,N}-\tilde{u}_{0,N}\|_{H^s}=N^{-\sigma}\to 0$ as $N\to \infty$.
Moreover, we have
\begin{align*}
\| u_N(t_N)-\tilde{u}_N(t_N)\| _{H^s}
&\geq N^s|b_{0,N}|\Big| e^{\omega (|a_{0,N}|^2+|b_{0,N}|^2)Nt_N} - e^{\omega (|\tilde{a}_{0,N}|^2+|\tilde{b}_{0,N}|^2)Nt_N} \Big| \\
&= \frac{r}{2}\Big| e^{\frac{\omega}{100|\omega |}(|a_{0,N}|^2+|b_{0,N}|^2)N^{2\sigma}}-e^{\frac{\omega}{100|\omega |}(|\tilde{a}_{0,N}|^2+|\tilde{b}_{0,N}|^2)N^{2\sigma}}\Big| .
\end{align*}
Now, it is easily checked that if $z_1,z_2\in \mathbf{C}$ satisfy $|z_1|,|z_2|\leq \frac{1}{10}$, then $|e^{z_1}-e^{z_2}|\geq \frac{1}{2}|z_1-z_2|$.
From this and
\begin{gather*}
\big| (|a_{0,N}|^2+|b_{0,N}|^2)N^{2\sigma}\big| \leq 2,\qquad \big| (|\tilde{a}_{0,N}|^2+|\tilde{b}_{0,N}|^2)N^{2\sigma}\big| \leq 5,\\
\big| (|a_{0,N}|^2+|b_{0,N}|^2)N^{2\sigma}- (|\tilde{a}_{0,N}|^2+|\tilde{b}_{0,N}|^2)N^{2\sigma}\big| =3,
\end{gather*}
we deduce that
\[ \| u_N(t_N)-\tilde{u}_N(t_N)\| _{H^s} \geq \frac{r}{2}\cdot \frac{1}{2}\cdot \frac{3}{100} >0\]
for any $N$.
This completes the proof.
%
\hfill $\Box$

\bigskip
Non-uniformly continuous dependence for the full equation is obtained as a corollary, by which we complete the proof of Proposition~\ref{prop:cont} (iii).
\begin{cor}\label{cor:notUCD}
Let $\alpha \in \mathbf{R}$ and $\beta <0$.
Then, the property \eqref{UCD} does not hold, both for the Cauchy problem \eqref{r-kdnls}-\eqref{ic} and for the Cauchy problem \eqref{kdnls}-\eqref{ic}, in $H^s(\mathbf{T})$ for any $s>1/2$.
\end{cor}

\noindent
\textbf{Proof}.
Recall that for $s>1/2$ there exists $\eta =\eta (s)>0$ such that the solution maps $S$ and $\overline{S}:=\mathcal{T}_{-\alpha}\circ S$ for the Cauchy problems \eqref{r-kdnls}--\eqref{ic} and \eqref{kdnls}--\eqref{ic} have been uniquely defined on $\{ u_0\in H^s(\mathbf{T})\,|\,\| u_0\|_{H^s}\leq \eta \}$.

To show that the solution map $S$ for \eqref{r-kdnls} does not satisfy \eqref{UCD}, we fix an arbitrary $0<r\leq \eta$ and take the same sequences $\{ u_{0,N}\} _N,\{ \tilde{u}_{0,N}\} _N$ of initial data as in the proof of Proposition~\ref{prop:notUCD} (setting $\sigma =\frac14$, for instance).
Write the corresponding solutions to \eqref{rkdnls} with $\alpha =0$ as $u_{N,lin}$ and $\tilde{u}_{N,lin}$, respectively, and $u_N:=S(u_{0,N})$, $\tilde{u}_N:=S(\tilde{u}_{0,N})$.
We have already seen that 
\begin{gather*}
\| u_{0,N}\|_{H^s}\sim \| \tilde{u}_{0,N}\|_{H^s}\sim r,\qquad \| u_{N,lin}(t_N)-\tilde{u}_{N,lin}(t_N)\|_{H^s}\geq 3cr,\\
t_N\sim N^{-\frac{1}{2}}\to 0,\qquad \| u_{0,N}-\tilde{u}_{0,N}\|_{H^s}\to 0\qquad (N\to \infty ),
\end{gather*}
where $c=\frac{1}{400}$, and the implicit constants are independent of $N$, $r$.
Now that $u_{N,lin}(t)=U_{\mu _N}(t)u_{0,N}$ with $\mu _N:=\frac{|\beta|}{2\pi}\| u_{0,N}\|_{L^2}^2$, we have
\begin{align*}
u_N(t)-u_{N,lin}(t)&=\int _0^t U_{\mu _N}(t-t')N_1[u_N(t'),u_{0,N};u_N(t')] \,dt' \\
&\quad +\int _0^t U_{\mu _N}(t-t')N_2[u_N(t'),u_N(t'),u_N(t')]\,dt'\\
&\quad +\int _0^t U_{\mu _N}(t-t')N_3[u_N(t'),u_N(t'),u_N(t')]\,dt',
\end{align*}
and similarly for $\tilde{u}_N-\tilde{u}_{N,lin}$.
Each term on the right side has been estimated in Propositions~\ref{prop:estN1}(ii), \ref{prop:estN2}, and \ref{prop:estN3}, respectively.
Noticing that $u_N\in \mathscr{B}_{u_{0,N}}^s(1)$, we see
\begin{align*}
&\| u_N-u_{N,lin}\| _{C([0,1];H^s)}\leq \| u_N-u_{N,lin}\| _{Z^{s}_{\mu _N}(1)}\\
&\quad \leq C\big( \delta +\| u_{0,N}\|_{H^s}\big) \| u_{0,N}\|_{H^s}+C\| u_{0,N}\|_{H^s}^2+C_s\| u_{0,N}\|_{H^s}^2\\
&\quad \leq (C\delta + C_sr)r,
\end{align*}
where $\delta$ is the constant determined in the proof of Proposition~\ref{prop:contraction}.
Possibly after replacing $\delta$ with a smaller one, and recalling that $r\leq \eta$ and that $\eta$ has been chosen to satisfy $\eta \ll \delta$, we confirm that
\begin{equation}\label{est:nonlerror}
\| u_N-u_{N,lin}\| _{C([0,1];H^s)}\leq cr,\qquad \| \tilde{u}_N-\tilde{u}_{N,lin}\| _{C([0,1];H^s)}\leq cr,
\end{equation}
from which it follows that $\| u_{N}(t_N)-\tilde{u}_{N}(t_N)\|_{H^s}\geq cr$, and hence the failure of \eqref{UCD} for $S$.

Concerning the solution map $\overline{S}$ for \eqref{kdnls}, we note that the solution $\overline{S}(u_{0,N})$ to \eqref{kdnls} and the corresponding solution $u_{N,lin}$ to \eqref{rkdnls} can be written as
\begin{gather*}
\overline{S}(u_{0,N})(t,x)=S(u_{0,N})\Big( t,x+\frac{\alpha}{\pi}\| u_{0,N}\|_{L^2}^2t\Big) ,\\
u_{N,lin}(t,x)=\big[ U_{\mu _N}(t)u_{0,N}\big] \Big( x+\frac{\alpha}{\pi}\| u_{0,N}\|_{L^2}^2t\Big) .
\end{gather*}
Hence, the estimates from the proof of Proposition~\ref{prop:notUCD} and the above estimate \eqref{est:nonlerror} show that (with $\tilde{\mu}_N:=\frac{|\beta|}{2\pi}\| \tilde{u}_{0,N}\|_{L^2}^2$)
\begin{align*}
&\| \overline{S}(u_{0,N})(t_N)-\overline{S}(\tilde{u}_{0,N})(t_N)\|_{H^s}\\
&\quad \geq \| u_{N,lin} (t_N) - \tilde{u}_{N,lin}(t_N) \|_{H^s} \\
&\qquad - \| \overline{S}(u_{0,N}) - u_{N,lin}\|_{C([0,1];H^s)} - \| \overline{S}(\tilde{u}_{0,N}) - \tilde{u}_{N,lin}\|_{C([0,1];H^s)}\\
&\quad = \| u_{N,lin} (t_N) - \tilde{u}_{N,lin}(t_N) \|_{H^s} \\
&\qquad - \| S(u_{0,N}) - U_{\mu _N}u_{0,N}\|_{C([0,1];H^s)} - \| S(\tilde{u}_{0,N}) - U_{\tilde{\mu}_N}\tilde{u}_{0,N}\|_{C([0,1];H^s)}\\
&\quad \geq 3cr - cr - cr =cr.
\end{align*}
This concludes the failure of \eqref{UCD} for $\overline{S}$.
\hfill $\Box$


\section{Proof of Theorem~\ref{global}}
\label{sec:proof2}

In this section, we shall prove Theorem~\ref{global}.
Let $\alpha \in \mathbf{R}$ and $\beta <0$ be fixed.

\medskip
\noindent
{\bf Proof of Theorem~\ref{global} (i)}.
Let $s>1/2$ and $u \in C([0,T];H^s(\mathbf{T}))$ be a solution to the KDNLS equation \eqref{kdnls}.
Let $P_{\leq N} :=\mathcal{F}^{-1}_k\chi _{\{ |k|\leq N\}}\mathcal{F}_x$, $P_{>N}:=1-P_{\leq N}$, and write $u^N:=P_{\leq N}u$.
Note that $u^N\to u$ in $C([0,T];H^s)$ as $N\to \infty$, and that $u^N$ is a smooth solution of
\[ \partial_t u^N -i \partial_x^2 u^N = \alpha P_{\leq N}\partial_x \big( |u|^2u\big) + \beta P_{\leq N}\partial_x \big[ H(|u|^2)u\big] ,\quad (t,x)\in [0,T]\times \mathbf{T}.\]
Taking the $L^2$ inner product of the equation with $u^N$, we have
\begin{align*}
\frac{d}{dt}\big\| u^N(t)\big\|_{L^2}^2 &= 2\mathrm{Re}\int _{\mathbf{T}}\partial_x F_N(u)\cdot \overline{u^N}\,dx\\
&\quad +2\mathrm{Re}\int _{\mathbf{T}}\Big\{ \alpha \partial_x \big( \big| u^N\big| ^2u^N \big) +\beta \partial_x \big[ H\big( \big| u^N\big| ^2\big) u^N\big] \Big\} \cdot \overline{u^N}\,dx\\
&=:I_1(t)+I_2(t),
\end{align*}
where 
\[ F_N(u):=\alpha \Big( P_{\leq N}(|u|^2u)-\big| u^N \big|^2u^N\Big) +\beta \Big( P_{\leq N}\big( H(|u|^2)u\big) -H\big( \big| u^N\big| ^2\big) u^N\Big) .\]
Integrating both sides, we have
\begin{equation}\label{L2N}
\big\| u^N(t)\big\| _{L^2}^2\ =\ \big\| u^N(0)\big\|_{L^2}^2+\int_0^t\Big( I_1(\tau )+I_2(\tau )\Big) \,d\tau ,\qquad t\in [0,T].
\end{equation}
On one hand, by the Sobolev inequality and the boundedness of the Hilbert transform $H$ on $H^s$, we can easily see that
\begin{align*}
&\big\| F_N(u)\big\|_{H^{\frac12}}\\
&\quad \leq \sum _{A\in \{ \alpha I\!d,\beta H\}}\Big( \big\| P_{>N}\big( A(|u|^2)u\big) \big\| _{H^\frac12} +\big\| A(|u|^2)u-A\big( \big| u^N\big| ^2\big) u^N\big\|_{H^\frac12}\Big) \\
&\quad \lesssim N^{-(s-\frac12)}\| u\|_{H^s}^3.
\end{align*}
Therefore, we have
\[ \int _0^T \big| I_1(\tau )\big| \,d\tau \lesssim TN^{-(s-\frac12 )}\| u\|_{C([0,T];H^s)}^4\to 0\quad (N\to \infty ).\]
On the other hand, using integration by parts, 
\begin{align*}
I_2&=\int _{\mathbf{T}}\Big\{ \frac{3}{2}\alpha \partial_x \big( \big| u^N\big| ^4\big)  + 2\beta D_x\big( \big| u^N\big| ^2\big) \big| u^N\big| ^2 + \beta \big( \big| u^N\big| ^2\big) \partial_x\big( \big| u^N\big| ^2\big) \Big\} \,dx\\
&=\beta \int _{\mathbf{T}}D_x\big( \big| u^N\big| ^2\big) \big| u^N\big| ^2\,dx =\beta \big\| D_x^{\frac12}\big( \big| u^N\big| ^2\big) \big\|_{L^2}^2.
\end{align*} 
Again from the Sobolev inequality, $I_2(t)$ converges to $\beta \| D^{1/2}_x(|u(t)|^2)\|_{L^2}^2$ uniformly in $t\in [0,T]$ as $N\to \infty$.
Hence, taking the limit $N\to \infty$ on both sides of \eqref{L2N}, we obtain the claim.
\hfill $\Box$

\bigskip
\noindent
{\bf Proof of Theorem~\ref{global} (ii)}.
Let $s>3/2$ and suppose that $u_1,u_2\in C([0,T];H^s(\mathbf{T}))$ are solutions of \eqref{kdnls} starting from $u_j(0)=u_{0,j}$, $j=1,2$.
Similarly to the proof of (i), we define $u_j^N:=P_{\leq N}u_j$, and set $w^N:=u_1^N-u_2^N$.
Note that $w^N$ is a smooth solution of
\begin{align*}
&\partial_t w^N -i \partial_x^2 w^N \\
&\quad = \partial_x F_N(u_1) - \partial_x F_N(u_2) + \sum _{A\in \{ \alpha I\!d,\beta H\}}\partial_x \Big[ A\big( \big| u_1^N\big| ^2\big) u_1^N -A\big( \big| u_2^N\big| ^2\big) u_2^N \Big] ,
\end{align*}
and then
\begin{align*}
\frac{d}{dt}\big\| w^N(t)\big\|_{L^2}^2 &= 2\mathrm{Re}\int _{\mathbf{T}}\Big( \partial_x F_N(u_1)-\partial_x F_N(u_2)\Big) \overline{w^N}\,dx \\
&\quad +\sum _{A\in \{ \alpha I\!d,\beta H\}} 2\mathrm{Re}\int_{\mathbf{T}} \partial_x \Big[ A\big( \big| u_1^N\big| ^2\big) u_1^N -A\big( \big| u_2^N\big| ^2\big) u_2^N\Big] \cdot \overline{w^N}\,dx.
\end{align*}
A similar argument as before using the Sobolev inequality implies that
\[ \big\| \partial_x F_N(u_j) (t) \big\|_{L^2}\lesssim N^{-(s-1)}M^3,\quad j=1,2,\quad t\in [0,T],\]
where
\[ M:=\max _{j=1,2}\| u_j\|_{C([0,T];H^s)}<\infty.\]
To handle the remaining terms, we first note that
\begin{align*}
&A\big( \big| u_1^N\big| ^2\big) u_1^N -A\big( \big| u_2^N\big| ^2\big) u_2^N\\
&\quad =A\big( \big| u_1^N\big| ^2\big) w^N + A\big( \big| w^N\big| ^2\big) u_2^N + A\big( 2\mathrm{Re}(u_2^N\overline{w^N})\big) u_2^N
\end{align*}
for $A\in \{ \alpha I\!d,\,\beta H\}$.
Applying integration by parts and the Sobolev inequality, the contribution from the first two terms on the right side is estimated (whenever $s>3/2$) as follows.
\[ \Big| 2\mathrm{Re}\int_{\mathbf{T}} \partial_x \Big[ A\big( \big| u_1^N\big| ^2\big) w^N + A\big( \big| w^N\big| ^2\big) u_2^N\Big] \cdot \overline{w^N}\,dx\Big| \ \lesssim \ M^2\big\| w^N(t)\big\|_{L^2}^2.\]
(The second term can be treated with the identity $\int Hf\cdot \overline{g}\,dx = \int f\cdot \overline{Hg}\,dx$.)
For the last term, we see that
\begin{align*}
&2\mathrm{Re}\int_{\mathbf{T}} \partial_x \Big[ A\big( 2\mathrm{Re}(u_2^N\overline{w^N})\big) u_2^N\Big] \cdot \overline{w^N}\,dx\\
&\quad =\int_{\mathbf{T}} A\big( 2\mathrm{Re}(u_2^N\overline{w^N})\big) 2\mathrm{Re}(\partial_x u_2^N\cdot \overline{w^N}) \,dx \\
&\qquad \qquad + \int_{\mathbf{T}}\partial_x A\big( 2\mathrm{Re}(u_2^N\overline{w^N})\big) \cdot 2\mathrm{Re}(u_2^N\overline{w^N}) \,dx\\
&\quad =\int_{\mathbf{T}} A\big( 2\mathrm{Re}(u_2^N\overline{w^N})\big) 2\mathrm{Re}(\partial_x u_2^N\cdot \overline{w^N})\,dx \\
&\qquad \qquad + \begin{cases}
0 &(A=\alpha I\!d),\\
\beta \big\| D_x^{1/2}\big( 2\mathrm{Re}(u_2^N\overline{w^N})\big) \big\|_{L^2}^2 &(A=\beta H).
\end{cases}
\end{align*}
The first integral on the right side is estimated by $M^2\| w^N(t)\|_{L^2}^2$ as before.
We now invoke the assumption $\beta \leq 0$ and dispose of the second term on the right side.
As a result, we have
\begin{align*}
\frac{d}{dt}\big\| w^N(t)\big\|_{L^2}^2 &\lesssim N^{-(s-1)}M^3\big\| w^N(t)\big\|_{L^2} + M^2 \big\| w^N(t)\big\|_{L^2}^2\\
&\lesssim M^2 \Big( \big\| w^N(t)\big\|_{L^2}^2 + N^{-2(s-1)}M^2\Big) ,\qquad t\in [0,T].
\end{align*}
The Gronwall inequality then implies that
\[ \big\| w^N(t)\big\|_{L^2}^2 \leq \Big( \big\| w^N(0)\big\|_{L^2}^2 + N^{-2(s-1)}M^2\Big) e^{CM^2t},\qquad t\in [0,T]. \]
Letting $N\to \infty$, we conclude that
\[ \| u_1(t)-u_2(t) \|_{L^2}^2 \leq \| u_{0,1} - u_{0,2} \|_{L^2}^2e^{CM^2t},\qquad t\in [0,T].\]
In particular, we have $u_1(t)=u_2(t)$ on $[0,T]$ if $u_{0,1}=u_{0,2}$, as desired.
\hfill $\Box$

\bigskip
\noindent
{\bf Proof of Theorem~\ref{global} (iii)}.
To prove non-existence of backward solutions in the case $\beta <0$, suppose for contradiction that there exists a solution $u\in C([-T,0];H^s(\mathbf{T}))$ for some $T>0$.
By the assumption on $u_0$, there exists $0<t_0\leq \min \{ T ,1\}$ such that $u(-t_0)\in H^s(\mathbf{T})$ and $\| u(-t_0)\|_{H^{s'}}\leq \eta (s')$.
From Theorem~\ref{local}, there is a solution $v\in C([-t_0,0];H^s(\mathbf{T}))\cap C^\infty ((-t_0,0]\times \mathbf{T})$ on $[-t_0,0]$ starting from $v(-t_0)=u(-t_0)$.
Then, unconditional uniqueness (ii) we have just proved implies that $u_0=v(0)\in C^\infty (\mathbf{T})$, which contradicts the assumption on $u_0$.
Therefore, such a solution does not exist for any $T>0$.
The non-existence of forward solutions for $\beta >0$ is verified by the transform $u(t,x)\mapsto \overline{u(-t,x)}$, which converts \eqref{kdnls} to the same equation but $(\alpha ,\beta )$ replaced with $(-\alpha ,-\beta)$.
\hfill $\Box$

\bigskip
\noindent
{\bf Proof of Theorem~\ref{global} (iv)}.
We only consider the case $s\in (1/2,1)$, since the case $s\geq 1$ is easier.
Assume $u_0\in H^s(\mathbf{T})$ and $0\neq \| u_0\|_{H^s}\leq \tilde{\eta}$.
If $\tilde{\eta}\leq \eta (s)$, from Proposition~\ref{prop:contraction} there is a solution $u\in Z^s_\mu (1)$ satisfying
\begin{gather*}
 \| u\| _{C([0,1];H^s)}\leq C\| u_0\| _{H^s},\\
\| u\| _{L^2((0,1);H^{s+\frac12})}\leq C\frac{\| u\|_{Z^s_\mu (1)}}{\| u_0\|_{L^2}}\leq C\frac{\| u_0\|_{H^s}}{\| u_0\|_{L^2}}.
\end{gather*}
In particular, there exists $t_0\in (0,1]$ such that $\| u(t_0)\|_{H^{s+\frac12}}\leq C\| u_0\|_{L^2}^{-1}\| u_0\|_{H^s}$.
By interpolation,
\[ \| u(t_0)\| _{H^1}\leq \| u(t_0)\|_{H^s}^{2s-1}\| u(t_0)\|_{H^{s+\frac12}}^{2(1-s)}\leq C\frac{\| u_0\|_{H^s}}{\| u_0\|_{L^2}^{2(1-s)}}.\]
Since the solution is smooth for $t\geq t_0$, the $H^1$ a priori estimate from Proposition~\ref{prop:apriori} is available.
Also using the $L^2$ bound from (i), we have
\[ \| u(t)\| _{H^1}\leq Ce^{C\| u(t_0)\|_{L^2}^2}\| u(t_0)\| _{H^1}\leq Ce^{C\| u_0\|_{L^2}^2}\frac{\| u_0\|_{H^s}}{\| u_0\|_{L^2}^{2(1-s)}},\qquad t\geq t_0.\]
Further interpolating it with the $L^2$ bound, and noting that $1-s_0=2(1-s)s_0$ for $s_0=\frac{1}{3-2s}$, we have
\begin{align*}
\| u(t)\| _{H^{s_0}}&\leq \| u(t)\| _{L^2}^{1-s_0}\| u(t)\|_{H^1}^{s_0}\leq Ce^{C\| u_0\|_{L^2}^2}\| u_0\| _{L^2}^{1-s_0}\frac{\| u_0\|_{H^s}^{s_0}}{\| u_0\|_{L^2}^{2(1-s)s_0}}\\
&=Ce^{C\| u_0\|_{L^2}^2}\| u_0\|_{H^s}^{s_0}, \qquad t\geq t_0.
\end{align*}
Hence, if $\tilde{\eta}$ is chosen so small that
\[ Ce^{C\tilde{\eta}^2}\tilde{\eta}^{s_0}\leq \eta (s_0),\]
then we can iterate the construction of the local $H^{s_0}$ solution in Theorem~\ref{local} and obtain a global smooth solution.

This is the end of the proof of Theorem~\ref{global}.
\hfill $\Box$


\appendix

\section{Proof of the linear estimates}
\label{sec:appendix}

In the Appendix, we describe the proof of Lemmas \ref{lem:linear1} and \ref{lem:inhom} in a slightly more general setting to make it clear what influence the modification of $e^{it \partial_x^2 -\mu t D_x}$ to $e^{it \partial_x^2 -\mu |t| D_x}$ has on the linear estimates, which may be useful for the future work.
We avoid a sharp cut-off $\chi_{\mathbf{R}_+}(t)$, which is incompatible with $X_\mu^{s, 1/2}$
.
In fact, we give a negative result (Lemma~\ref{counterex} below) about the estimate of $\chi_{\mathbf{R}_+} f$ in $H^{1/2}$ related to the inequality mentioned in \cite[the remark about the case $b = 1/2$ on the last line on page 1989 to the first line on page 1990]{MR}, which is of interest on its own.
But it does not matter for the proof of the main theorem in \cite{MR}, because it is correctly reformulated in \cite{MV1} and \cite{MV2}.

Let $\mathcal{M}$ be either $\mathbf{R}^d$ or $\mathbf{T}^d$, $d\geq 1$, and denote the space of the corresponding Fourier variable, $\mathbf{R}^d$ or $\mathbf{Z}^d$, by $\hat{\mathcal{M}}$.
We consider the following general linear dispersive-parabolic equation
\[ \partial_t u - ia(D) u + p(D) u = 0,\qquad t>0,\quad x\in \mathcal{M},\]
where $a(\cdot )$ and $p(\cdot )$ are arbitrary functions on $\hat{\mathcal{M}}$ taking values in $\mathbf{R}$ and $[0,\infty )$, respectively, and $a(D)$, $p(D)$ denote the corresponding Fourier multipliers.
Define $\overline{U}(t)=e^{ita(D)-|t|p(D)}=\mathcal{F}^{-1}_\xi e^{ia(\xi )t-p(\xi )|t|}\mathcal{F}_x$, and define the Banach spaces $\overline{X}^{s,b}$, $\overline{Y}^{s,b}$ corresponding to the linear propagator $\overline{U}(t)$ by the restriction of tempered distributions $\mathscr{S}'(\mathbf{R}\times \mathcal{M})$ to those of which the following norms are finite:
\begin{align*}
\| u\| _{\overline{X}^{s,b}}&:=\Big\| \langle \xi \rangle ^s\left\langle i(\tau -a(\xi ))+p(\xi)\right\rangle ^b\tilde{u}(\tau ,\xi )\Big\| _{L^2_{\tau ,\xi}(\mathbf{R}\times \hat{\mathcal{M}})},\\
\| u\| _{\overline{Y}^{s,b}}&:=\Big\| \langle \xi \rangle ^s\big\| \left\langle i(\tau -a(\xi ))+p(\xi )\right\rangle ^b\tilde{u}(\tau ,\xi )\big\| _{L^1_\tau (\mathbf{R})} \Big\| _{L^2_\xi (\hat{\mathcal{M}})}.
\end{align*}

In the following estimates, $\psi$ is a fixed function in $C^\infty_0(\mathbf{R};[0,1])$ satisfying $\psi =1$ on $[-1,1]$ and $\mathrm{supp}\,\psi \subset [-2,2]$.

\begin{lem}[{cf.~\cite[Proposition~2.1]{MR}}]\label{lem:linear1-abs}
For any $s\in \mathbf{R}$, we have 
\begin{equation*}
\| \psi (t)\overline{U}(t)\phi \| _{\overline{X}^{s,\frac12}\cap \overline{Y}^{s,0}}\lesssim \| \phi \| _{H^s}, 
\end{equation*}
where the implicit constant is independent of $s$, $a(\cdot )$, and $p(\cdot )$.
\end{lem}

\bigskip
\noindent
\textbf{Proof}.
We follow the argument in the proof of \cite[Proposition~2.1]{MR}.
Observe that
\begin{align*}
&\left\| \psi (t)\overline{U}(t)\phi \right\|_{\overline{X}^{s,\frac{1}{2}}}\\
&\quad = \left\| \langle \xi \rangle ^s\langle i(\tau -a(\xi ))+p(\xi )\rangle ^{\frac{1}{2}}\mathcal{F} _t\left[ \psi (t)e^{-p(\xi )|t|}\right] (\tau -a(\xi ))\hat{\phi}(\xi ) \right\| _{L^2_{\tau,\xi}}\\
&\quad \lesssim \left\| \langle \xi \rangle ^s\hat{\phi}(\xi ) \left\| \langle \tau \rangle ^{\frac{1}{2}}\mathcal{F}_t\big[ \psi(t)e^{-p(\xi )|t|}\big] (\tau )\right\| _{L^2_\tau}\right\| _{L^2_\xi}\\
&\qquad +\left\| \langle \xi \rangle ^s\hat{\phi}(\xi )p(\xi )^{\frac{1}{2}} \left\| \mathcal{F}_t\big[ \psi(t)e^{-p(\xi )|t|}\big] (\tau )\right\| _{L^2_\tau}\right\| _{L^2_\xi},
\end{align*}
and similarly,
\[ \left\| \psi (t)\overline{U}(t)\phi \right\|_{\overline{Y}^{s,0}}=\left\| \langle \xi \rangle ^s\hat{\phi}(\xi ) \left\| \mathcal{F}_t\big[ \psi(t)e^{-p(\xi )|t|}\big] (\tau )\right\| _{L^1_\tau}\right\| _{L^2_\xi}.\]
It then suffices to show
\begin{align*}
\left\| \langle \tau \rangle ^{\frac{1}{2}}\mathcal{F}_t\big[ \psi(t)e^{-p(\xi )|t|}\big] (\tau )\right\| _{L^2_\tau} + p(\xi )^{\frac{1}{2}} \left\| \mathcal{F}_t\big[ \psi(t)e^{-p(\xi )|t|}\big] (\tau )\right\| _{L^2_\tau} &\\
+\left\| \mathcal{F}_t\big[ \psi(t)e^{-p(\xi )|t|}\big] (\tau )\right\| _{L^1_\tau} &\lesssim 1,
\end{align*}
uniformly in $\xi$ and $p(\cdot )$.
This estimate holds trivially when $p(\xi )=0$.
If $p(\xi )>0$, the scaling properties
\begin{gather*}
\left\| \mathcal{F}_t[e^{-p(\xi )|t|}]\right\| _{L^1_\tau}=\left\| \mathcal{F}_t[e^{-|t|}]\right\| _{L^1_\tau},\qquad
\left\| e^{-p(\xi )|t|}\right\| _{\dot H_t^{\frac{1}{2}}}=\left\| e^{-|t|}\right\| _{\dot H_t^{\frac{1}{2}}},\\
p(\xi )^{\frac{1}{2}}\left\| e^{-p(\xi )|t|}\right\| _{L^2_t}=\left\| e^{-|t|}\right\| _{L^2_t}
\end{gather*}
imply that
\begin{align*}
&\left\| \langle \tau \rangle ^{\frac{1}{2}}\mathcal{F}_t\big[ \psi(t)e^{-p(\xi )|t|}\big] (\tau )\right\| _{L^2_\tau}\\
&\qquad \lesssim \left\| \langle \tau \rangle ^{\frac{1}{2}}\hat{\psi}\right\| _{L^2_\tau}\left\| \mathcal{F}_t[e^{-p(\xi )|t|}]\right\| _{L^1_\tau} + \left\| \hat{\psi}\right\| _{L^1_\tau}\left\| |\tau |^{\frac{1}{2}}\mathcal{F}_t[e^{-p(\xi )|t|}]\right\| _{L^2_\tau} \lesssim 1,\\[5pt]
&p(\xi )^{\frac{1}{2}} \left\| \mathcal{F}_t\big[ \psi(t)e^{-p(\xi )|t|}\big] (\tau )\right\| _{L^2_\tau}\leq p(\xi )^{\frac{1}{2}} \left\| \psi \right\|_{L^\infty_t}\left\| e^{-p(\xi )|t|}\right\|_{L^2_t} \lesssim 1,\\
&\left\| \mathcal{F}_t\big[ \psi(t)e^{-p(\xi )|t|}\big] (\tau )\right\| _{L^1_\tau}\lesssim \left\| \hat{\psi}\right\| _{L^1_\tau}\left\| \mathcal{F}_t[e^{-p(\xi )|t|}]\right\| _{L^1_\tau} \lesssim 1.
\end{align*}
Therefore, we obtain the claimed estimate.
\hfill $\Box$

\begin{lem}[{cf.~\cite[Proposition~2.3(a)]{MR}}]\label{lem:inhom-abs}
For any $s\in \mathbf{R}$, we have
\begin{equation*}
\Big\| \psi (t)\int _0^t \overline{U}(t-t')F(t')\,dt' \Big\| _{\overline{X}^{s,\frac12}\cap \overline{Y}^{s,0}}\lesssim \| F\| _{\overline{X}^{s,-\frac12}\cap \overline{Y}^{s,-1}}, 
\end{equation*}
where the implicit constant is independent of $s$, $a(\cdot )$, and $p(\cdot )$.
\end{lem}

\bigskip
\noindent
\textbf{Proof}.
We basically follow the argument in \cite[Proposition~2.3]{MR}, excepting that we avoid using the claim
\begin{align*}
   \Big\| &\chi _{\mathbf{R}_+}(t) \psi(t) \int_0^t e^{ia(\xi )(t-t') - p(\xi )|t-t'|} f(t') \ dt' \Big\| _{H^{\frac{1}{2}}(\mathbf{R})}   \\
   &\lesssim \Big\| \psi(t) \int_0^t e^{ia(\xi )(t-t') - p(\xi )|t-t'|} f(t') \ dt' \Big\| _{H^{\frac{1}{2}}(\mathbf{R})},  
\end{align*}
which seems not clear to hold.
We will instead use the fact that
\begin{equation}\label{claim:bdd1+}
f\in H^{1}(\mathbf{R}),\quad f(0)=0\quad \Longrightarrow \quad \| \chi _{\mathbf{R}_+}f\| _{H^{1}}\leq \| f\| _{H^{1}}.
\end{equation}
A similar estimate to \eqref{claim:bdd1+} is known to hold if $H^1$ is replaced with $H^s$ for any $\frac{1}{2}<s<\frac{3}{2}$ (cf.~Jerison and Kenig \cite[Lemmas~3.7, 3.8]{JK}), though it fails for $s=1/2$, as we will see in Lemma~\ref{counterex} below.
We only need the above elementary case $s=1$.

We begin with observing the equality (for $F\in \mathscr{S}(\mathbf{R}\times \mathcal{M})$)
\begin{align*}
&\mathcal{F}_{t,x}\left[ \psi (t)\int _0^t\overline{U}(t-t')F(t')\,dt'\right] (\tau , \xi )\\
&\quad =\mathcal{F}_t\left[ \psi (t)\int _0^te^{ia(\xi )(t-t')-p(\xi )|t-t'|}\hat{F}(t',\xi )\,dt'\right] (\tau )\\
&\quad =\mathcal{F}_t\left[ \psi (t)\int _0^te^{-p(\xi )|t-t'|}\int _{\mathbf{R}}e^{i\lambda t'}\tilde{F}(\lambda +a(\xi ),\xi )\,d\lambda \,dt'\right] (\tau -a(\xi ))\\
&\quad =\mathcal{F}_t\Big[ \Psi _\xi \big[ \tilde{F}(\cdot +a(\xi ),\xi )\big] \Big] (\tau -a(\xi )),
\end{align*}
where, after the computation
\begin{align*}
\int _0^te^{i\lambda t'-p(\xi )|t-t'|}\,dt'&=\frac{e^{i\lambda t}-e^{-p(\xi )|t|}}{i\lambda +\mathrm{sgn}(t)p(\xi )}=\frac{-i\lambda +\mathrm{sgn}(t)p(\xi )}{-i\lambda +p(\xi )}\cdot \frac{e^{i\lambda t}-e^{-p(\xi )|t|}}{i\lambda +p(\xi )}
\end{align*}
for $t\neq 0$, we have introduced the notation
\begin{align*}
\Psi _\xi [f](t)&:=\psi (t)\int _{\mathbf{R}}\Big( \int _0^te^{i\lambda t'}e^{-p(\xi )|t-t'|} \,dt'\Big) f(\lambda )\,d\lambda \\
&\;=\Phi _\xi \Big[ \frac{-i\lambda}{-i\lambda +p(\xi )}f\Big] (t)+\mathrm{sgn}(t)\Phi _\xi \Big[ \frac{p(\xi )}{-i\lambda +p(\xi )}f\Big] (t),\\[5pt]
\Phi _\xi [g](t)&:=\psi (t)\int _{\mathbf{R}}\frac{e^{i\lambda t}-e^{-p(\xi )|t|}}{i\lambda +p(\xi )}g(\lambda )\,d\lambda .
\end{align*}
The claimed estimate is then reduced to the following:
\begin{equation}\label{est:inhom}
\begin{split}
&\left\| \langle i\tau +p(\xi )\rangle ^{\frac{1}{2}}\mathcal{F}_t\big[ \Psi _\xi [f]\big] (\tau )\right\| _{L^2_\tau}+\left\| \mathcal{F}_t\big[ \Psi _\xi [f]\big] (\tau )\right\| _{L^1_\tau}\\
&\qquad \lesssim \left\| \frac{f(\lambda )}{\langle i\lambda +p(\xi )\rangle ^{1/2}}\right\| _{L^2_\lambda}+\left\| \frac{f(\lambda )}{\langle i\lambda +p(\xi )\rangle}\right\| _{L^1_\lambda},
\end{split}
\end{equation}
uniformly in $\xi \in \hat{\mathcal{M}}$ and $p(\cdot )$.

We will prove several estimates on $\Phi _\xi$ in Lemma~\ref{lem:Phi} below.
Let us admit these estimates for a moment and proceed to the proof of \eqref{est:inhom}.
The contribution from the first term of $\Psi _\xi$ can be treated directly by \eqref{est:Phi1} and \eqref{est:Phi3}, because $\big|\frac{-i\lambda}{-i\lambda +p(\xi )}\big| \leq 1$.
To estimate the second term of $\Psi _\xi$ with $\mathrm{sgn}(t)$, we first see that 
\[ \text{L.H.S. of \eqref{est:inhom}}~\lesssim \langle p(\xi )\rangle ^{-\frac{1}{2}}\left\| \langle i\tau +p(\xi )\rangle \mathcal{F}_t\big[ \Psi _\xi [f]\big] (\tau )\right\| _{L^2_\tau},\]
which follows from the Cauchy-Schwarz inequality (see the proof of Lemma~\ref{L1}).
We then get the bound 
\begin{align*}
&\langle p(\xi )\rangle ^{-\frac{1}{2}}\left\| \langle i\tau +p(\xi )\rangle \mathcal{F}_t\Big[ \mathrm{sgn}(t)\Phi _\xi \Big[ \frac{p(\xi )}{-i\lambda +p(\xi )}f\Big] \Big] \right\| _{L^2_\tau}\\
&\leq \langle p(\xi )\rangle ^{-\frac{1}{2}}\left\| \mathrm{sgn}(t)\Phi _\xi \Big[ \frac{p(\xi )}{-i\lambda +p(\xi )}f\Big] \right\| _{H^1_t}\!\!\!\!+\langle p(\xi )\rangle ^{\frac{1}{2}}\left\| \mathrm{sgn}(t)\Phi _\xi \Big[ \frac{p(\xi )}{-i\lambda +p(\xi )}f\Big] \right\| _{L^2_t}.
\end{align*}
Now, we invoke \eqref{claim:bdd1+} and obtain the bound
\begin{align*}
&\langle p(\xi )\rangle ^{-\frac{1}{2}}\left\| \Phi _\xi \Big[ \frac{p(\xi )}{-i\lambda +p(\xi )}f\Big] \right\| _{H^1_t}+\langle p(\xi )\rangle ^{\frac{1}{2}}\left\| \Phi _\xi \Big[ \frac{p(\xi )}{-i\lambda +p(\xi )}f\Big] \right\| _{L^2_t}\\
&\quad \sim \langle p(\xi )\rangle ^{-\frac{1}{2}}\left\| \langle i\tau +p(\xi )\rangle \mathcal{F}_t\Big[ \Phi _\xi \Big[ \frac{p(\xi )}{-i\lambda +p(\xi )}f\Big] \Big] \right\| _{L^2_\tau},
\end{align*}
which is, by \eqref{est:Phi2}, evaluated by
\[ \langle p(\xi )\rangle ^{-\frac{1}{2}}\left\| \frac{p(\xi )}{-i\lambda +p(\xi )}f(\lambda )\right\| _{L^2_\lambda}\lesssim \left\| \frac{1}{\langle i\lambda +p(\xi )\rangle ^{1/2}}f(\lambda )\right\| _{L^2_\lambda},\]
as desired.
\hfill $\Box$


\begin{lem}[{cf.~\cite[Proposition~2.2]{MR}}]\label{lem:Phi}
Let $p:\hat{\mathcal{M}}\to [0,\infty )$ and $\xi \in \hat{\mathcal{M}}$.
Then, we have the following estimates, with the implicit constants independent of $p(\cdot )$ and $\xi$:
\begin{gather}
\left\| \langle i\tau \!+\!p(\xi )\rangle ^{\frac{1}{2}}\mathcal{F}_t\big[ \Phi _\xi [g]\big] (\tau )\right\| _{L^2_\tau}\lesssim \left\| \frac{g(\lambda )}{\langle i\lambda \!+\!p(\xi )\rangle ^{1/2}}\right\| _{L^2_\lambda}+\left\| \frac{g(\lambda )}{\langle i\lambda \!+\!p(\xi )\rangle}\right\| _{L^1_\lambda}, \label{est:Phi1} \\[5pt]
\left\| \langle i\tau +p(\xi )\rangle \mathcal{F}_t\big[ \Phi _\xi [g]\big] (\tau )\right\| _{L^2_\tau}\lesssim \| g(\lambda )\| _{L^2_\lambda}, \label{est:Phi2} \\[5pt]
\left\| \mathcal{F}_t\big[ \Phi _\xi [g]\big] (\tau )\right\| _{L^1_\tau}\lesssim \left\| \frac{g(\lambda )}{\langle i\lambda +p(\xi )\rangle}\right\| _{L^1_\lambda}. \label{est:Phi3}
\end{gather}
\end{lem}

\bigskip \noindent
\textbf{Proof}.
Although the proof is almost the same as that for \cite[Proposition~2.2]{MR}, we provide the detail for completeness.
We will focus on the case $p(\xi )\neq 0$, though the argument can be easily adapted to the case $p(\xi )=0$.

Divide $\Phi _\xi [g]$ as follows:
\begin{align*}
\Phi _\xi [g](t)&=\psi (t)\int e^{i\lambda t}\frac{\chi_{\{|\lambda |+p(\xi )\geq 1\}}g(\lambda )}{i\lambda +p(\xi )}\,d\lambda \\
&\quad -\psi (t)e^{-p(\xi )|t|}\int \frac{\chi_{\{|\lambda |+p(\xi )\geq 1\}}g(\lambda )}{i\lambda +p(\xi )}\,d\lambda \\
&\quad +\psi (t)\int \frac{e^{i\lambda t}-1}{i\lambda +p(\xi )}\chi_{\{|\lambda |+p(\xi )\leq 1\}}g(\lambda )\,d\lambda \\
&\quad -\psi (t)\big( e^{-p(\xi )|t|}-1\big) \int \frac{\chi_{\{|\lambda |+p(\xi )\leq 1\}}}{i\lambda +p(\xi )}g(\lambda )\,d\lambda \\
&=:\Phi ^1_\xi [g](t)-\Phi ^2_\xi [g](t)+\Phi ^3_\xi [g](t)-\Phi ^4_\xi [g](t).
\end{align*}

\medskip
{\bf Estimate on $\Phi ^1_\xi$}.
Since
\[ \left| \mathcal{F}_t\big[ \Phi ^1_\xi [g]\big] \right| = \Big| \hat{\psi} * \Big[ \frac{\chi_{\{|\tau |+p(\xi )\geq 1\}}g(\tau )}{i\tau +p(\xi )}\Big] \Big| \lesssim \big| \hat{\psi} \big| * \frac{|g(\tau )|}{\langle i\tau +p(\xi )\rangle} \]
and $\langle i\tau +p(\xi )\rangle \leq \langle \tau -\tau' \rangle +\langle i\tau' +p(\xi )\rangle$, we have
\begin{align*}
&\left\| \langle i\tau +p(\xi )\rangle ^{b}\mathcal{F}_t\big[ \Phi ^1_\xi [g]\big] \right\| _{L^2_\tau}\\
&\quad \leq \| \psi \| _{H^{b}_t}\left\| \frac{g(\tau )}{\langle i\tau +p(\xi )\rangle}\right\| _{L^1_\tau}+\| \hat{\psi}\| _{L^1_\tau}\left\| \frac{g(\tau )}{\langle i\tau +p(\xi )\rangle^{1-b}}\right\| _{L^2_\tau},\qquad 0\leq b\leq 1,\\
&\left\| \mathcal{F}_t\big[ \Phi ^1_\xi [g]\big] \right\| _{L^1_\tau}\leq \| \hat{\psi} \| _{L^1_\tau}\left\| \frac{g(\tau )}{\langle i\tau +p(\xi )\rangle}\right\| _{L^1_\tau}.
\end{align*}
These inequalities and the Cauchy-Schwarz inequality yield \eqref{est:Phi1}--\eqref{est:Phi3} in this case.

\medskip
{\bf Estimate on $\Phi ^2_\xi$}.
We have
\[ \left| \mathcal{F}_t\big[ \Phi ^2_\xi [g]\big] \right| = \Big| \int \frac{\chi_{\{|\lambda |+p(\xi )\geq 1\}}g(\lambda )}{i\lambda +p(\xi )}\,d\lambda \Big| \big| \hat{\psi} * \mathcal{F}_t\big( e^{-p(\xi )|t|}\big) \big| ,\]
We have already seen in the proof of Lemma~\ref{lem:linear1} that
\[ \left\| \langle i\tau +p(\xi )\rangle ^{\frac{1}{2}}\hat{\psi} * \mathcal{F}_t\big( e^{-p(\xi )|t|}\big) \right\| _{L^2_\tau}+\left\| \hat{\psi} * \mathcal{F}_t\big( e^{-p(\xi )|t|}\big) \right\| _{L^1_\tau}\lesssim 1,\]
which verifies \eqref{est:Phi1} and \eqref{est:Phi3}.
Moreover, a similar argument shows
\[ \left\| \langle i\tau +p(\xi )\rangle \hat{\psi} * \mathcal{F}_t\big( e^{-p(\xi )|t|}\big) \right\| _{L^2_\tau}\lesssim 1+p(\xi )^{\frac{1}{2}}.\]
We apply the Cauchy-Schwarz inequality to obtain
\[ \left\| \langle i\tau +p(\xi )\rangle \mathcal{F}_t\big[ \Phi ^2_\xi [g]\big]\right\| _{L^2_\tau}\lesssim \langle p(\xi )\rangle ^{\frac{1}{2}}\left\| \frac{g(\lambda )}{\langle i\lambda +p(\xi )\rangle}\right\| _{L^1_\lambda} \lesssim \| g\| _{L^2},\]
which shows \eqref{est:Phi2}.

\medskip
{\bf Estimate on $\Phi ^3_\xi$}.
By the Taylor expansion, we see
\[ \mathcal{F}_t\big[ \Phi ^3_\xi [g]\big] (\tau )=\sum _{n\geq 1}\frac{\mathcal{F}_t[t^n\psi ](\tau )}{n!}\int \frac{(i\lambda )^n\chi_{\{ |\lambda |+p(\xi )\leq 1\}}}{i\lambda +p(\xi )}g(\lambda )\,d\lambda .\]
On one hand, the restriction $|\lambda |+p(\xi )\leq 1$ implies (for $n\geq 1$) that
\[ \Big| \int \frac{(i\lambda )^n\chi_{\{ |\lambda |+p(\xi )\leq 1\}}}{i\lambda +p(\xi )}g(\lambda )\,d\lambda \Big| \lesssim \left\| \chi_{\{ |\lambda |+p(\xi )\leq 1\}}g(\lambda )\right\| _{L^1_\lambda}\lesssim \left\| \frac{g(\lambda)}{\langle i\lambda +p(\xi )\rangle}\right\| _{L^1_\lambda }.\]
On the other hand, the support property of $\psi$ implies (under the condition $p(\xi )\leq 1$ and for $0\leq b\leq 1$) that
\[ \left\| \langle i\tau +p(\xi )\rangle ^b\mathcal{F}_t[t^n\psi ]\right\| _{L^2_\tau}+\left\| \mathcal{F}_t[t^n\psi ]\right\| _{L^1_\tau}\lesssim \| t^n\psi \| _{H^1_t}\lesssim n2^n.\]
These estimates yield that
\begin{align*}
&\left\| \langle i\tau +p(\xi )\rangle ^b\mathcal{F}_t\big[ \Phi ^3_\xi [g]\big] \right\| _{L^2_\tau}+\left\| \mathcal{F}_t\big[ \Phi ^3_\xi [g]\big] \right\| _{L^1_\tau}\lesssim \Big( \sum _{n\geq 1}\frac{n2^n}{n!}\Big) \left\| \frac{g(\lambda)}{\langle i\lambda +p(\xi )\rangle}\right\| _{L^1_\lambda },
\end{align*}
verifying \eqref{est:Phi1} and \eqref{est:Phi3} in this case.
The estimate \eqref{est:Phi2} then follows from the Cauchy-Schwarz inequality.

\medskip
{\bf Estimate on $\Phi ^4_\xi$}.
This is very similar to the estimate on $\Phi ^3_\xi$ above.
In fact, since 
\begin{align*}
\big| \mathcal{F}_t\big[ \Phi ^4_\xi [g]\big] (\tau )\big| 
&=\bigg| \sum _{n\geq 1}\frac{\mathcal{F}_t[|t|^n\psi ](\tau )}{n!}\int \frac{(-p(\xi ))^n\chi_{\{ |\lambda |+p(\xi )\leq 1\}}}{i\lambda +p(\xi )}g(\lambda )\,d\lambda \bigg| \\
&\lesssim \sum _{n\geq 1}\frac{|\mathcal{F}_t[|t|^n\psi ](\tau )|}{n!}\left\| \frac{g(\lambda)}{\langle i\lambda +p(\xi )\rangle}\right\| _{L^1_\lambda }
\end{align*}
and
\[ \big\| |t|^n\psi \big\| _{H^1_t}=\| t^n\psi \| _{H^1_t}\lesssim n2^n,\]
the same argument as before can be applied.
%

This completes the proof of Lemma~\ref{lem:Phi}.
\hfill $\Box$

\bigskip
The following lemma is concerned with the trace of functions in $H^{1/2}(\mathbf{R})$.
\begin{lem} \label{counterex} 
There exists a sequence $\{ f_n \} \subset H^{1/2}(\mathbf{R})$ such that $f_n$ is continuous on $\mathbf{R}$, $f_n(0) = 0$,  $\sup_{n \geq 1} \| f_n \|_{H^{1/2}(\mathbf{R})} < \infty$ and
\begin{align*}
   \bigl \| \chi_{\mathbf{R}_+} f_n \bigr \|_{H^{1/2}(\mathbf{R})} \ \longrightarrow \ \infty \quad (n \to \infty).
\end{align*}
In particular, there is no positive constant $C$ such that for any $f \in H^{1/2}(\mathbf{R}) \cap C(\mathbf{R})$ with $f(0) = 0$,
\begin{equation}\label{est:sharpcutoff}
   \bigl \| \chi_{\mathbf{R}_+} f \bigr \|_{H^{1/2}(\mathbf{R})}\leq C \| f \|_{H^{1/2}(\mathbf{R})}.
\end{equation}
\end{lem}

\begin{rem}
The estimate \eqref{est:sharpcutoff} does not hold for general functions $f$ in $H^{1/2}(\mathbf{R})$, and hence multiplication by $\chi_{\mathbf{R}_+}$ is not bounded in $H^{1/2}(\mathbf{R})$.
Indeed, for $f\in H^{1/2}(\mathbf{R})$ it is known that $\chi_{\mathbf{R}_+}f\in H^{1/2}(\mathbf{R})$ if and only if $\int _0^\infty t^{-1}|f(t)|^2\,dt <\infty$ (see Strichartz \cite[Theorem~3.1 on page 1054]{Str67} and Lions and Magenes \cite[Theorem~11.7 on page 66]{LM}), while the last integral clearly diverges if we take, for instance, $f\in C^\infty_0(\mathbf{R})$ such that $f=1$ on $[-1,1]$.
In the above lemma, however, we aim to find a counterexample from among continuous functions vanishing at the origin.
\end{rem}

\noindent
\textbf{Proof of Lemma~\ref{counterex}}.
We first recall the following characterization of the $H^{1/2}(\mathbf{R})$ norm (see Strichartz \cite[Statement 2.3 on page 1035]{Str67}).
\begin{gather*}
\| f\| _{H^{1/2}}\ \sim \ \| f\| _{L^2}+\| S_{1/2}f\|_{L^2},\\
(S_{1/2}f)(t)=\bigg[ \int _0^\infty \Big( \int _{-1}^1|f(t+ru)-f(t)|\,du\Big) ^2\,\frac{dr}{r^2}\bigg] ^{1/2}.
\end{gather*}
From this, we can easily deduce that
\begin{equation}\label{est:counterex}
\int _0^\infty \frac{|f(t)|^2}{t}\,dt \ \lesssim \ \big\| \chi _{\mathbf{R}_+}f\big\| _{H^{1/2}}^2
\end{equation}
for any function $f$ on $\mathbf{R}$ such that the right side is finite (see \cite[the proof of Theorem~3.1 on page 1054]{Str67}).

We define functions $f_n$, $n=1,2,\dots $ by
\[ \hat{f}_n(\tau )\ =\ \frac{\psi _n(\tau )}{(e+|\tau |)\log (e+|\tau |)},\qquad \tau \in \mathbf{R},\]
where $\psi _n\in C^\infty _0(\mathbf{R})$ is chosen so that
\begin{gather*}
\psi _n(-\tau )=\psi _n(\tau ),\qquad \psi _n(\tau )=1~~(|\tau |\leq n), \\
-1\leq \psi _n\leq 1,\qquad |\psi_n'(\tau )|\lesssim (1+|\tau |)^{-1}
\end{gather*}
and 
\[ \big( f_n(0)=\big) \ \int _{\mathbf{R}}\frac{\psi _n(\tau )}{(e+|\tau |)\log (e+|\tau |)}\,d\tau \ =\ 0.\]
Such a function $\psi_n$ can be constructed as follows.
First, we fix $\phi \in C_0^\infty (\mathbf{R})$ even, $0\leq \phi \leq 1$, $\phi =1$ on $[-1,1]$, and monotone on $[0,\infty)$.
Then, for each $n$, there exists $N>n$ such that 
\[ \int _{\mathbf{R}}\frac{\phi (\tau /N)}{(e+|\tau |)\log (e+|\tau |)}\,d\tau \ =\ 2 \int _{\mathbf{R}}\frac{\phi (\tau /n)}{(e+|\tau |)\log (e+|\tau |)}\,d\tau ,\]
because $\int _{\mathbf{R}}\frac{d\tau}{(e+|\tau |)\log (e+|\tau |)}=\infty$.
Now, we define $\psi _n(\tau ):=2\phi (\tau /n)-\phi (\tau /N)$, which satisfies the desired properties.

Clearly, $\{ f_n\}$ is bounded in $H^{1/2}(\mathbf{R})$.
Since $\hat{f}_n\in L^1(\mathbf{R})$, we have $f_n\in C(\mathbf{R})$ and
\[ f_n(t) = \int _{\mathbf{R}}\frac{\psi _n(\tau )e^{it\tau }}{(e+|\tau |)\log (e+|\tau |)}\,d\tau = 2\int _0^{\infty}\frac{\psi _n(\tau )\cos (t\tau )}{(e+\tau )\log (e+\tau )}\,d\tau ,\qquad t\in \mathbf{R}.\]
%
%
For $n^{-1}\leq t\leq 1$, we have
\[ \int _0^{t^{-1}}\frac{\psi _n(\tau )\cos (t\tau )}{(e+\tau )\log (e+\tau )}\,d\tau \ \gtrsim \ \log \log (e+t^{-1}),\]
while an integration by parts shows, for $0<t\leq 1$, that
\begin{align*}
&\Big| \int _{t^{-1}}^\infty \frac{\psi _n(\tau )\cos (t\tau )}{(e+\tau )\log (e+\tau )}\,d\tau \Big| \\
&\quad \lesssim \frac{1}{t}\Big| \Big[ \frac{\psi _n(\tau )\sin (t\tau )}{(e+\tau )\log (e+\tau )}\Big] _{t^{-1}}^\infty \Big| +\frac{1}{t}\int _{t^{-1}}^\infty \Big| \Big( \frac{\psi _n(\tau )}{(e+\tau )\log (e+\tau )}\Big) '\Big| \,d\tau \\
&\quad \lesssim \frac{1}{\log (e+t^{-1})},
\end{align*}
since the last inequality is verified by observing that
\[ \Big| \Big( \frac{\psi _n(\tau )}{(e+\tau )\log (e+\tau )}\Big) '\Big| \ \lesssim \ \frac{1}{(e+\tau )^2\log (e+\tau )}\qquad (\tau >0).\]
Consequently, for some $0<\delta \ll 1$ we have
\[ f_n(t)\ \gtrsim \ \log \log (e+t^{-1}) \ \gtrsim \ 1 \qquad (n^{-1}\leq t\leq \delta ),\]
and thus
\[ \int _0^\infty \frac{|f_n(t)|^2}{t}\,dt \ \gtrsim \ \int _{n^{-1}}^\delta \frac{1}{t}\,dt \ \sim \ \log n\qquad (n\gg 1).\] 
The claim then follows from \eqref{est:counterex}.
\hfill $\Box$


\bigskip
{\bf Acknowledgements}.
The first author N.K is partially supported by JSPS KAKENHI Grant-in-Aid for Young Researchers (B) (16K17626) and Grant-in-Aid for Scientific Research (C) (20K03678).
The second author Y.T is partially supported by JSPS KAKENHI Grant-in-Aid for Scientific Research (B) (17H02853).



\bigskip
\begin{thebibliography}{15}
\addcontentsline{toc}{section}{References}

\bibitem{BP} H.~Bahouri and G.~Perelman, 
Global well-posedness for the derivative nonlinear Schr\"odinger equation,
preprint, December 2020. 
\url{https://arxiv.org/abs/2012.01923}

\bibitem{Bour} J.~Bourgain, 
Fourier transform restriction phenomena for certain lattice subsets and applications to nonlinear evolution equations. II. The KdV-equation,
\textit{Geom. Funct. Anal.} {\bf 3} (1993), 209--262.

\bibitem{DP} K.~B.~Dysthe and H.~L.~P\'ecseli, 
Non-linear Langmuir wave modulation in collisionless plasma,
\textit{Plasma Physics} {\bf 19} (1977), 931--943.

\bibitem{G} Z.~Guo, 
Local well-posedness and a priori bounds for the modified Benjamin-Ono equation,
\textit{Adv. Differential Equations} {\bf 16} (2011), 1087--1137.

\bibitem{HKV} B.~Harrop-Griffiths, R.~Killip and M.~Vi\c{s}an, 
Large-data equicontinuity for the derivative NLS,
preprint, June 2021. 
\url{https://arxiv.org/abs/2106.13333}


\bibitem{H} S.~Herr, 
On the Cauchy problem for the derivative nonlinear Schr\"odinger equation with periodic boundary condition,
\textit{Int. Math. Res. Not.} {\bf 2006}, Art. ID 96763, 33 pages.

\bibitem{IK} A.~D.~Ionescu and C.~E.~Kenig, 
Global well-posedness of the Benjamin-Ono equation in low-regularity spaces,
\textit{J. Amer. Math. Soc.} {\bf 20} (2007), 753--798.

\bibitem{IKT} A.~D.~Ionescu, C.~E.~Kenig and D.~Tataru, 
Global well-posedness of the KP-I initial-value problem in the energy space,
\textit{Invent. Math.} {\bf 173} (2008), 265--304.

\bibitem{JK} D.~Jerison and C.~E.~Kenig, 
The inhomogeneous Dirichlet problem in Lipschitz domains, 
\textit{J. Funct. Anal.} {\bf 130} (1995), 161--219.

\bibitem{KPV} C.~E.~Kenig, G.~Ponce and L.~Vega, 
Small solutions to nonlinear Schr\"odinger equations,
\textit{Ann. Inst. H. Poincar\'e Anal. Non Lin\'eaire} {\bf 10} (1993), 255--288.

\bibitem{K09} N.~Kishimoto, 
Low-regularity bilinear estimates for a quadratic nonlinear Schr\"odinger equation,
\textit{J. Differential Equations} {\bf 247} (2009), 1397--1439.

\bibitem{KT18} N.~Kishimoto and Y.~Tsutsumi, 
Ill-posedness of the third order NLS equation with Raman scattering term,
\textit{Math. Res. Lett.} {\bf 25} (2018), 1447--1484.

\bibitem{KT20b} N.~Kishimoto and Y.~Tsutsumi, 
Ill-posedness of the Third Order NLS with Raman Scattering Term in Gevrey Spaces, 
in \textit{Mathematics of Wave Phenomena, Trends in Mathematics}, Birkh\"auser, Basel, 2020, 219--233.

\bibitem{KT20} N.~Kishimoto and Y.~Tsutsumi, 
A priori bounds for the kinetic DNLS, 
in \textit{2019-20 MATRIX Annals}, Springer, Cham, 2021, 797--803.

\bibitem{LM} J.~L.~Lions and E.~Magenes, 
Non-Homogeneous Boundary Value Problems and Applications, Vol.~1, 
Springer-Verlag, Berlin, Heidelberg, New York, 1972.

\bibitem{LiuZ} T.-P.~Liu and Y.~Zeng,
Shock Waves in Conservation Laws with Physical Viscosity,
\textit{Mem. Amer. Math. Soc.}, {\bf 234} (2015), no.~1105.

\bibitem{MT18} T.~Miyaji and Y.~Tsutsumi, 
Local well-posedness of the NLS equation with third order dispersion in negative Sobolev spaces, 
\textit{Differential Integral Equations} {\bf 31} (2018), 111--132.

\bibitem{MW} E.~Mj\o lhus and J.~Wyller, 
Nonlinear Alfv\'en waves in a finite-beta plasma,
\textit{J. Plasma Physics} {\bf 40} (1988), 299--318.

\bibitem{MPV19} L.~Molinet, D.~Pilod and S.~Vento, 
On unconditional well-posedness for the periodic modified Korteweg-de Vries equation, 
\textit{J. Math. Soc. Japan} {\bf 71} (2019), 147--201.

\bibitem{MR} L.~Molinet and F.~Ribaud, 
On the low regularity of the Korteweg-de Vries-Burgers equation,
\textit{Int. Math. Res. Not.} {\bf 2002}, no.~37, 1979--2005.

\bibitem{MR04a} L.~Molinet and F.~Ribaud, 
Well-posedness results for the generalized Benjamin-Ono equation with small initial data,
\textit{J. Math. Pures Appl. (9)} {\bf 83} (2004), 277--311.

\bibitem{MV1} L.~Molinet and S.~Vento, 
Sharp ill-posedness and well-posedness results for the KdV-Burgers equation: the real line case,
\textit{Ann. Sc. Norm. Super. Pisa Cl. Sci. (5)} {\bf 10} (2011), 531--560.

\bibitem{MV2} L.~Molinet and S.~Vento, 
Sharp ill-posedness and well-posedness results for the KdV-Burgers equation: the periodic case,
\textit{Trans. Amer. Math. Soc.} {\bf 365} (2013), 123–-141.

\bibitem{NTT10} K.~Nakanishi, H.~Takaoka and Y.~Tsutsumi, 
Local well-posedness in low regularity of the mKdV equation with periodic boundary condition,
\textit{Discrete Contin. Dyn. Syst.} {\bf 28} (2010), 1635--1654.

\bibitem{MP} R.~Peres de Moura and A.~Pastor, 
On the Cauchy problem for the nonlocal derivative nonlinear Schr\"odinger equation,
\textit{Commun. Math. Sci.} {\bf 9} (2011), 63--80.

\bibitem{Str67} R.~S.~Strichartz, 
Multipliers on fractional Sobolev spaces, 
\textit{J. Math. Mech.} {\bf 16} (1967), 1031--1060.

\bibitem{Tak1} H.~Takaoka, 
Well-posedness for the one-dimensional nonlinear Schr\"odinger equation with the derivative nonlinearity, 
\textit{Adv. Differential Equations} {\bf 4} (1999), 561--580.

\bibitem{TT04} H.~Takaoka and Y.~Tsutsumi, 
Well-posedness of the Cauchy problem for the modified KdV equation with periodic boundary condition, 
\textit{Int. Math. Res. Not.} {\bf 2004}, no.~56, 3009--3040.

\bibitem{Tsug} K.~Tsugawa, 
Parabolic smoothing effect and local well-posedness of fifth order semilinear dispersive equations on the torus, 
preprint, July 2017. 
\url{https://arxiv.org/abs/1707.09550}
\end{thebibliography}
\end{document}